\documentclass[11pt]{amsart}
\usepackage{geometry}

\usepackage[latin1]{inputenc}
\usepackage{amsfonts}
\usepackage[toc,page,title,titletoc,header]{appendix}
\usepackage{graphicx,psfrag,epsfig,multirow,caption,subcaption}
\usepackage{amssymb,amsmath,amscd,amsthm,amssymb,verbatim,setspace}
\numberwithin{equation}{section}
\usepackage{mathrsfs,wrapfig}
\usepackage{indentfirst}
\usepackage{extarrows}

\usepackage[colorlinks=true]{hyperref}

\usepackage{marginnote}
\usepackage{color}
\numberwithin{equation}{section}
\newtheorem{Thm}{Theorem}[section]
\newtheorem{Lm}[Thm]{Lemma}
\newtheorem{SLm}[Thm]{Sublemma}
\newtheorem{Prop}[Thm]{Proposition}

\newtheorem{Def}[Thm]{Definition}
\newtheorem{Not}[Thm]{Notation}
\newtheorem{Rk}[Thm]{Remark}
\newtheorem{Conjecture}{Conjecture}

\newcommand{\al}{\alpha}
\newcommand{\Id}{\mathrm{Id}}

\newcommand{\cS}{\mathcal{S}}
\newcommand{\cR}{\mathcal{R}}

\newcommand{\cM}{\mathcal{M}}

\newcommand{\cD}{\mathcal{D}}

\newcommand{\bL}{\mathbb{L}}
\newcommand{\bG}{\mathbb{G}}
\newcommand{\bP}{\mathbb{P}}

\newcommand{\R}{\mathbb{R}}

\newcommand{\cC}{\mathcal{C}}

\newcommand{\bu}{\mathbf{u}}
\newcommand{\bl}{\mathbf{l}}

\newcommand{\dt}{\delta}

\newcommand{\eps}{\epsilon}

\def\sx{\mathsf{x}}

\def\sC{\mathsf{C}}

\def\bx{\mathbf{x}}
\def\bp{\mathbf{p}}

\title[Super-hyperbolic orbits and noncollision singularities]{Super-hyperbolic orbits and noncollision singularities  in a four-body problem}
\begin{document}
\author{Guan Huang and Jinxin Xue}
\email{huangguan@tsinghua.edu.cn}
\address{Yau Mathematical Sciences Center, Jinchunyuan West Building 304, Tsinghua University, Beijing, China, 100084}
\email{jxue@tsinghua.edu.cn}
\address{Department of Mathematics, Jingzhai 310, Tsinghua University, Beijing, China, 100084}

\maketitle
\begin{abstract}
In this paper, we prove the existence of super-hyperbolic orbits in four-body problem, which solves a conjecture of  Marchal-Saari. We also prove the existence of noncollision singularities in the same model, which solves a conjecture of Anosov. Moreover, the two type of solutions coexist for certain mass ratios. We also make a conjecture on the classification of final motions of $N$-body problem, in which the superhyperbolic type of orbits is one of the building blocks.\end{abstract}

\tableofcontents
\addtocontents{toc}{\setcounter{tocdepth}{1}} 
\section{Introduction}

In this paper, we construct a model of four-body problem in which there exist both noncollision singularities and another kind of exotic solutions called superhyperbolic solutions, where for the latter the size of the system grows super-linearly with respect to time, i.e. $\frac{\max_i|Q_i(t)|}{t}\to\infty$ as $t\to\infty$. The study of noncollision singularities, i.e. finite time blowup solutions in $N$-body problem was motivated by the Painlev\'e conjecture. Two models of four-body problems with noncollision singularities were constructed in \cite{X1,GHX} to which we refer readers for backgrounds on noncollision singularities. In this paper, we  focus mainly on superhyperbolic solutions. This type of solutions was conjectured in \cite{MS}, which is defined globally in time, thus is not a singularity. However, since a superhyperbolic solution has velocity grows to infinity by definition, repeated accelerations are needed to construct such a solution. Therefore, we can consider superhyperbolic solutions as a slow-down version of noncollision singularities. 

Our result in this paper solves the existence problem of superhyperbolic solutions. It is worthy of putting the result into a broader framework. Indeed, superhyperbolic solutions have an important position in the classification of final motions (asymptotic behaviors as $t\to\infty$)  in $N$-body problem for $N\geq 4$.
For three-body problem, it is well-known that there is a
\index{Chazy classification} classification of the final motions given by Chazy (c.f. \cite{AKN}) a hundred years ago, in which all logical possible combinations of the following building blocks
\begin{equation}\label{EqBuildingBlocks3}\mathcal H,\mathcal P, \mathcal{B},\mathcal{OS}\end{equation} is realized by an orbit, where $\mathcal{H}$ means \emph{hyperbolic}, i.e. $|Q_k(t)|\to\infty,|\dot Q_k(t)|\to c_k>0$ for $k=1,2,3$, $\mathcal P$ means \emph{parabolic}, i.e. $|Q_k(t)|\to \infty,\ |\dot Q_k(t)|\to 0$ for $k=1,2,3$, $\mathcal B$ means \emph{bounded}, i.e. $\sup_t\max_k |Q_k(t)|<\infty$, $\mathcal{OS}$ means \emph{oscillatory}\index{oscillatory orbit}, i.e. $\limsup|Q(t)|=\infty,\ \liminf|Q(t)|<\infty$, and $\mathcal{HP}$ means \emph{hyperbolic-parabolic} i.e. some body is hyperbolic and some other body is parabolic asymptotically, similarly for $\mathcal{HB}$ and $\mathcal{PB}$ (c.f. \cite[Chapter 2.3.4 ]{AKN} for more details). 

The classification of final motions in $N$-body problem for $N\geq 4$ is far less studied. We make the following conjecture, which claims that in addition to the above building blocks \eqref{EqBuildingBlocks3} for three-body problem, we only need  to  add one new block $\mathcal{SH}$, representing super-hyperbolic orbits. 
\begin{Conjecture}\label{ConSH}
	For 	$N$-body problem for $N\geq 4$,  all logical possible combinations of $$\mathcal H,\mathcal P,\mathcal B,\mathcal{OS},\mathcal{SH}$$ is realized by an orbit and there is no other possibilities. 
\end{Conjecture}
The conjecture is strongly supported by the following theorem of Marshal-Saari \cite{MS}. 

\begin{Thm}\label{ThmMS}
	For $N(>3)$-body problem as $t\to\infty$, there is a dichotomy: either it is a superhyperbolic orbit; or it decouples into several subsystems moving apart linearly and each subsystem grows at most like $O(t^{2/3})$. 
\end{Thm}
This theorem leads to a possible approach to Conjecture \ref{ConSH} by analysing the growth rates. Indeed, letting $I(t)=\sum m_i|Q_i(t)|^2$ be the momentum of inertia of the system, Theorem  \ref{ThmMS} implies that $\sqrt I(t)/t=o(1)$ can only be at most $\sqrt I(t)=O(t^{2/3})$ as $t\to\infty$. Thus,  Theorem  \ref{ThmMS} reduces Conjecture \ref{ConSH} to the following one.

\begin{Conjecture}
	For 	$N$-body problem for $N\geq 4$,  suppose a orbit satisfies $\sqrt I(t)=O(t^{2/3})$ as $t\to\infty$, then it is a combination of $\mathcal P,\mathcal B,\mathcal{OS}.$ Conversely, any logical possible  combination of $\mathcal P,\mathcal B,\mathcal{OS}$ is realized by an orbit satisfying $\sqrt I(t)=O(t^{2/3})$ as $t\to\infty$. 
\end{Conjecture}
%Note that this is true for three-body problem by Chazy's classification. 

%It seems the following conjecture is also true. 
%\begin{Conjecture}
	%For 	$N$-body problem for $N\geq 4$,  suppose we have $\sqrt I(t)=o(t^{2/3})$, the solution is bounded $\mathcal B$. 
%\end{Conjecture}

%From the above discussion, we see that superhyperbolic orbit is the only possible new final motions in $N>3$ case than $N=3$ case, if the conjectures hold. In this paper, we address the very existence problem of superhyperbolic orbits. 

As the complexity of possible final motions grows exponentially as $N$ gets large, we believe that a classification as asserted in Conjecture \ref{ConSH} gives a clean picture worthy of pursuing. %The situation can be understood by analogy as follows. For a positive series $\sum a_n$ that decays faster than exponential $a_n\leq e^{-n}$, it is easy to know that it is convergent. However, for a se

The following theorem is analogous to the classical von Zeipel theorem for noncollision singularities (see Appendix of \cite{X2}). 
\begin{Thm}
	Let $(Q_i(t),\dot Q_i(t))_{i=1}^N$ be a superhyperbolic solution of the $N$-body problem. Then we have as $t\to\infty$
	$$\limsup_{t\to\infty}\max_{i,j}|Q_i(t)-Q_j(t)|=\infty,\quad \liminf_{t\to\infty}\min_{i,j}|Q_i(t)-Q_j(t)|=0.$$
\end{Thm}
From this theorem, we can view a superhyperbolic orbit as a slow version of a noncollision singularity. The proof if in fact very simple. First, the $\limsup$ equation follows trivially from the definition of a superhyperbolic orbit. To prove the liminf one, we assume the liminf is $a>0$. Then in the Hamiltonian, the potential energy is bounded, thus by the energy conservation, the kinetic energy, hence the velocities are also bounded, so the orbit cannot be superhyperbolic, a contradiction. 
\subsection{The configuration}
As we have discussed in the beginning of the paper, superhyperbolic orbits can be considered as a slow-down version of noncollision singularities. However, constructing such an orbit is by no means easy, since it requires more delicate control of the return times between two close encounters than noncollision singularities. As a result, it seems that superhyperbolic orbits are in some sense more delicate and rarer than noncollision singularities. In fact,  we are not able to find superhyperbolic orbits in the models \cite{X1,GHX}. Instead, we use here a new model and it is remarkable that in this model, superhyperbolic orbits and noncollision singularities coexist.

\begin{figure}[ht]
\begin{center}
\includegraphics[width=0.4\textwidth]{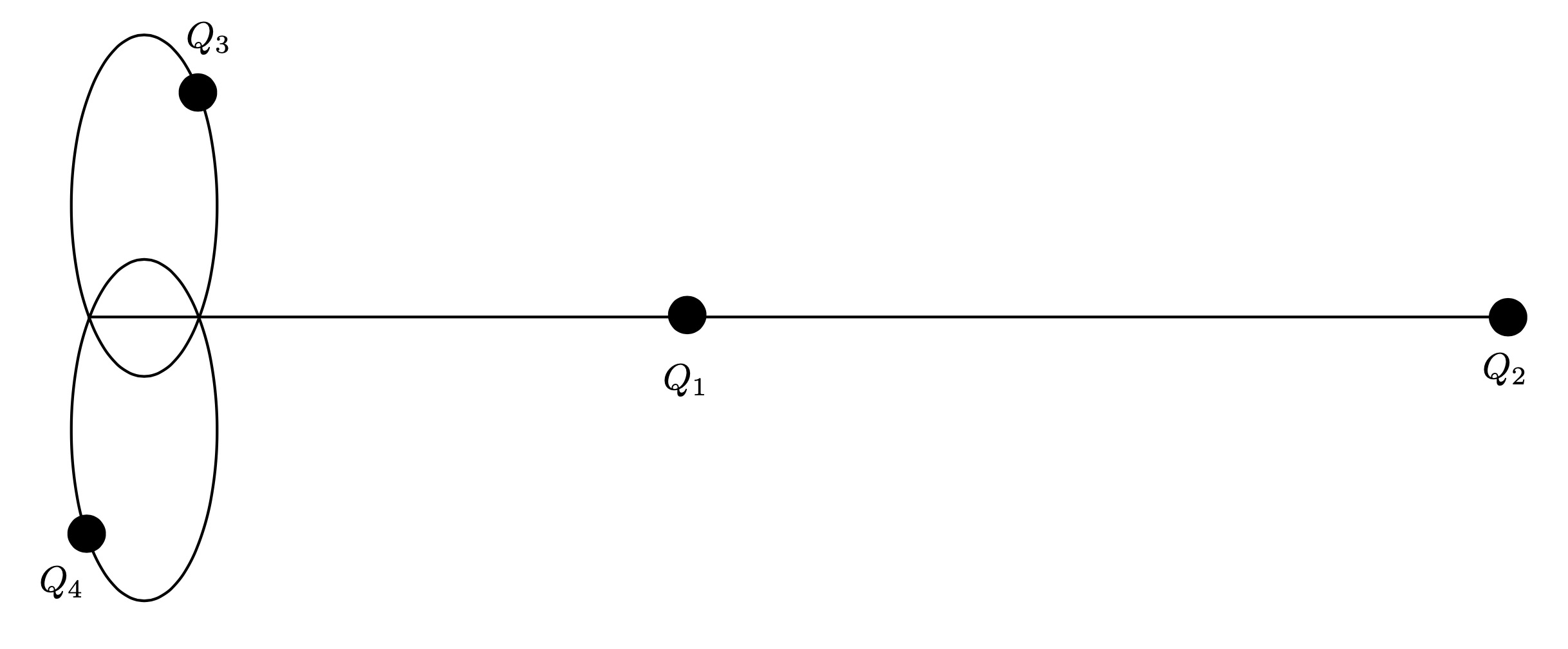}
\end{center}
\caption{The configuration of the four-body problem}
\label{config1}
\end{figure}

The configuration is as follows (see Figure \ref{config1}): there is  a particle  $Q_2$ with mass $m_2$ moving to the right nearly along the $x$-axis and a binary $Q_3$-$Q_4$ with masses $m_3=m_4$ performing nearly Kepler elliptic motion when it is far away from other particles. The mass center of the binary moves to the opposite direction of $Q_2$ nearly along the $x$-axis. There is another particle $Q_1$ with mass $m_1$ traveling back and forth between the binary and $Q_2$ nearly along the $x$-axis.  When $Q_1$ approaches the binary, a triple is formed and the triple will pass close to a triple collision after which all the three particles are significantly accelerated and the particle $Q_1$ is ejected to the right direction.  It will quickly catch up with $Q_2$. Through two-body interaction, a large amount of the kinetic  energy of $Q_1$ would  transfer to $Q_2$. Then $Q_1$ would change the  direction of movement and slowly  get close to the binary.  We need to carefully choose the  mass ratios and the rate of acceleration to make sure this configuration would repeat itself infinitely many times and exists globally in time, so that the super-hyperbolic orbits could exhibit in this model. 

%The main mechanism of acceleration is from the near triple collision dynamics. There is a case of triple collision analyzed by Devaney called the isosceles three-body problem \cite{D}. This is a configuration with two equal masses having the same distances to the third mass. Assume the relative position of the two equal masses is parallel to the $y$-axis and the third mass lies on the $x$-axis, then it was classified by Devaney that there are two types of final motions after approaching near triple collision: either the two equal masses move apart along the $y$-axis while the third body oscillates on the $x$-axis, or the third mass moves along the $x$-axis towards infinity with arbitrarily large velocity and the binary moves opposite to the third body, experiencing repeated binary collisions.  In our model, we ask our near-triple-collision orbits to shadow the latter type.

\subsection{The statements of the results.}
We first consider the isosceles four-body problem by requiring $|Q_i-Q_j|=|Q_{7-i}-Q_j|$, $i=3,4,j=1,2$, where $Q_1$, $Q_2$ and the mass center of $Q_3$-$Q_4$ have always horizontal velocities. In this case, double collisions of the binary $Q_3$-$Q_4$ and of $Q_1$ and $Q_2$ are not avoidable.  The collisions of these types can be  regularized as elastic collisions. For this system, we have the following result. 
\begin{Thm}\label{ThmI4BP}
There is an open set $\mathcal{A}\subset\mathbb{R}$ that contains 1 such that for each $m_1\in\mathcal{A}$ there exists $\lambda_0>1$ such that for each $\lambda\geqslant\lambda_0$,  there exists an open set $\mathcal{B}_{\lambda}\subset\mathbb{R}$  such that if $m_2\in\mathcal{B}_{\lambda}$, then for the isosceles four-body problem with masses $(m_1,m_2)$ and $m_3=m_4=1$, we have 
\begin{enumerate}
	\item there exists a nonempty set $\Sigma_{sh}$ of initial conditions such that for all $x\in \Sigma_{sh}$, the orbit starting from $x$ exists globally in time upon the regularization of double collisions and 
	$$\frac{1}{t}\max\{|Q_i(t)|, i=1,2,3,4\}\to\infty, \mathrm{\ as\ } t\to\infty.$$
	\item there exists a non-empty set $\Sigma_{ncs}$ of initial conditions such that for all $x\in\Sigma_{ncs}$, there exists $t_x<+\infty$ such that as $t\to t_x$, $\max\{|Q_i(t)|, i=1,2,3,4\}\to\infty$. %Moreover, there exists $x\in\Sigma$ such that $\lim_{n\to\infty}\frac{|\dot Q_1(t_{n+2})|}{|\dot Q_1(t_n)|}=\infty$, where $t_n<t_x$ is a sequence of times when $Q_1$ visits a fixed neighborhood of the origin for the $n$-th time. 
\end{enumerate}
 %Moreover, for any sequence of positive numbers $\{a_n\to\infty\}$, there exists $x\in\Sigma$ and a sequence $0<t_1<t_1'<t_2<t_2'<\cdots<t_x$ such that for the  orbit starting from $x$, we have
%$$\frac{|\dot{Q}_3(t_n')+\dot{Q}_4(t_n')|}{|\dot{Q}_3(t_n)+\dot{Q}_4(t_n)|}>a_n,$$
%where $t_n<t_n'$ are two consecutive moments when the mass center of the pair $Q_3$-$Q_4$ passes through the origin.
%Moreover, for any sequence of numbers $\{a_n\}\to \infty$, we have $\frac{1}{2}|\dot Q_3(t_n)+\dot Q_4(t_n)|\geq a_n$ where $t_n(<t_x)$ is the $n$-th time when the mass center of $Q_3$ and $Q_4$ passes through the origin.
\end{Thm}

Perturbing slightly away from the above solutions, we shall allow the binary to gain certain nonzero but very tiny angular of momentum. In this way, we avoid double collision. We obtain the following result on the existence of  super-hyperbolic orbits and non-collision singularities  in the   model.
\begin{Thm}\label{ThmF4BP}
There is an open set $\mathcal{A}\subset \R$ that contains 1  such that for each $m_1\in\mathcal{A}$ there exists $\lambda_0>1$ such that for each $\lambda\geqslant\lambda_0$, there exists an open set $\mathcal{B}_{\lambda}\subset\mathbb{R}$  such that if $m_2\in\mathcal{B}_{\lambda}$, then for the full four-body problem with masses $(m_1,m_2)$ and $m_3=m_4=1$, we have 
\begin{enumerate}
	\item   there exists a nonempty set $\Sigma_{sh}$ of initial conditions such that for all $x\in \Sigma_{sh}$, the orbit starting from $x$ exists globally in time without the occurrence of any collision between the particles  and 
	$$\frac{1}{t}\max\{|Q_i(t)|, i=1,2,3,4\}\to\infty, \mathrm{\ as}\ t\to\infty.$$
	\item  there exists a non-empty set $\Sigma_{ncs}$ of initial conditions such that for all $x\in\Sigma_{ncs}$, there exists $t_x<+\infty$ such that as $t\to t_x$, $\max\{|Q_i(t)|, i=1,2,3,4\}\to\infty$ and no collision happens before $t_x$. %Moreover, there exists $x\in\Sigma$ such that $\lim_{n\to\infty}\frac{|\dot Q_1(t_{n+2})|}{|\dot Q_1(t_n)|}=\infty$, where $t_n<t_x$ is a sequence of times when $Q_1$ visits a fixed neighborhood of the origin for the $n$-th time. 
	\end{enumerate}
 \end{Thm}
 In fact, for only the existence of non-collision singularities a much wider range of mass ratios is allowed, see Remark \ref{mass-ratio-remark} and the discussion in Subsection \ref{1.3} below.

In the last two theorems, by interpolating between noncollision singularities and superhyperbolic orbits, we can indeed get superhyperbolic orbits with $\frac{1}{t}\max\{|Q_i(t)|, i=1,2,3,4\}$ growing to infinity with arbitrarily high speed. 

We discussed in \cite{GHX} that the main theorem therein solves an analogue of a conjecture of Anosov claiming that a noncollision singularity can be found in a neighborhood of the singularities of \cite{MM}. The model that we use in the present paper can be considered to be close to the model of \cite{MM}. Thus, it is reasonable to consider statement (2) of  Theorem \ref{ThmF4BP} as a solution to the conjecture of Anosov.

In \cite{GHX} we prove the existence of noncollision singularities in a similar model. There we consider a binary moving back and forth between two particles moving off to infinity in opposite directions. The mechanism of acceleration is similar to here, i.e. Devaney's isosceles three-body problem. However, there are some important differences between the model of \cite{GHX} and the model here: 
\begin{enumerate}
	\item The model in the present paper does not admit noncollision singularies with four equal masses in contrast to \cite{GHX}. 
	\item During each return, the model of \cite{GHX} accelerates twice. The rate of acceleration is so large that we do not know how to slow it down thus  it is not clear to us whether it is possible to construct superhyperbolic orbit in the model of \cite{GHX}.  
	\item A main difficulty in the proofs of the theorems in the present paper is to exclude double collision between $Q_1$ and $Q_2$, which is not present in the model of \cite{GHX}. 
\end{enumerate}

\subsection{Backgrounds and implications}\label{1.3}

%Superhyperbolic orbits are related intimately to noncollision singularities. We may think superhyperbolic orbits as slow version of noncollision singularities. However, 

%As we have discussed above, constructing a superhyperbolic orbit requires more delicate control than that constructing a noncollision singularity. 

As we have discussed above, whereas near triple collision accelerates very intensively such that we can find noncollision singularities that accelerates arbitrarily fast, to slow down the orbit becomes a subtle problem. %As a result, for the model studied in \cite{GHX}, we are not able to construct superhyperbolic orbit therein. 

The heuristic idea goes as follows. In Figure \ref{config1}, the velocity of $Q_1$ after the triple collision can be accelerated arbitrarily and we fix such a ratio $\lambda$ throughout. When $Q_1$ comes close to $Q_2$, by momentum conservation, there is an energy transfer between the two bodies. Suppose the mass of $Q_2$ is much smaller than that of $Q_1$, then both will move towards the direction opposite to the binary so $Q_1$ will never return to interact with the binary for the next time. On the other hand, if the mass of $Q_1$ is much smaller than that of $Q_2$, then $Q_1$ will always return and a noncollision singularity may occur, which is the case for the model in \cite{X1}. Thus, we expect that if we arrange the mass ratio of $Q_2$ and $Q_1$ carefully, we can control the return speed of $Q_1$ so small and its return time so long that to complete infinitely many returns, it takes infinite amount of time, which gives a superhyperbolic orbit. 

This consideration leads to the following problem. It seems there is a critical mass ratio $\nu^*$ of $Q_2$ and $Q_1$ such that if $m_2/m_1>\nu^*$, then there is no superhyperbolic orbit. 

{\bf Problem: } \emph{In this model, what is the critical ratio and what happens to the critical mass ratio, i.e. does superhyperbolic orbit exist or not?}

Some readers may complain the conjectured picture of possible final motions for $N(>3)$-body problem in Conjecture \ref{ConSH} remains too complicated when the number of bodies grows. To pursue a simpler picture, it is reasonable to consider generic initial conditions or almost every initial data. Recall that it is conjectured that noncollision singularies has zero Lebesgue measure and first categroy. It is also natural to conjecture the same is true for superhyperbolic orbits. Finally, we recapitulate the following conjecture from \cite{X2} and refer readers to find more discussions therein. 
\begin{Conjecture} For generic initial condition, the solution of the $N$-body problem is globally defined on $\R$, and as  $t\to \pm\infty$, 
	each body in the system approaches either a linear motion with constant velocity or a Kepler elliptic motion around a center moving linearly with constant velocity. 
\end{Conjecture}
The paper is organized as follows. In Section \ref{SCoord}, we introduce the coordinates: Jacobi coordinates and blowup coordinates, and the preliminaries on the isosceles three-body problem. In Section \ref{SMM}, we study the isosceles four-body probem and proof of Theorem \ref{ThmI4BP}. In Section \ref{SJacobi}, we give the Hamiltonians of the four-body problem and proves the key return time estimate. In Section \ref{SMainProof}, we prove the main theorem \ref{ThmF4BP}. Finally, we have four appendices. In Appendix \ref{App-delaunay}, we give some formulas for Delaunay coordinates.   We give In Appendix \ref{L-Derivatives} the proof of Proposition \ref{PropDG1} on the derivative calculation for the $L,\ell$ parts and In Appendix \ref{SG0g0} Proposition \ref{PropDG2} on the derivative estimates for the $G$, $g$ parts . This calculation is very similar to that in \cite{GHX}, thus we mainly focus on the differences.    At last, in Appendix \ref{app-numeric}, we present a numeric verification of the non-degenerate conditions.  

\section*{Acknowledgment}
The authors are supported by grant NSFC (Significant project No.11790273) in China. J.X. is in addition supported by the Xiaomi endowed professorship of Tsinghua University.
\section{Preliminaries}\label{SCoord}
In this section we introduce several systems of coordinates for the four-body problem that we will use later.
Without loss of generality, we assume $m_3=m_4=1$ throughout the paper and we always assume the binary $Q_3$-$Q_4$  is on the left, $Q_2$ is on the right and $Q_1$ travels between the binary and $Q_2$. 
\subsection{The Jacobi-Cartesian coordinates}\label{subsection-jacobi}
The first step is to remove the translation invariance. We introduce $q_i=Q_i-Q_4,\ i=3,1,2$. So we get that the following symplectic form is preserved by assuming $\sum_{i=1}^4 P_i=0$
$$\sum_{i=1}^4 dP_i\wedge dQ_i=\sum_{i=1}^3 dP_i\wedge dq_i.$$
We next assume that the particle $Q_1$ is closer to the binary than $Q_2$ and introduce is the following, which we call the {\it left Jacobi coordinates},
\begin{equation}\label{EqJacobi}
\begin{cases}
\bx_0&=q_3\\
\bx_1&=q_1-\frac{q_3}{2}\\
\bx_2&=q_2-\frac{m_1 q_1+ q_3 }{m_1+2}
\end{cases},\quad
\begin{cases}
\bp_0&=P_3+\frac{ P_1}{2}+\frac{ P_2}{2}\\
\bp_1&=P_1+\frac{ m_1 P_2}{m_1+2}\\
\bp_2&=P_2
\end{cases}.
\end{equation}
It can be easily checked that we have the following reduced symplectic form preserved by the coordinates change
$$\sum_{i=3,1,2} dP_i\wedge dq_i=\sum_{i=0,1,2}d\bp_i\wedge d\bx_i.$$
Clearly, $(\bx_0,\bp_0)$ describes the relative motion between the binary $Q_3$-$Q_4$, $(\bx_1,\bp_1)$ does that between $Q_1$ and the mass center of the binary, and $(\bx_2, \bp_2)$ does that between $Q_2$ and the mass center of the triple $Q_1$-$Q_3$-$Q_4$.

\subsection{The isosceles three-body problem (I3BP)}
\begin{figure}[ht]
\begin{center}
\includegraphics[width=0.8\textwidth]{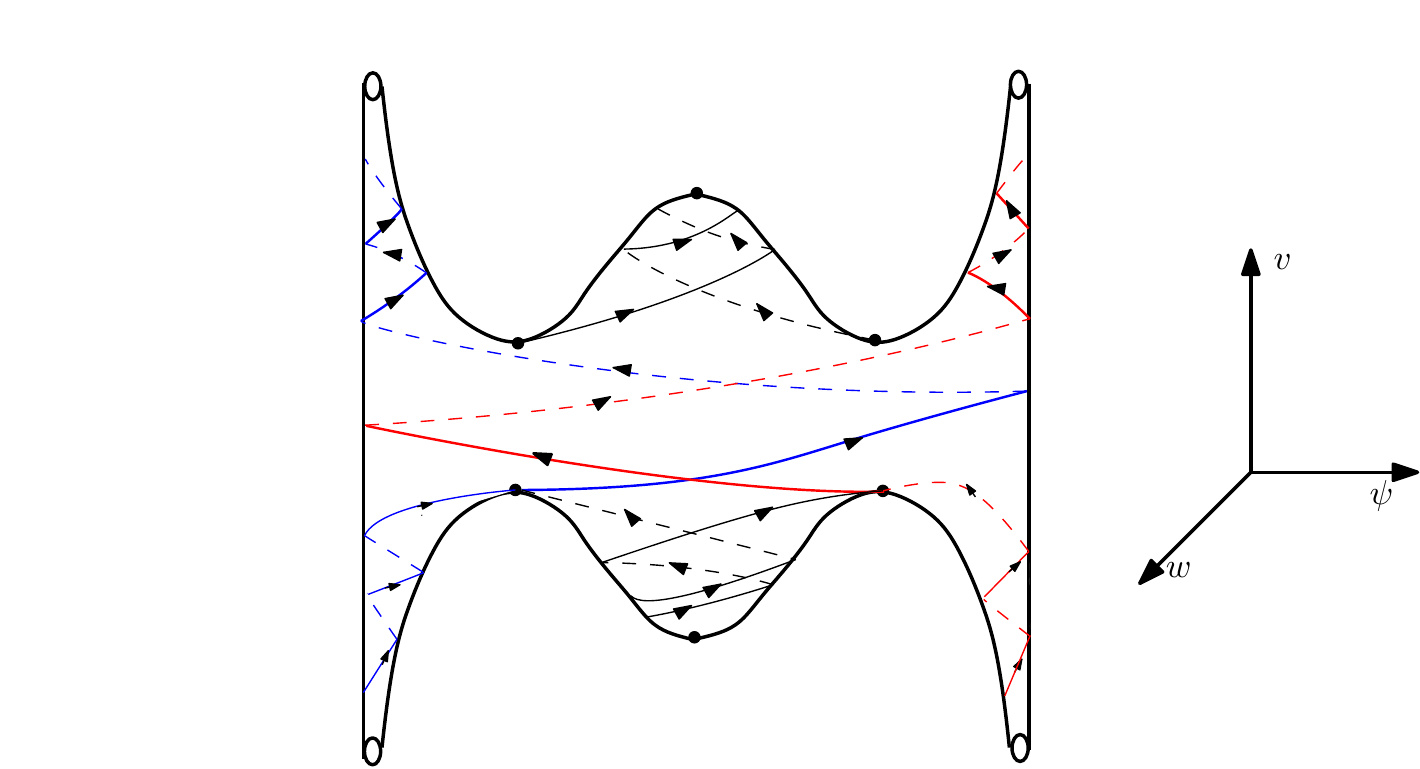}
\label{fig-collision}
\end{center}
\caption{The collision manifold $\cM_0$}
\label{collision-manifold}
\end{figure}

 In this section, we analyze the case when the triple $Q_1$-$Q_3$-$Q_4$ is close to triple collision and $Q_2$ is far apart. As a first approximation, we ignore $Q_2$ and focus only on the $Q_1$-$Q_3$-$Q_4$ three-body problem near triple collision. The problem is called the isosceles three-body problem first studied by Devaney \cite{D}. In Jacobi coordinates  we have
\[\bx_0=(0,x_0)^t,\quad \bx_1=(x_1,0)^t,\quad \bp_0=(0,p_0)^t,\quad \bp_1=(p_1,0)^t,\]
We denote $\bx=(M_0^{1/2}\bx_0,M_1^{1/2}\bx_1)\in \R^4$ and $\bp=(M_0^{-1/2}\bp_0,M_1^{-1/2}\bp_1)\in \R^4$ and by $|\cdot|$ the Euclidean norm, and by \begin{equation}\label{reduced-mass}M_0=\frac12,\  M_1=\frac{2m_1}{m_1+2}\end{equation} the reduced masses.
We next introduce the blowup coordinates
%\begin{equation}\label{EqBlowup}
\[\begin{aligned}
\begin{cases}
r&=|\bx|= \sqrt {M_0|\bx_0|^2 + M_1|\bx_1|^2},\\
v&=r^{-1/2}(\bx\cdot\bp) ,\\
\psi&=\arctan\frac{\sqrt{M_1}x_1}{\sqrt{M_0}x_0},\\
w& =r^{-1/2}(\sqrt{M_0/M_1} x_0 p_1-\sqrt{M_1/M_0} x_1p_0).
\end{cases}
\end{aligned}\]
%\end{equation} 
The physical meanings of the variables are as follows. The variable $r$ measures the size of the I3BP, $v$ is the projection of the rescaled momentum $r^{1/2}\bp$ to the radial componnet $\bx$, $w$ is the scalar part of $\mathbf w$  where $\mathbf w=r^{1/2}\bp-v r^{-1}\bx$ is the projection of $r^{1/2}\bp$ to the tangential of the sphere $|r^{-1}\bx|=1$ and $\psi$ measures the relative size of the positions $x_0$ and $x_1$.
All variables except $r$ are rescaling invariant and the triple collision corresponds to $r=0$. The coordinate change is accompanied by a time reparametrization $dt= r^{3/2}d\tau$. Equations of motion can be derived in terms of these coordinates
\begin{equation}\label{triple-iso}
\begin{cases}
r'&=rv\\
v'&=w^2+\frac{1}{2}v^2+\bar V(\psi),\\
\psi'&=w,\\
w'&=-\frac{1}{2} vw-\partial_\psi \bar V(\psi),
\end{cases}
\end{equation}
where \begin{equation}\label{iso-potential}\bar V(\psi)=-\frac{1}{\sqrt{2}\cos\psi}-\frac{4m_1}{\sqrt{2\cos^2\psi+\frac{2(2+m_1)}{m_1}\sin^2\psi}}.\end{equation}
and the energy relation becomes
\begin{equation}\label{EqEnergy}
rE=\frac{1}{2}v^2+\frac{1}{2}w^2+\bar V(\psi),
\end{equation}
where $E$ is the constant value of the Hamiltonian.
 In particular $r$-equation is $\frac{dr}{d\tau}=rv$, so we see that $\{r=0\}$ is an invariant submanifold. All the energy level sets are the same in the limit $r\to 0$, denoted by $\mathcal M_0$, has two dimensions. We call $\mathcal{M}_0$ the collision manifold and the dynamics on $\mathcal{M}_0$ is   illustrated in Figure \ref{collision-manifold}.

\begin{Thm}[\cite{D}] \label{collision-m1}
\begin{enumerate}
\item The collision manifold $\mathcal M_0$ is a topological 2-sphere with four punctures $($called arms$)$ with each puncture corresponding to $\psi=\pm\pi/2$. The manifold $\cM_0$ is symmetric with respect to $\psi\mapsto -\psi$ or $w\mapsto -w$.
\item The variable  $v$ is a Lyapunov function in the sense that $\frac{d}{d\tau}v\geq 0$ along the flow on $\mathcal M_0$ and $\frac{d}{d\tau}v=0$ iff at the fixed points.
\item There are six fixed points on $\mathcal M_0$: the two with $(v,\psi,w)=(\pm v_\dagger,0,0)$ correspond to Euler $($collinear$)$ central configurations and the four with $(v,\psi,w)=(\pm v_*,\pm\psi_*,0)$ correspond to Lagrange $($equilateral$)$ central configurations, where $v_\dagger,v_*$ are numbers depending only on masses and $\psi_*=\arctan\left(\sqrt\frac{3m_1}{2+m_1}\right)$. $($In the following, since we consider only the lower Lagrange point, we fix the convention $v_*<0$$)$. 
\item The Euler fixed point with $v>0$ is a sink and with $v<0$ is a source.  For some masses, the eigenvalues are real and for others complex.  The Lagrange fixed points are saddles.
%\item For a Lagrange fixed point with $v<0$, one of the stable manifolds is homo-clinic to the Euler fixed point with $v<0$; one stable manifold  comes along one lower arm and one unstable manifold  climbs to one upper arm.
\end{enumerate}
\end{Thm}
We are only interested in the Lagrange fixed points in the lower half space $v<0$. The main observation is that, we can arrange the orbit to stay close to the saddle for as long time as we wish by selecting  the initial condition sufficiently close to its stable manifold. Since in the lower half space we have $v_*<0$,  from the equation $\frac{d}{d\tau}r=rv$ we get that $r$ will decrease to as small as we wish before leaving a neighborhood of triple collision.
%We denote by $W^s_0$ the stable manifold of left lower Lagrange fixed point on $\mathcal M_0$ and by $W^s_{iso}$ its stable manifold in the isosceles three-body problem. We see that $W^s_0$ has one-dimension and is a submanifold of $W^s_{iso}$.
We shall consider the situation that the $Q_1$ comes from the right to have a near triple collision with the binary $Q_3$-$Q_4$. We hope that after the near triple collision the binary moves to the left and $Q_1$ moves to the right. In this case the relative position $x_1$ from the mass center of the binary to $Q_1$ has a positive sign before and after the near triple collision, therefore the variable $\psi$ should also carry a positive sign before and after the near triple collision.  Thus we need the right lower Lagrange fixed point to have a stable manifold coming from the right lower arm and an unstable manifold escaping to the right upper arm.  The existence of the stable manifold coming from the right lower arm follows directly from the fact that $v$ is a Lyapunov function. However, the existence of an unstable manifold escaping to the right upper arm is nontrivial, which depends on the mass ratio $m_1/m_3=m_1$. This was studied by \cite{SM}. Here we only cite the relevant statements.
\begin{Thm}[\cite{SM}]\label{collision-m2} Assume $m_1<55/4$ in which case the Euler fixed points are respectively sink and source with complex eigenvalues. Then there exist $b_1\simeq 0.379$ and $b_2\simeq 2.662$ such that the following holds for the right lower Lagrange fixed point:
\begin{enumerate}
\item when $b_1<m_1<b_2$, the two  unstable manifolds escape from the two upper arms respectively;
\item when $m_1>b_2$, one of the unstable manifold escapes from the right upper arm and the other dies at the upper Euler fixed point.
\end{enumerate}
\end{Thm}

%We are only interested in the Lagrange fixed points in the lower half space $v<0$. The main observation is that if the initial condition is sufficiently close to the stable manifold, then the orbit will stay close to the saddle for as long time as we wish.
%Since in the lower half space we have $v_*<0$, so from the equation $r'=rv$ we get that $r$ will decrease to as small as we wish.

\begin{Not}\label{NotWs}
\begin{enumerate}
\item Denote by $O$ the right lower Lagrange fixed point on the collision manifold of $Q_1$-$Q_3$-$Q_4$.
\item Denote  by $W^s_{\therefore}(O)$ the stable manifold of $O$, in the isosceles three-body problem $Q_1$-$Q_3$-$Q_4$. We see that $W^s_{\therefore}(O)$ has have dimension two.
\item Denote by $W^u_0(O)$ the branch of the unstable manifold of $O$ on the collision manifold $\mathcal{M}_0$ that escapes from the right upper arm. 
\end{enumerate}
\end{Not}
Note the fact the an orbit starting sufficiently close to the stable manifold $W^s_\therefore(O)$ of the Lagrange fixed point $O$ will follow the unstable manifold of that point very far. There are two branches of the unstable manifold on $\mathcal M_0$ and there are two sides of the two dimensional manifold $W^s_\therefore(O)$ in the three dimensional energy level set. In order for the exiting orbit to follow the correct branch of the unstable manifold $W^u_0(O)$ defined above, we have to choose initial condition on the correct side of $W^s_\therefore(O)$.  In the rest of the paper we use the phrase ``{\it correct side}" of the stable manifold $W^s_{\therefore}(O)$ for this meaning. 
\section{The isosceles four-body problem and proof of Theorem~\ref{ThmI4BP}}\label{SMM}
In this section we consider the isosceles four-body problem (I4BP) and prove Theorem~\ref{ThmI4BP}. %We assume that $Q_1$ and $Q_2$ move on the $x$-axis  and $Q_3$-$Q_4$ stays on collision-ejection orbit perpendicular to the $x$-axis. %Thus we have $G_0=g_0-\pi/2=G_1=g_1=G_2=g_2=0$ and the system has three degrees of freedom. 
We first introduce the Poincar\'e sections $\mathcal S_\pm$ to separate the near triple collision regime and the perturbed Kepler motions regime.
\begin{Def}\label{DefSection}
	Let us fix $\epsilon>0$ to be a sufficiently small number  whose meaning will be reserved throughout the paper. We introduce a Poincar\'e section
	$\mathcal S_-:=\{r=\eps^{-1}\}$ before the near triple collision and $\cS_+=\{v=\eps^{-1/2}\}$  after the near triple collision. \end{Def}

\subsection{The local and global maps}
Using the sections $\cS_\pm$, we introduce the return maps $\bL$ and $\bG$ called local map and global map respectively. 
\begin{Def}[Local and global maps]\label{DefLocalGlobal}
We define the \emph{local map} $\mathbb L$  to be the Poincar\'e map going from $\cS_- $ to $\cS_+$, and \emph{global map} $\mathbb G$  to be the Poincar\'e map going from $\cS_+ $ to $\cS_-$.% By default, we use the blowup coordinates for the local map and Delaunay coordinates for the global map, if not otherwise mentioned.
\end{Def}
For the global map piece of the orbit, we treat $(\bx_0,\bp_0)$ as a degenerate Kepler elliptic motion and introduce Delaunay coordinates $(L_0,\ell_0)$ for it, and treat $(\bx_1,\bp_1)$ and $(\bx_2,\bp_2)$ as degenerate Kepler hyperbolic motions and introduce Delaunay coordinates $(L_1,\ell_1)$ and $(L_2,\ell_2)$ respectively for them. The isosceles assumption guarantees that $G_0=g_0-\frac{\pi}{2}=G_1=g_1=G_2=g_2=0$. 
For the local map piece of orbit, we convert $(\bx_0,\bp_0,\bx_1,\bp_1)$ into blowup coordinates $(r,v,\psi,w)$ and $(\bx_2,\bp_2)$ into Delaunay coordinates $(L_2,\ell_2)$. 

\begin{Not} We use the super-script $i$ and $f$ to stand for the corresponding variables on the \emph{initial} and \emph{final} sections respectively.
\end{Not}

The following lemma says that the local map is well approximated by the dynamics of the isosceles three-body problem $Q_1$-$Q_3$-$Q_4$.

\begin{Prop}\label{PropLocalMM} For the I4BP,   
there exist $C_1>1,D>0, \chi_0:=\chi_0(m_1)\gg 1$ and $\delta_0:=\delta_0(m_1)>0$ such that for any $\chi>\chi_0$, $\delta\in[\chi^{-D},\delta_0]$ and $\mu\in(0,1)$ the following holds: Let  $\zeta:\ [0,1]\to \cS_-$ be a smooth curve on the section $\cS_-$ satisfying
\begin{enumerate}
\item
 $\zeta$  is on the correct side of $W^s_{\therefore}(O)$;
 \item  $\zeta(0)\in W^s_{\therefore}(O)$ and $\mathrm{dist}(\zeta(1),W^s_{\therefore}(O))=\epsilon^2$;
\item on the curve $\zeta$,  $|\bx_2|\geq \chi$ and  $|\bp_2|=O(1)$.
\end{enumerate}
Then \begin{enumerate}
\item  there exists a subsegment $\bar \zeta$ of $\zeta$, such that  for any point $\mathbf{x}$ on the image of  the curve $\bar\zeta$ under the local map, its $r$-coordinate satisfies
$r^f\in[(1-\mu)\delta,(1+\mu)\delta ]$.
\item the oscillation of $L_2$ is estimated as $|L_2^f-L_2^i|\leq \frac{C_ 1}{ \chi^3}$, and
\item the travel time for points on $\bar\zeta$ between the two sections satisfies $C_1^{-1}\log\dt^{-1}<|\tau^f-\tau^i|<C_1\log\dt^{-1}$, and in  the original time scale $t$, we have $|t^f-t^i|\leqslant C_1$.
\end{enumerate}
\end{Prop}
We refer the reader to \cite[Section 6.1]{GHX} for a proof of this statement. Here we only outline the main idea. For the I3BP, if a curve $\zeta$ is transverse to the stable manifold $W^s_{\therefore}(O)$, the closer of a point on $\zeta$ to $W^s_{\therefore}(O)$, the longer its orbit stays  in a neighborhood of $O$. From the $r'=rv$ equation, we get that the longer the orbit stays close to $O$ where $v=v_*<0$, the more $r$ decreases. So it is possible to select a subsegment on $\zeta$ such that $r^f$ lies in the given window $[(1-\mu)\dt,(1+\mu)\dt]$. Since $Q_2$ is far, its perturbation to the I3BP can be estimated to be small. 

Immediately after the local map, on the sections $\cS_+$, we have that $r^f$ is small and the velocities are large. So we introduce the renormalization map to zoom in the spatial variables and to slow down the velocities as follows. 
\begin{Def}[Renormalization map]\label{renormalize-map}  Let  $\lambda:=\lambda(\delta)=\frac{\epsilon^{-1}}{2\delta}$, where $\delta\leqslant\delta_0$ is as in Proposition \ref{PropLocalMM}.  The renormalization map $\cR_\lambda$ on the section $\cS_+$ is as follows,
$$(r, v,\psi,w, L_2,\ell_2)\mapsto (\lambda r, v,\psi,w, \sqrt\lambda L_2,\ell_2).$$
The effect of the renormalization on the Delaunay coordinates and the Cartesian coordinates are as follows respectively
\begin{equation}
\begin{aligned}
&(L_0,\ell_0,L_1,\ell_1, L_2,\ell_2)\mapsto (\sqrt\lambda L_0,\ell_0,\sqrt\lambda L_1,\ell_1,\sqrt\lambda L_2,\ell_2);\\
&(x_0,p_0;x_1,p_1;x_2,p_2)\mapsto (\lambda x_0,\lambda^{-1/2} p_0;\lambda x_1,\lambda^{-1/2} p_1;\lambda x_2,\lambda^{-1/2}  p_2).
\end{aligned}
\end{equation}
Moreover, we also make the time change $t\mapsto t\lambda^{3/2}$ and the change of Hamiltonian $H\mapsto H/ \lambda$ when applying the renormalization. 
\end{Def}

The main effect of the renormalization is to rescale the semimajor of the binary's elliptic motion to order 1 size. Indeed, since $r$ is of order $\dt$ after the local map, the renormalization stretches $r$ to order $1/\eps$. Since on $\cS_+$, we have $v=\eps^{-1/2}$. Using the definition of $v$, we see that $\bp$ is of order 1, thus both subsystems $(\bx_0,\bp_0)$ and $(\bx_1,\bp_1)$ have order 1 energies. 
Applying  the renormalization to the outcome of local map in Proposition \ref{PropLocalMM}, we have the following.

\begin{Prop}\label{PropRenorm} 
After applying the renormalization on the section $\cS_+$, for the image $\mathcal{R}_{\lambda(\delta)}\mathbb{L}\bar\zeta$ where $\bar\zeta$ is the resulting segment in Proposition \ref{PropLocalMM},  we have
\begin{equation}\label{EqRenorm}\{- E_0,E_1\}\subset (1-\mu,1+\mu),\qquad r_0,\;r_1\eps,\; L_2/\sqrt{\lambda(\delta)}=O_{\dt\to 0}(1). 
\end{equation}
Here $E_0:= -\frac{M_0}{2L_0^2}$, $E_1:=\frac{2m^2_1M_1}{L^2_1}$ are the energy for $(\bx_i,\bp_i)$, $i=0,1$ as defined in \eqref{energy-formula1} and $M_i$, $i=0,1$ are the reduced masses as in \eqref{reduced-mass}.
\end{Prop}
Readers can find a proof of this statement in \cite[Section 6.2]{GHX}. 

After the local map and the renormalization, $Q_1$ would move to the right with speed of order 1, and $Q_2$ moves with speed of order $\lambda^{-1/2}$. When $Q_1$ catches up with $Q_2$, we want that  after the two-body interaction between $Q_1$ and $Q_2$, a large part of the kinetic energy of $Q_1$ would transfer to $Q_2$, such that $Q_2$ continues moving to the right with speed of order 1 and $Q_1$ is bounced back to the left with  much smaller velocity: the relative speed between $Q_1$ and the mass center of the binary $Q_3$-$Q_4$ is of order $\lambda^{-1/2}$. This could be achieved by carefully choosing the mass ratio of $Q_1$-$Q_2$ and the initial speed of $Q_2$ as follows. 
\begin{Def} \label{mass-admissible}Let  $m_1\in(b_1,2)$, $0<\mu\ll 1$ and $\lambda\gg1$. We say that $m_2$ is $(m_1,\mu,\lambda)$-\emph{admissible} if  we have
\begin{equation}\label{admissible-m2}
\sqrt\mu\leq \frac{2m_2-m_1(m_1+m_2+2)}{(m_1+2)(m_1+m_2)}\lambda^{1/2}- \frac{8m_1m_2(m_1+m_2+2)}{(m_1+2)^2(m_1+m_2)^2}\leq 10\sqrt{\mu}.\end{equation}
\end{Def}
In the limit $\mu\to 0$ and $\lambda\to\infty$, we have $m_2\to \frac{m_1^2+2m_1}{2-m_1}$. Thus solution exists for $m_1\in(b_1,2)$ and $\mu,\lambda^{-1}$ small enough and the admissible $m_2$'s form an interval of length $O(\sqrt{\mu \lambda^{-1}})$.  Note that for $m_1\geqslant2$ the above condition will never be satisfied. This is natural, since if $m_1\geqslant2$, then when ejected from near triple collision, the absolute velocity of $Q_1$ would be smaller than that of  the mass center of the binary. In such scenario, $Q_1$ would never catch up with the binary after being reflected by $Q_2$. 

If the mass ratio is chosen as in the last definition, then we have the following estimate of the distance and travel time for $Q_1$ to reach the next $\cS_-$ section. 
\begin{Prop} \label{choiceofm2} 
Let   $m_1\in (b_1,2)$ $0<\mu\ll1$, and $\lambda_*,m_2$ be such that $m_2$ is  $(m_1,\mu,\lambda_*)$-admissible $($respectively, in addition $\lambda\geq 2\lambda_*)$. We choose $\eps$ in the definition of the sections $\cS_\pm$ to be smaller than $\mu$.  Suppose  on the section $\cS_+$ before  renormalization $\cR_{\lambda(\delta)}$ we have 
	\begin{equation}\label{p2-initial1}|\bp_2|\in\frac{2m_2\sqrt{2M_1}}{m_1+m_2}(1-2\mu,1+2\mu),
	\end{equation} and \eqref{EqRenorm} are satisfied after  renormalization $\cR_{\lambda_*}$ $($respectively $\cR_{\lambda})$. Then 
\begin{enumerate}
	\item the particle $Q_1$ ejected from near triple collision as in Proposition \ref{PropRenorm} would come to a two-body interaction with $Q_2$ after which $Q_1$ would return to the section $\cS_-$. 
	\item After renormalization $\cR_{\lambda_*}$ $($respectively $\cR_{\lambda})$, the timespan of these orbits from section $\cS_+$ to $\cS_-$ is $O_{\chi^{-1}\to0}(\mu^{-1/2}\lambda_*^{1/2}\chi)$ $($respectively $O_{\chi^{-1}\to0}(\lambda_*^{1/2}\chi))$ and when reaching $\cS_-$, \eqref{p2-initial1} is satisfied and the distance between $Q_1$ and $Q_2$ is $O_{\chi^{-1}\to0}(\mu^{-1/2}\lambda_*^{1/2}\chi)$ $($respectively $O_{\chi^{-1}\to0}(\lambda_*^{1/2}\chi))$, where $\chi$ is the initial distance between $Q_1$ and $Q_2$ on the section $\cS_+$ after  renormalization $\cR_{\lambda_*}$ $($respectively $\cR_{\lambda})$. 
\end{enumerate}

\end{Prop}
We will prove this proposition in Section \ref{proofpm2}. %For the global map $\mathbb{G}$: $\cS_+\to\cS_-$ we always refer to the orbits as in the above proposition. 

\subsection{The key derivative estimate}

On the sections $\cS_{\pm}$, we perform the standard energetic reduction to eliminate $(L_1,\ell_1)$ from Delaunay coordinates by fixing the total energy and treating the variable $\ell_1$ as the new time.  So we use  the four variables $(L_0,\ell_0,L_2,\ell_2)$ to parametrize the sections $\cS_{\pm}$. 
We have the following estimates of the derivatives of the local and global maps. In both cases, we use Delaunay coordinates. Denote by $\cM$ the phase space of the isosceles four-body problem (I4BP), which has six dimensions.

\begin{Prop}[Proposition 3.11 of \cite{GHX}] \label{PropDL1}
Let $\gamma:\ [0,T]\to \cM$ be an orbit of the I4BP with $\gamma(0)\in \cS_-$ and $\gamma(T)\in \cS_+$ as in Proposition \ref{PropLocalMM}. 
Then for $\dt$ sufficiently small, we have $T=O(\log\delta^{-1})$ and the following  estimate of the derivatives of the local map along $\gamma$
\begin{equation}
d\mathcal Rd\mathbb L=e^{-Tv_*} \mathbf u_1\otimes \mathbf l_1+O(1), %O\left[\begin{array}{cccc}\eps^2\chi&\eps^4\chi&\beta\eps^2\chi&\frac{\eps^2}{\beta^2}\\
%\chi&\eps^2\chi&\beta\chi&\frac{1}{\beta^2}\\
%1&\eps^2&\beta&\frac{1}{\beta^2\chi}\\
%\chi&\eps^2\chi&\beta\chi&\frac{1}{\beta^2}\end{array}\right].
\end{equation}
where using the variables $(L_0,\ell_0,L_2,\ell_2),$ we have $$\bu_1=(1,o(1),o(1),o(1)),\ \bl_1=(o(1),1,o(1),O(1))$$ as $\dt\to0$ and $\chi\to\infty$.
\end{Prop}
We refer readers to \cite{GHX} for a proof. The proposition implies that the only dominant term in $d\mathcal Rd\mathbb L$ is given by $\frac{\partial L^f_0}{\partial \ell_0^i}$. Since we can parametrize the curve $\zeta$ in Proposition \ref{PropLocalMM} by $\ell_0$ as a substitute of $\psi$, and $L_0^2$ has the meaning of the semimajor of the elliptic motion, the proposition means that by changing the initial point on $\zeta$ slightly, we can get a significant change of the semimajor of the binary on the section $\cS_+$, which follows from a similar analysis to Proposition \ref{PropLocalMM} of the hyperbolic dynamics near the Lagrange fixed point. 

We next estimate the derivative of the global map. 
\begin{Prop}\label{PropDG1}
There exists $\chi_0\gg1$ such that the following holds:  Let $\gamma:\ [0,T]\to \cM$ be an orbit of the I4BP as in Proposition \ref{choiceofm2} with $\gamma(0)\in \cS_+$ and $\gamma(T)\in \cS_-$ and  the initial condition $\gamma(0)$ satisfying \eqref{EqRenorm}  and in addition
$|\bx_2|= \chi\geq \chi_0.$
Then we have the following derivative estimate of the derivatives of the global map along $\gamma$ \begin{equation}\label{DGG}
d\mathbb G=f_{\mu,\lambda}\cdot \chi^2\bar\bu_1\otimes \bar\bl_1+O_{\chi^{-1}\to0}((\eps^2+\lambda^{-1/2})f_{\mu,\lambda}\cdot\chi), %O\left[\begin{array}{cccc}\eps^2\chi&\eps^4\chi&\beta\eps^2\chi&\frac{\eps^2}{\beta^2}\\
%\chi&\eps^2\chi&\beta\chi&\frac{1}{\beta^2}\\
%1&\eps^2&\beta&\frac{1}{\beta^2\chi}\\
%\chi&\eps^2\chi&\beta\chi&\frac{1}{\beta^2}\end{array}\right].
\end{equation}
where $f_{\mu,\lambda}=\mu^{-1/2}\lambda_*^{1/2}$ if $\lambda=\lambda_*$ and $=\lambda_*^{1/2}$ if $\lambda>2\lambda_*$ corresponding to Proposition \ref{choiceofm2}, $\bar\bu_1=(0,1,0,O(1))$ and $ \bar\bl_1=(1,0,0,0)$.
\end{Prop}
The  proof of this proposition will be given in Appendix \ref{L-Derivatives}.  The main dominant entry that will play an important role in the proof is $\frac{\partial \ell_0^f}{\partial L_0^i}$, which means that a small change in the semimajor of the binary on $\mathcal S_+$ will be stretched to a huge phase difference $\ell_0$ on the section $\cS_-$ after the global map. %This is an important mechanism used in \cite{X1,GHX}. We refer readers to \cite[Section 7.22]{GHX}  for a heuristic explanation. 
%The key point in this estimate is the $(2,1)$-entry $\frac{\partial \ell_0}{\partial L_0}=O(\mu^{-3/2}\lambda^{3/2}\chi)$, which implies that a small change in the initial $L_0$ will lead to a huge change in the final phase variable $\ell_0$.

From the last two propositions, it is clear that we have the transversality condition $\bl_1\cdot \bar\bu_1\neq 0$ and $\bar\bl_1\cdot \bu_1\neq 0$.  We define $\mathbb P:=\mathbb G\cR\mathbb L:\ \cS_-\to \cS_-$ the Poincar\'e map. Thus we have the following cone preservation property. 

\begin{Not}[Cone]
	\begin{enumerate}
		%\item Through out the whole text, we always assume  $0<\delta_0\ll\epsilon\ll\mu\ll1$. 
		\item Let $\bu_1,\ldots,\bu_n$ be a tuple of linearly independent vectors, we denote by $\mathcal C_\eta(\bu_1,\ldots,\bu_n)$ the \emph{$\eta$-cone} around $\bu_1,\ldots,\bu_n$ that is the set of vectors forming an angle at most $\eta$ with the plane span$\{\bu_1,\ldots,\bu_n\}$. 
	\end{enumerate}
\end{Not}

We fix a small number $\eta>0$ and choose $\dt$ small and $\chi$ large accordingly. 
\begin{Prop}\label{poincare1}Let $\bx\in \mathcal S_-$ be the initial condition for an orbit of the I4BP such that $\mathbb P(\bx)\in \mathcal S_-$. Suppose the assumptions of  Proposition \ref{PropDL1} and \ref{PropDG1} are satisfied along the orbit, then we have 
$$(d_\bx\mathbb P)\mathcal C_\eta(\bar\bu_1)\subsetneq \mathcal C_\eta(\bar\bu_1).$$
Moreover, for each $v\in \mathcal C_\eta(\bar\bu_1)$, we have $|d_\bx\mathbb P(v)|\geq f_{\mu,\lambda}\cdot e^{-Tv_*}\chi^2 |v|$. 
\end{Prop}

\subsection{Proof of Theorem \ref{ThmI4BP}}

We now give the proof of Theorem  \ref{ThmI4BP} assuming the propositions in the previous subsections.
We first show that there exists a Cantor set of  initial conditions such that the map $\mathbb P$ can be iterated for infinitely many steps with prefixed renormalization maps $\cR_{\lambda(\delta_i)}$, $i\in\mathbb{N}$. Then we show that such  initial conditions indeed lead to  super-hyperbolic orbits or non-collision singularities  (with double collisions regularized) depending on the choice of $\{\delta_i\}$.

\begin{proof}[Proof of Theorem \ref{ThmI4BP}]
From now on we only consider the phase points on $\cS_+$ such that the conditions in Proposition \ref{choiceofm2} and  \eqref{EqRenorm} are satisfied. Clearly, they form an open set on the section $\cS_+$.  Fix a small number $\eta>0$,
for any such a phase point $X=(L_0,\ell_0,L_2,\ell_2)\in \cS_+$, we introduce the cone $\mathcal C_\eta(\bu_1)\subset T_X\cS_+$.

We have the following non-degeneracy property for the global map.

\begin{Lm}\label{LmNondeg}
Let $\xi$ be an initial segment of length $\frac{1}{\chi}$ on the section $\cS_+$ with all its tangent vectors lying inside the cone $\mathcal C_\eta(\bar\bu_1)$  for each point $X$ on $\xi$.  Then its image under the global map  on the section $\cS_-$  wind  around the cylinder formed by $(L_0,\ell_0)$ along the direction $\ell_0$ many times.  In particular, it intersects transversely with the stable manifold $W_{\therefore}^s(O)$ when projected to $\cM_\therefore$.
 \end{Lm}
 We refer the reader to \cite[Section 5.5]{GHX} for an argument working exactly the same for proving this statement.

 {\bf Step 1:} {\it Construction of the Cantor set.}
 
Recall that we   fix $m_1\in(b_1,2)$, $0<\mu\ll1$,  $\dt_*(\leq \dt_0)$ sufficiently small, $\lambda_*=\epsilon^{-1}(2\delta_*)^{-1}$ and  $m_2$ to be $(m_1,\mu,\lambda_*)$-admissible. 
Let us choose a sequence  $\{\delta_i\}$ such that 
\begin{equation}\label{delta-i}
(\chi_0\prod_{k=0}^{i-1}\lambda_k)^{-D}\leqslant \delta_i\leqslant \delta_*,\; \lambda_0=1,\; \lambda_i=\lambda(\delta_i)=\epsilon^{-1}(2\delta_i)^{-1},\; i=1,2,\dots,
\end{equation}
where $\chi_0$ and $D$ are as in Proposition \ref{PropLocalMM}, and we assume either $\dt_i=\dt_*$ or $\dt_i\leq \dt_*/2$. Then the assumptions in Proposition \ref{PropLocalMM} is satisfied for $\delta=\delta_i$ and $\chi_{i-1}=\chi_0\prod_{k=0}^{i-1}\lambda_k$ after the application of the the $i$-th renormalization map $\cR_{\lambda_i}$, $i=1,2,3,\dots$. 

We start with an initial segment $\zeta_{0}$ on the section $\cS_-$ from Proposition \ref{PropLocalMM} with the chosen $\delta=\delta_1$.

We know  that for  the image $\mathbb L\zeta_0$ of the curve $\zeta_{0}$ under  the local map to  the section $\cS_+$, we can select a subsegment $\bar \zeta_0\subset \zeta_0$ such that each point in the image $\mathbb L\bar\zeta_0$ has $r$-component ranging from $(1-\mu)\delta_1$ to $(1+\mu)\delta_1$.  

Applying the renormalization map to $\mathbb L \bar\zeta_{0}$. Then the segment is rescaled to
 length of order one and the $r$-component has values in the interval  $(\frac{1}{2}(1-\mu)\eps^{-1},\frac{1}{2}(1+\mu)\eps^{-1})$ and the phase points on it satisfy the condition \eqref{EqRenorm}.  The renormalization map stretches the $L_0$ and $L_2$ variables
  with the same ratio $\lambda_1^{1/2}$ and leaves $\ell_0$, $\ell_2$  untouched. Therefore  after the renormalization, the segment $\mathcal R\mathbb L\bar\zeta_{0}$ has tangent vectors lying in the cone $\mathcal C_\eta(\bu_1)$ by Proposition \ref{PropDL1}, provided the rescaling factor $\lambda_1\gg \eta^{-1}$.

Next, we apply the global map to obtain the strong expansion in Proposition  \ref{poincare1}. The strong expansion shows that after applying the global map $\mathbb G$ and arriving at the section  $\cS_-$, the resulting curve
$\mathbb P\bar\zeta_0$ winds around the $\ell_0$-circle for many times. Then on the section $\cS_-$, we are in a position to pick from the image  $\mathbb P\bar\zeta_{0}$ segments with phase points satisfying the assumptions of Proposition  \ref{PropLocalMM} with $\delta=\delta_2$ as in  \eqref{delta-i}.   This involves deleting many open intervals from the  segment $\mathcal{R}\mathbb L\bar\zeta_{0}$, and correspondingly in the original segment $\bar\zeta_{0}$. We then repeat the above procedure to  the newly picked segments on the section $\cS_-$.  Repeating this procedure for infinitely many steps,   in the limit, we get a Cantor set $\mathcal{I}_{\zeta_{0},\{\delta_i\}}$ on the curve $\zeta_{0}$ as a result of deleting open intervals for infinitely many steps.

{\bf Step 2:} {\it Time estimate for super-hyperbolic orbit or non-collision singularity.}

We next show that each initial condition in the Cantor set $\mathcal{I}_{\zeta_{0},\{\delta_i\}}$ leads to a super-hyperbolic orbit or a non-collision singularity depending on the choice of $\{\delta_i\}$. %Note that here  we allow double collisions between the binary $Q_3$-$Q_4$ and those between $Q_1$ and $Q_2$ but does not allow any other type of double or triple collision.

Consider $\sx_0\in\mathcal{I}_{\zeta_{0},\{\delta_i\}}$ and assume that the initial distance between the particles $Q_2$ and the mass center of the binary $Q_3$-$Q_4$ is $\chi_0$ as in Proposition \ref{PropLocalMM}. 
Suppose $Q_1$ comes to a near triple collision with   the binary.
Then after the local map and the first renormalization $\cR_{\lambda_1}$, we have
the distance between $Q_1$ and $Q_2$ is $\chi_1= O(\lambda_1\chi_0)$ and $Q_1$ and $Q_2$ are moving to the right with  $|\dot Q_1|=O(1)$,  $|\dot Q_2|=O(\lambda_1^{-1/2})$ and  the binary $Q_3$-$Q_4$ moves to the left with speed of order~$O(1)$  by Proposition \ref{PropRenorm}. After catching up  and experiencing  a two-body interaction with $Q_2$, the particle $Q_1$ would be reflected and would get close to the binary, reaching    the section $\cS_-$, on which, by Proposition \ref{choiceofm2}, the distance between $Q_1$ and $Q_2$ is
$\chi_1=O(f_{\mu,\lambda_1}\lambda_1\chi_0),$
and   the time $t_1$  that (after the renormalization) the orbit spend between the sections $\cS_{+}$ and $\cS_-$, is 
$ t_1=O(f_{\mu,\lambda_1}\lambda_1\chi_0).$

In the original spacetime scale without the renormalization, the distance between $Q_1$ and $Q_2$ is 
$\bar\chi_1=O(f_{\mu,\lambda_1}\chi_0),$ and 
 the traveling time is  
 $
\bar  t_1=\begin{cases}
 	&O(\mu^{-1/2}\chi_0), \mathrm{\ if\ } \lambda_1=\lambda_*,\\
 	&O((\lambda_*/\lambda_1)^{1/2}\chi_0),\mathrm{\ if\ } \lambda_1\geq 2\lambda_*.
 \end{cases}
$
Then $Q_1$ and the binary $Q_3$-$Q_4$ would come to another near triple collision, after which the second renormalization $\cR_{\lambda_2}$ is performed. We need to update $\chi_0$ to $\chi_1$ and $\lambda_1$ to $\lambda_2$ for estimating $\chi_2$ and $t_2$ and update $\lambda_1$ to $\lambda_1\lambda_2$ for estimating $\bar\chi_2$ and $\bar t_2$, and repeat the argument.
 %Now at the section $\cS_+$, the distance between $Q_1$ and $Q_2$ is $\chi_2=O(\mu^{-1/2}\lambda_*^{1/2}\lambda_1\lambda_2\chi_0)$. Again by Proposition \ref{choiceofm2}, $Q_1$ would return to $\cS_-$ after a timespan $t_2=O(\mu^{-1}\lambda_*\lambda_1\lambda_2\chi_0)$. The distance between $Q_1$ and $Q_2$ when reaching $\cS_-$ is $O(\mu^{-1}\lambda_*\lambda_1\lambda_2\chi_0)$. In the original spacetime (with any renormalization), we have the time traveling time is $\bar t_2=O(\mu^{-1}\lambda_*\lambda_1^{-1/2}\lambda_2^{-1/2}\chi_0)$ and the distance between $Q_1$ and $Q_2$ is $O(\mu^{-1}\lambda_*\chi_0)$. 

We then perform the induction. Suppose that during the first $i$-steps, there are $j(i)$ times with $\lambda_k=\lambda_*$, $1\leq k\leq i$. Then we get from Proposition \ref{choiceofm2}, we see that in the original spacetime scale, the timespan between the $i$-th renormalization $\cR_{\lambda_i}$ and the $(i+1)$-th renormalization $\cR_{\lambda_{i+1}}$ is \begin{equation}\label{i-th-timespan}\bar t_i=O(\mu^{-j(i)/2}\mu^{\dagger}\prod_{k=1}^i(\lambda_*/\lambda_{k})^{1/2}\chi_0),
\end{equation}  
and the distance $\bar\chi_i=O(\mu^{-j(i)/2}\mu^{\dagger}\lambda_*^{i/2}\chi_0),$
where $\dagger=-\frac{1}{2}$ if $\lambda_i=\lambda_*$ and $\dagger=0$ if $\lambda_i\ge 2\lambda_*$.

 Therefore the total timespan before the $(i+1)$-th renormalization is estimated as
 \begin{equation}\label{time-i}
T_i=O(\sum_{\ell=1}^i\mu^{-j(\ell)/2}\mu^{\dagger}\prod_{k=1}^\ell(\lambda_*/\lambda_{k})^{1/2}\chi_0).
\end{equation} %and the distance between $Q_1$ and $Q_2$ is $D_i=O(\mu^{-i/2}\lambda_*^{i/2}\chi_0)$. 

So, we have the followings two cases:
\begin{enumerate}
	\item[ i)] if $\lim_{i\to\infty}T_i<+\infty$, then we obtain  non-collision singularities; 
		\item[ ii)] If $\lim_{i\to\infty}T_i=+\infty$, then we have  super-hyperbolic orbits.
\end{enumerate}

For i), an easy choice is to take $\lambda_k>2\lambda_*$ for all $k=1,2,\ldots.$ Thus we have $j(k)=0$ for all $k=1,2,\ldots$. We thus see that the sequence $T_i$ converges exponentially, since $O$ is uniform. We may also take $\lambda_k\to\infty$. The velocity of the particle $Q_1$ is estimated as $O(\sqrt{\prod\lambda_i})$. In this way, we make $\frac{|\dot Q_1(t_{k+2})|}{|\dot Q_1(t_{k})|}\to\infty. $
This gives the  statement ii) in Theorem \ref{ThmI4BP}. 
% taking $\delta_i=\mu\delta_*^{3/2}$. Recalling $\lambda_i=\epsilon^{-1}\delta_i^{-1}/2$,  then  in this case, $\lambda_*^{1/2}\lambda_i^{-1/2}=\delta^{1/4}_1$ and  by \eqref{time-i}, we have 
%$T_i=O(\sum_{k=1}^i\delta_*^{k/4}\chi_0)$. Clearly $\lim_{i\to\infty}T_i<\infty$ if $\delta_*$ is small enough. 
%By \eqref{delta-i}, we can choose $\{\delta_i\}$ such that $\lim_{i\to\infty}\delta_i=0$. This will leads to non-collision singularities with arbitrary  large acceleration rates. 

For ii), the choice of $\{\delta_i\}$ is more subtle. We give two valid selections here. The first one is to take $\delta_i=\delta_*$, hence, $\lambda_i=\lambda_*$. Then by \eqref{time-i}, $T_i=O(\sum_{k=1}^i\mu^{-k/2}\chi_0)\to\infty$ as $i\to\infty$ and since at time $T_i$ the distance between $Q_1$ and $Q_2$ is $O(\mu^{-i/2}\lambda_*^{i/2}\chi_0)$, we have $$\frac{\max\{|Q_j(t)|, j=1,2,3,4\}}{t}\sim\lambda_*^{i/2}=\lambda_*^{O(\frac{1}{|\log\mu|}\log t)},\mathrm{\ as\ }t\to\infty.$$ In fact this is the slowest rate of growths for the super-hyperbolic orbits constructed here.  

A second choice of $\{\delta_i\}$ is of the form $\delta_i=\mu^{N_i}\delta_*$, satisfying  \eqref{delta-i}. Here $N_i\in \mathbb{N}$, $N_1\geqslant1$ and if for some $i$,  $N_i\geqslant1$, then $N_k=0$, $k=i,\dots i+N_i$ and $N_{i+N_i+1}\geqslant1$.  Let $\{i_k,k=1,\dots\}$ be the increasing sequence such that $N_{i_k}\geqslant1$ and $N_i=0$ if $i\not\in\{i_k\}$. %Then by \eqref{i-th-timespan} the timespan between the $i_k$-th and $i_{k+1}$-th renormalization is
%$O(\sum_{j=0}^{N_{i_k}}\sqrt\mu^{N_{i_k}-j})=O(1)$. 
Then by \eqref{i-th-timespan}, we get  $\lim_{i\to\infty}T_i\to\infty$ since $T_i$ grows by a multiple of $\sqrt\mu$ if $i\not\in\{i_k\}$. %Moreover, in this case, we have estimate $\frac{\max\{|Q_j(t)|, j=1,2,3,4\}}{t}=O(\mu^{-f(t)}\lambda_*^{f(t)})$ where as $t\to\infty$, $f(t)=O(i_{\lfloor t\rfloor})$.   It is not difficult to see that with suitable choice of $\{i_k\}$, we can make $f(t)$ larger than any given function as $t\to\infty$.
  
We thus prove the assertion of Theorem \ref{ThmI4BP}. 
\end{proof} 
%{\color{red}Here I delete some arguments on arbitrarily fast growth, since there is a constraint by sojourn time, the issue is a bit complicated. }

\begin{Rk}\label{mass-ratio-remark} A straightforward analysis of the construction we can see that while the non-collision singularity is allowed for all $m_2>\frac{m_1^2+2m_1}{2-m_2}$, the super-hyperbolic orbit is possible only for $0<m_2-\frac{m_1^2+2m_1}{2-m_2}\ll1$. The latter restriction is essential in our construction for slowing down the the speeds of the particles. In \cite{SX}, the existence of super-hyperbolic orbits was presented for a much wider range of mass ratios in the Mather-McGehee model of collinear four-body problem. Infinitely many double collisions are not avoidable there. % just as Mather-McGehee model
\end{Rk}
 \section{Estimate of the return time}\label{SJacobi}
 
In this section, we first write down the Hamiltonian of the F4BP and study the transformation between the left  and right Jacobi coordinates in the first two subsections.  
After that, we give the proof  of Proposition \ref{choiceofm2}. %, which will be provided in this section. The point of the proof is that by choosing the masses $m_1$ and $m_2$ carefully ($(m_1,\mu,\lambda_*)$-admissible), we shall arrange that the returning momentum $\bp_1$ to be of order $(\mu/\lambda_*)^{1/2}$, which is made to be so small that $Q_1$ should take very long time to return to the next near triple collision. 

\subsection{The Hamiltonian systems in Jacobi coordinates}
The original Hamiltonian has the form
$$H(P,Q)=\frac{P_1^2}{2m_1}+\frac{P_2^2}{2m_2}+\frac{P_3^2}{2}+\frac{P_4^2}{2}-\frac{m_1m_2}{|Q_1-Q_2|}-\sum_{i=1,2}^{j=3,4}\frac{m_i}{|Q_i-Q_j|}-\frac{1}{|Q_3-Q_4|}.$$
In the (left) Jacobi coordinates, defined in Section \ref{subsection-jacobi}, the Hamiltonian reads
\begin{equation}\label{hamilton-left1}
	\begin{aligned}
		&H_{:\cdot\cdot}(\bx,\bp)=\left(\frac{\bp_1^2}{2M_1}+\frac{\bp_0^2}{2M_0} -\frac{m_1}{|\bx_1-\frac{\bx_0}{2}|}-\frac{m_1}{|\bx_1+\frac{\bx_0}{2}|}-\frac{1}{|\bx_0|}\right) \\
		&+\left(\frac{\bp_2^2}{2M_2}  -\frac{m_1m_2}{|\bx_2-\frac{2 \bx_1}{m_1+2} |}-\frac{m_2}{|\bx_2+\frac{m_1\bx_1}{m_1+2}-\frac{\bx_0}{2}|}-\frac{m_2}{|\bx_2+\frac{m_1\bx_1}{m_1+2}+\frac{\bx_0}{2}|}\right)
	\end{aligned}
\end{equation}
where
$M_0=\frac12,\  M_1=\frac{2m_1}{m_1+2},\  M_2=\frac{m_2(m_1+2)}{m_1+m_2+2}.$ The first parenthesis is a three-body problem in Jacobi coordinates and the second one is a perturbed two-body problem. We shall used this Hamiltonian when the orbits of the system lie between the two sections $\cS_-$ and $\cS_+$ and close to triple collision. 
%
% We will work with two different forms of Hamiltonians. When $Q_1$ is close to the binary $Q_3$-$Q_4$  (local map), we write the Hamiltonian as
%\begin{equation}\label{EqHamLoc}
%\begin{aligned}
%H_{:\cdot\cdot}(\bx,\bp)&=H_\therefore+H_2\\
%H_\therefore&=\frac{\bp_1^2}{2M_1}+\frac{\bp_0^2}{2M_0} -\frac{m_1}{|\bx_1-\frac{\bx_0}{2}|}-\frac{m_1}{|\bx_1+\frac{\bx_0}{2}|}-\frac{1}{|\bx_0|}:=\frac{\bp_1^2}{2M_1}+\frac{\bp_0^2}{2M_0} +V,\\
% H_2&=\frac{\bp_2^2}{2M_2} -\frac{k_2}{|\bx_2|}+U_2, \quad k_2=(m_1+2)m_2,\\
%   U_2&=\frac{k_2}{|\bx_2|}-\frac{m_1m_2}{|\bx_2-\frac{2 \bx_1}{m_1+2} |}-\frac{m_2}{|\bx_2+\frac{m_1\bx_1}{m_1+2}-\frac{\bx_0}{2}|}-\frac{m_2}{|\bx_2+\frac{m_1\bx_1}{m_1+2}+\frac{\bx_0}{2}|}.
%\end{aligned}
%\end{equation}
%Note that the subsystem $H_\therefore$ is a three-body problem in Jacobi coordinates and $H_2$ is a Kepler problem perturbed by $U_2$.

When $Q_1$ is far away from the binary  and $Q_2$, i.e. the orbit has not reached the section $\cS_-$ and have exited the section $\cS_+$, we choose a different way of grouping terms and use the following form of Hamiltonian system, which is a perturbation of three Kepler problems
\begin{equation} \label{EqHamGlob}
	\begin{aligned}
		H_{:\cdot\cdot}(\bx,\bp)&=\left(\frac{\bp_0^2}{2M_0}-\frac{k_0}{|\bx_0|}\right)+\left(\frac{\bp_1^2}{2M_1}-\frac{k_1}{|\bx_1|}\right)+\left(\frac{\bp_2^2}{2M_2}-\frac{k_2}{|\bx_2|}\right) +U_{01}+U_2\\
		k_0&=1,\ k_1=2m_1,\; k_2=(m_1+2)m_2,\\
		U_{01} &=\left(\frac{k_1}{|\bx_1|}-\frac{m_1}{|\bx_1-\frac{\bx_0}{2}|}-\frac{m_1}{|\bx_1+\frac{\bx_0}{2}|}\right),\\
		U_2&=\frac{k_2}{|\bx_2|}-\frac{m_1m_2}{|\bx_2-\frac{2 \bx_1}{m_1+2} |}-\frac{m_2}{|\bx_2+\frac{m_1\bx_1}{m_1+2}-\frac{\bx_0}{2}|}-\frac{m_2}{|\bx_2+\frac{m_1\bx_1}{m_1+2}+\frac{\bx_0}{2}|}.
	\end{aligned}
\end{equation}
We denote  the energy \begin{equation}\label{energy-formula1}
E_i=\frac{|\bp_i|^2}{2M_i}-\frac{k_i}{|\bx_i|}\end{equation}
 for the Kepler problem $(\bx_i,\bp_i)$,  $ i=0,1,2.$

When $Q_1$ is much  closer to $Q_2$ than the binary $Q_3$-$Q_4$, we use the right Jacobi coordinates.  
\begin{equation}\label{EqJacobiR}
\begin{cases}
	\bx_0^R=q_3\\
	\bx_1^R=q_1-q_2\\
	\bx_2^R=\frac{m_1q_1+m_2q_2}{m_1+m_2}-\frac{1}{2}q_3,
\end{cases}\begin{cases}
	\bp_0^R=P_3+\frac{P_1}2+\frac{P_2}2\\
	\bp_1^R=-\frac{m_1}{m_1+m_2}P_2+\frac{m_2}{m_1+m_2}P_1\\
	\bp_2^R=P_2+P_1.
\end{cases}	
\end{equation}
Here, $(\bx_0^R,\bp_0^R)$ is the same as $(\bx_0,\bp_0)$, the variables $(\bx_1^R,\bp_1^R)$ describes the relative motion between $Q_1$ and $Q_2$, and $(\bx_2^R,\bp_2^R)$ does that between the mass center of $Q_1$-$Q_2$ and the mass center of the binary $Q_3$-$Q_4$. 

For the above used left Jacobi coordinates \eqref{EqJacobi} as well as the constants $M_i,k_i,m_i$, we shall put a superscript $L$ to avoid confusion with the right Jacobi coordinates \eqref{EqJacobiR} and relevant constants below. However, if there is no danger of confusion, we omit the superscripts. 

To decide which set of coordinates to use, we introduce a middle section $\cS^M=\{|\bx_1^R|=\beta\chi\}$ where $\beta=\lambda^{-1/2}$ and $\lambda$ is the constant defining the renormalization. To the left (respectively right) of the section $\cS^M$, we use the left (respectively right) Jacobi coordinates. %We change the coordinates from the left one to the right one on the section $\cS^M$ when $Q_1$ arrives at $\cS^M$ from the left side and vise versa. 

%in which the Hamiltonian reads
%\begin{equation}
%\begin{aligned}
%&H_{:\cdot\cdot}(\bx^R,\bp^R)=\frac{1}{2M_1^R}|\bp_1^R|^2+\frac{1}{2M_0^R}|\bp_0^R|^2+\frac{1}{2M_2^R}|\bp_2^R|^2\\
%&-\frac{m_1m_2}{|\bx_1^R|}-\frac{1}{|\bx_0^R|}-\frac{m_1}{|-\frac{1}{2}\bx_0^R-\bx_2^R+\frac{m_2}{m_1+m_2}\bx_1^R|}-\frac{m_1}{|\frac{1}{2}\bx_0^R-\bx_2^R+\frac{m_2}{m_1+m_2}\bx^R_1|}\\
%&-\frac{m_2}{|-\frac{1}{2}\bx_0^R-\bx_2^R-\frac{m_1}{m_2+m_2}\bx_1^R|}-\frac{m_2}{|\frac{1}{2}\bx_0^R-\bx_2^R-\frac{m_1}{m_2+m_1}\bx_1^R|},
%\end{aligned}
%\end{equation}where $M_0^R=\frac{1}{2}$, $M_1^R=\frac{m_1m_2}{m_1+m_2}$, $M_2^R=\frac{2m_2(m_1+m_2)}{m_1+m_2+2}$.
%
In the right case, $Q_1$ and $Q_2$ form a close pair and they are far away from the pair $Q_3$-$Q_4$, we write the Hamiltonian as the following %to emphasis this structure,
\begin{equation}\label{3pairHam}
	\begin{aligned}
		H_{:\cdot\cdot}(\bx^R,\bp^R)&=\Big(\frac{|\bp_0^R|^2}{2M_0^R}-\frac{k_0^R}{|\bx^R_0|}\Big)+\Big(\frac{|\bp_1^R|^2}{2M_1^R}-\frac{k_1^R}{|\bx_1|}\Big)+\Big(\frac{|\bp^R_2|^2}{2M_2^R}-\frac{k_2^R}{|\bx_2^R|}\Big)+U_{12},\\
		k_0^R&=1,\;k_1^R=m_1m_2,\; k_2^R=2(m_1+m_2),\\
		U_{12}&=\frac{k_2^R}{|\bx_2^R|}-\frac{m_1}{|-\frac{1}{2}\bx^R_0-\bx^R_2+\frac{m_2}{m_1+m_2}\bx^R_1|}-\frac{m_1}{|\frac{1}{2}\bx^R_0-\bx^R_2+\frac{m_2}{m_1+m_2}\bx^R_1|}\\
		&\quad\quad-\frac{m_2}{|-\frac{1}{2}\bx^R_0-\bx^R_2-\frac{m_1}{m_2+m_2}\bx^R_1|}-\frac{m_2}{|\frac{1}{2}\bx^R_0-\bx^R_2-\frac{m_1}{m_2+m_1}\bx^R_1|}.
	\end{aligned}
\end{equation}
Note that this is a system of three almost decoupled Kepler problems perturbed by the potential $U_{12}$.  

%We then prove Proposition \ref{choiceofm2} in Section \ref{proofpm2} and Lemma \ref{lm-avoid-dc} in Section \ref{avoid-dc}. 
\subsection{Transition from the left to the right}\label{STransition}
We now give the linear transform changing the left Jacobi coordinates to the right Jacobi coordinates and vise versa.   
Let us denote $(\bx_0^L,\bx_1^L,\bx_2^L)^T=\mathcal{L}(q_3,q_1,q_2)^T$, \ $(\bx_0^R,\bx_1^R,\bx_2^R)^T=\mathcal{R}(q_3,q_1,q_2)^T$, $(\bp_0^L,\bp_1^L,\bp_2^L)=\mathcal{L}^{-T}(P_3,P_1,P_2)^T$, and
$(\bp_0^R,\bp_1^R,\bp_2^R)=\mathcal{R}^{-T}(P_3,P_1,P_2)^T$, where $-T$ means \emph{transpose inverse}.
Then
\[\mathcal{L}=\left(\begin{array}{ccc}1&0&0\\
	-\frac{1}{2}&1&0\\
	\frac{-1}{m_1+2}& -\frac{m_1}{m_1+2}&1\end{array}\right),\quad \mathcal{R}=\left(\begin{array}{ccc}1&0&0\\
	0&1&-1\\
	-\frac{1}{2}&\frac{m_1}{m_1+m_2}&\frac{m_2}{m_1+m_2}\end{array}\right).\]
Therefore,
$$\mathcal{L}\mathcal{R}^{-1}=\left(\begin{array}{ccc}1&0&0\\
	0&\frac{m_2}{m_1+m_2}&1\\
	0&-\frac{m_1(2+m_1+m_2)}{(2+m_1)(m_1+m_2)}&\frac{2}{2+m_1}\end{array}\right),
\quad\mathcal{R}\mathcal{L}^{-1}=\left(\begin{array}{ccc}1&0&0\\0&\frac{2}{m_1+2}&-1\\0&\frac{m_1(m_1+m_2+2)}{(m_1+m_2)(2+m_1)}&\frac{m_2}{m_1+m_2}\end{array}\right),$$
Hence we have the following explicit formulas for the transitions of coordinates, 
$$(\bx^L_0,\bx^L_1,\bx^L_2)^T=\mathcal L\mathcal R^{-1}(\bx^R_0,\bx^R_1,\bx^R_2)^T, \quad (\bp^L_0,\bp^L_1,\bp^L_2)^T=(\mathcal R\mathcal L^{-1})^T(\bp^R_0,\bp^R_1,\bp^R_2)^T$$
\begin{equation}\label{R-L-transit}\begin{cases}\bx_1^L=\frac{m_2}{m_1+m_2}\bx_1^R+\bx_2^R,\\
		\bx_2^L=-\frac{m_1(2+m_1+m_2)}{(2+m_1)(m_1+m_2)}\bx_1^R+\frac{2}{2+m_1}\bx_2,\end{cases}\quad
	\begin{cases}\bp_1^L=\frac{2}{2+m_1}\bp_1^R+\frac{m_1(2+m_1+m_2)}{(m_1+2)(m_1+m_2)}\bp_2^R,\\
		\bp_2^L=-\bp_1^R+\frac{m_2}{m_1+m_2}\bp_2^R,\end{cases}\end{equation}
and
\begin{equation}\label{L-R-transit}\begin{cases}
		\bx_1^R=\frac{2}{m_1+2}\bx_1^L-\bx_2^L,\\
		\bx_2^R=\frac{m_1(m_1+m_2+2)}{(m_1+m_2)(2+m_1)}\bx_1^L+\frac{m_2}{m_1+m_2}\bx_2^L,\end{cases}\quad\begin{cases}\bp_1^R=\frac{m_2}{m_1+m_2}\bp_1^L-\frac{m_1(m_1+m_2+2)}{(m_1+2)(m_1+m_2)}\bp_2^L,\\
		\bp_2^R=\bp_1^L+\frac{2}{2+m_1}\bp_2^L.\end{cases}\end{equation}
\subsection{Proof of Proposition \ref{choiceofm2}}\label{proofpm2}
We consider the scenario where 
the triple $Q_1$-$Q_3$-$Q_4$ is ejected from the triple collision such that $|\bp_0|$ and $|\bp_1|$ are sufficiently large and $|\bp_2|\sim1$. Suppose the distance between $Q_1$ and $Q_2$ is $\chi$, which is as large as we wish. After renormalization $\cR_{\lambda}$, $\bp_0$ and $\bp_1$ are rescaled to size of order $O(1)$ and $|\bp_2|=O(\beta)$, where we denote $\beta=\lambda^{-1/2}$.  Then $Q_1$ would catch up with $Q_2$. 
For orbits traveling between  sections $\cS_\pm$ and $\cS^M$, we have the following lemma. 

\begin{Lm} \label{L01-bound}There exists $C>1$ such that  
	\begin{enumerate}
		\item 
		for  the I4BP orbits $\gamma(t):[0,T]\to\cM$ with $\gamma(0)\in\cS_+$ and $\gamma(T)\in\cS^M$ such that \[\{L_1(0),\; L_0(0)\; L_2(0)\beta\}\subset[C^{-1},C],\]
		we have $T=O(\chi)$ and 
		\[\max_{t\in[0,T]}\{ |L_1(t)-L_1(0)|,\; |L_0(t)-L_0(0)|, \;| L_2(t)\beta-L_2(0)\beta|\}=O(\epsilon^2);\]
		\item for I4BP orbits $\gamma(t):[0,T]\to\cM$ with $\gamma(0)\in \cS^M$ and $\gamma(T)\in \cS^M$ such that along $\gamma$ there exists a double collision between $Q_1$ and $Q_2$, and   
		$$\{L_0^R(0), L_1^R(0),  L_2^R(0)\}\subset[{C}^{-1},C],$$ we have $T=O(\beta \chi)$ and 
		$\max_{t\in[0,T]}\{|L_i^R(t)-L^R_i(0)|, i=0,1,2\}=O(\frac{1}{\chi}). $
	\end{enumerate}
\end{Lm}
We refer the reader to Lemma 6.6 of \cite{X1} for  a proof. The equations of motion are given in \eqref{EqHamLimit}. 
The lemma follows from integrating equation \eqref{EqHamLimit}  combined with  a simple bootstrap argument. The argument goes as follows. We start by assuming  the oscillations of $L_i,\ i=0,1,2$ are bounded by a generous constant $C\gg \eps^2$. Then this implies that $r_0$ is bounded, and $r_1$ and $\ell_1$ grows at least linearly in time. Then integrating the estimates in Appendix \eqref{EqHamLimit} over time of order $\chi$, we get that the oscillations of $L_i$ are indeed $O(\eps^2)$, which is the value of $(\ell_1^i)^{-2}=\int_{\ell_1^i}^\chi\frac{1}{\ell_1^3}$ evaluated on the section $\cS_+$. The same  argument works for statement (2) with the estimates \eqref{RightHam} in Appendix \ref{L-Derivatives}. 

We are now ready to prove Proposition \ref{choiceofm2}. 
\begin{proof}[Proof of Proposition \ref{choiceofm2}]
	By assumption, after the renormalization $\cR_{\lambda_\delta}$, $\lambda=\lambda(\delta) $, on the section $\cS_+$, we have
	$$|\bp_1^L|\in\sqrt{2M_1^L}(1-\mu,1+\mu),\; \quad |\bp_2^L|\in\frac{1}{\sqrt{\lambda}}\sqrt{2M_1^L}\frac{2m_2}{m_1+m_2}(1-2\mu, 1+2\mu).$$
	By Lemma \ref{L01-bound}, we know that when the orbits reach  the middle section $\cS^M:=\{|\bx_1^R|=\beta\chi\}$, we have 
	\begin{equation}\label{p-bound}
	\begin{split} &|\bp_1^L|\in\sqrt{2M_1^L}(1-\mu,1+\mu)+O(\epsilon^2),\\
	& |\bp_2^L|\in\frac{1}{\sqrt{\lambda}}\sqrt{2M_1^L}\frac{2m_2}{m_1+m_2}(1-2\mu, 1+2\mu)+O(\frac{\epsilon^2}{\sqrt{\lambda}}).
	\end{split}
	\end{equation}
	On the section $\cS^M$,  using \eqref{L-R-transit}, we change from the left Jacobi coordinates to the right Jacobi coordinates and have
	$\bp_1^R=O(1)$  and 
	$\bp_2^R=O(1)$. Both are pointing to the right. 
	Then by Lemma \ref{L01-bound}, when the orbits return to the section $\cS^M$, the values $|\bp_1^R|$ and $|\bp_2^R|$ undergoes an oscillation of order $O(\eps^2)$ with $\bp_1^R$ changing  direction from right to left.  Switching back to the left Jacobi coordinates, using \eqref{L-R-transit} and \eqref{R-L-transit}, we then obtain on the section $\cS^M$, 
	\begin{equation*}
		\begin{aligned}
	\bar\bp_2^L&=\frac{2m_2}{m_1+m_2}\bp_1^L-\frac{m_1(m_1+m_2+2)-2m_2}{(m_1+2)(m_1+m_2)}\bp_2^L,\\
	\bar\bp_1^L&=\frac{m_1(m_1+m_2+2)-2m_2}{(m_1+2)(m_1+m_2)}\bp_1^L+\frac{4m_1(m_1+m_2+2)}{(m_1+2)^2(m_1+m_2)}\bp_2^L,
		\end{aligned}
	\end{equation*}
	where $\bp_1^L$ and $\bp_2^L$ are as in \eqref{p-bound} up to an order $O(\eps^2)$ oscillation). Also note that now we have  $|\bx_1^L|=O(\chi)$ and by assumption $\eps<\mu$.  Therefore we have on the section $\cS^M$, $|\bar\bp_2^L|\in\frac{2m_2\sqrt{2M_1^L}}{m_1+m_2}(1-2\mu,1+2\mu)$ with direction pointing to the right and $\bar\bp_1^L$ points to the left (thus we have $m_1(m_1+m_2+2)-2m_2<0$)
	 with length
	\begin{equation*}
		\begin{aligned}
	|\bar\bp_1^L|&=\sqrt{2M_1^L}\left(\frac{2m_2-m_1(m_1+m_2+2)}{(m_1+2)(m_1+m_2)}-\lambda^{-1/2}\frac{8m_1m_2(m_1+m_2+2)}{(m_1+2)^2(m_1+m_2)^2}\right)\\
	&+O(\mu+\eps^2) \frac{m_1(m_1+m_2+2)-2m_2}{(m_1+2)(m_1+m_2)}+O(\lambda^{-1/2}(\mu+\eps^2)).		
		\end{aligned}
	\end{equation*}
Now we have two cases from the assumption: 
\begin{enumerate}
	\item $\lambda=\lambda_*$ and  $m_2$ is $(m_1,\mu,\lambda_*)$-admissible,
	\item  $\lambda\geq 2\lambda_*$ and  $m_2$ is $(m_1,\mu,\lambda_*)$-admissible, but not $(m_1,\mu,\lambda)$-admissible. 
\end{enumerate}

	In case (1), by equation \eqref{admissible-m2}), we get that the first row of the above $|\bar\bp_1^L|$ lies in the interval $[1,10]\cdot\sqrt{\mu\lambda_*^{-1}}\sqrt{2M_1^L}$. Note also that \eqref{admissible-m2} implies that $\frac{m_1(m_1+m_2+2)-2m_2}{(m_1+2)(m_1+m_2)}=O(\lambda_*^{-1/2})$, thus we get the second row of the above $|\bar\bp_1^L|$  is estimated as $O(\mu\lambda_*^{-1/2})$. Thus we get  $|\bar\bp_1^L|\in[\frac{1}{2}, 11]\cdot\mu^{1/2}\lambda_*^{-1/2}\sqrt{2M_1^L}$.  
	
	In case (2), we have that $\frac{2m_2-m_1(m_1+m_2+2)}{(m_1+2)(m_1+m_2)}$ is $O(\sqrt{\mu\lambda^{-1}})$ close to $\lambda_*^{-1/2}\frac{8m_1m_2(m_1+m_2+2)}{(m_1+2)^2(m_1+m_2)^2}$ since $m_2$ is $(m_1,\mu,\lambda_*)$-admissible. Thus we get that $|\bar\bp_1^L|$ is of order $\lambda_*^{-1/2}$ for all $\lambda>2\lambda_*$. 
	
In both cases, with the estimate of the velocity, we get the corresponding estimate of time spans and distances.	The assertion of the proposition now follows. 
 \end{proof}

\section{The full four-body problem and proof of Theorem \ref{ThmF4BP} }\label{SMainProof}
In this section, we prove Theorem \ref{ThmF4BP} for the full four-body problem(F4BP). The definitions of the sections $\cS_\pm$, the local, global and renormalization maps for the F4BP are similar to the isosceles case. We will show in this section that the super-hyperbolic orbit of the F4BP is obtained by slightly perturbing the I4BP orbit from Theorem~\ref{ThmI4BP}. The main issue is to control the extra variables
$(G_0,g_0,G_1,g_1,G_2,g_2)$ as well as excluding the possibility of collisions.  
\subsection{The coordinates and the diagonal form}
For each pair of Jacobi coordinates $(\bx_i,\bp_i),\ i=0,1,2,$ we convert it into Delaunay coordinates $(L_i,\ell_i,G_i,g_i)$. 
We always assume the total angular momentum is vanishing, i.e.
$G_0+G_1+G_2=0.$ Since the Hamiltonian does not depend on each individual angle $g_i$ but depends only on the relative angles, this reduces the number of degrees of freedom by 1.  By further restricting to the zero energy level and to the Poincar\'e sections, we can reduce two more variables  when defining the local and global maps as we did in the last section. We define the Poincar\'e sections $\cS_\pm$ in the same way as Definition \ref{DefSection}. The variables to be removed are chosen to be $L_1$ and $\ell_1$ on the sections $\cS_\pm$.  Thus we have the following set of variables to describe the dynamics  $$X:=(L_0,\ell_0,L_2,\ell_2), \ Y:=(G_0,\mathsf g_{01}:=g_0-g_1-\pi/2,G_2,\mathsf g_{21}:=g_2-g_1),$$
the variables on the sections $\cS_\pm$. 
 We shall estimate the matrices $d\mathbb G=\frac{\partial (X,Y)|_{\cS_-}}{\partial (X,Y)|_{\cS_+}}$, and $\ d\mathbb L=\frac{\partial (X,Y)|_{\cS_+}}{\partial (X,Y)|_{\cS_-}}$. As we state in the next lemma,  the derivative matrix has a product structure in the isosceles limit, i.e. we evaluate the derivative matrix in the limit $|Y|\to 0$ (along an I4BP orbit). 
\begin{Lm}\label{LmDiagonal} The derivative matrix $\frac{\partial (X,Y)|_{\cS_-}}{\partial (X,Y)|_{\cS_+}}$  evaluated along an I4BP orbit has the following structure $$\frac{\partial (X,Y)|_{\cS_-}}{\partial (X,Y)|_{\cS_+}}=\left[\begin{array}{ccc}
\frac{\partial X|_{\cS_-}}{\partial X|_{\cS_+}}&0\\
0&\frac{\partial Y|_{\cS_-}}{\partial Y|_{\cS_+}}\\
\end{array}\right].$$
The same diagonal form of the derivative matrix also holds for the local map hence for the Poincar\'e return map. 
\end{Lm}
We refer the reader to \cite[Section 9.1]{GHX} for a proof of this statement. 

%The first diagonal block  here is estimated in Proposition \ref{PropDL1} for local map and \ref{PropDG1} for global map. We next show that how to estimate the second diagonal block. 

\subsection{The I3BP near infinity}
%We shall also consider solutions of the isosceles three-body problem with the binary and the third body moving off to infinity in opposite directions. We use another coordinate change of McGehee.

%We may view the entire $y$-axis as a normally hyperbolic invariant manifold.

%In this coordinates, we get that the infinity is a normally hyperbolic invariant manifold topologically a cylinder parametrized by $(L_0,\ell_0)$.
We will need the following special orbits $\gamma_{I}$ and $\gamma_O$ as guides of our orbits when experiencing near triple collisions. Follow our discussion before (see Proposition \ref{choiceofm2}), after the two-body interaction with $Q_2$,  we need that the relative speed between $Q_1$ and the binary $Q_3$-$Q_4$ is of order $\lambda_*^{-1/2}$.  Suppose the triple $Q_1$-$Q_3$-$Q_4$ goes to triple collision, then  we expect that for  $\lambda_*$ sufficiently large, the orbit approaching triple collision stays very close to certain orbit $\gamma_I$ on the stable manifold $W^s_{\therefore}(O)$ with the feature that it comes from the infinity with almost zero initial velocity of $Q_1$. Therefore, to guide our incoming orbits of near triple collision,  we choose the orbit $\gamma_I$ to be an orbit on $W^s_\therefore(O)$ approaching infinity parabolically as $t\to-\infty$, whose existence has been well studied (c.f. \cite{SM}).

\begin{Prop}[\cite{SM}] \label{ThmGammaI} There exists an orbit of the I3BP $Q_1$-$Q_3$-$Q_4$ lying on the stable manifold $W_\therefore^s(O)$ of the Lagrange fixed point $O$ and satisfying the boundary condition $\bx^L_1(t)\to \infty$ and $\bp^L_1(t)\to0$  as $t\to-\infty$.
\end{Prop}

Next, if the triple $Q_1$-$Q_3$-$Q_4$ stays close to the right lower Lagrange fixed point for a long time and escapes to the right upper arm shadowing the unstable manifold on $\cM_0$, we expect $r$ remains small for a long time. So we use $\gamma_O$, that is, the unstable manifold of the right lower Lagrange fixed point on the collision manifold $\cM_0$, to guide the orbits leaving the triple collision. 
\begin{Not}
	\begin{enumerate}
		\item We will call the orbit in Proposition \ref{ThmGammaI} $\gamma_I$ and the unstable manifold $W^u_0$ in Notation \ref{NotWs} $\gamma_O$.
		\item The orbits $\gamma_I$ and $\gamma_O$ are invariant under time translation, thus for convenience, for $\gamma_I$, we assume $\gamma_I(0)\in \cS_-$ and for $\gamma_O$ we assume $\gamma_O(0)\in \cS_+$, after time translations if necessary.
	\end{enumerate}
\end{Not}

\begin{Not}
	\begin{enumerate}
		\item Denote by $\cM$ the phase space of  the F4BP. 
		\item Denote by $\cM_\therefore$ the phase space of the three-body problem $($F3BP$)$ $Q_1$-$Q_3$-$Q_4$, which has six dimensions after reducing the momentum conservation and angular momentum conservation.
		\item Denote by $\pi_\therefore$ the projection to the $\cM_\therefore$ subsystem.
	\end{enumerate}
\end{Not}
The next proposition shows that $\gamma_I$ indeed serves as a good guide for the incoming orbits. 
\begin{Prop}[Proposition 3.17 of \cite{GHX}]\label{PropDevaneyOrbit}
	Let $\dt'$ and $\chi$  satisfy $\dt'>\chi^{-D}$ for some constant $D>1$.
	Let $\gamma:\ [0,T]\to \cM$ be an orbit of the I4BP such that
	\begin{enumerate}
		\item  $|\bx_1(\gamma(0))|\sim\chi$, $|\bp_1(\gamma(0))|\leq \delta'$;
		\item
		$\gamma(T)\in \cS_-$ satisfies $\mathrm{dist}(\pi_\therefore(\gamma(T)),W^s_\therefore(O))<\dt'$.
	\end{enumerate}
	Then as  $\dt'\to 0$, we have
	$\mathrm{dist}(\pi_\therefore(\gamma(T)),\gamma_I(0))\to 0$.
\end{Prop}
%Since  we only use $\gamma_I$ to approximate the real incoming orbits between the segment  $\cS_-$  and a small neighborhood of the Lagrangian fixed point $O$, by  Proposition \ref{PropDevaneyOrbit}, this approximation can be made as precise as we want for sufficient small $\delta=O(\mu^{1/2}\lambda^{-1/2})$. 
\subsection{The local map}\label{SSLoc}
In this section, we study the derivative of the local map. %The setting is the same as that in \cite{GHX}, thus we only present the main statements and refer readers to \cite[Section 8]{GHX} for the proofs. 
When the subsystem $Q_1$-$Q_3$-$Q_4$ is near triple collision, we shall treat the F4BP as a small perturbation of the product system of the F3BP and a Kepler problem $(\bx_2,\bp_2)$.  Near triple collision of the  $Q_1$-$Q_3$-$Q_4$ subsystem, we update the meanings of $X$ and $Y$ as  $X=(r,v,\psi,w)$ the blowup coordinates that we used in the I3BP, and we introduce in addition the new scaling invariant variables $Y:=(w_0=r^{-1/2}G_0,\mathsf g_{01},w_2=r^{-1/2}G_2,\mathsf g_{21})$. 

We first consider the F3BP $Q_1$-$Q_3$-$Q_4$. Assuming that the total angular momentum $G_0+G_1=0$, the system has three degrees of freedom, and we introduce the variables $w_0,\mathsf g_{01}$ in addition to the blowup coordinates $(r,v,\psi,w)$ of the I3BP. The diagonal structure Lemma~\ref{LmDiagonal} also holds for the F3BP, so we can talk about the second diagonal block that is $2\times 2$ by projecting to the $(w_0,\mathsf g_{01})$-plane. %The phase space of the I3BP is an invariant submanifold of codimension 2 of that of the F3BP with $\mathsf g_{01}=\pi/2$ and $w_0=0$. 
\begin{Lm}\label{LmDF3BPLoc}
When linearized at each Lagrange fixed point, the equations of motion of the F3BP projected to the $(w_0,\mathsf g_{01})$-plane has  one  positive $\mu_u>0$ and one negative $\mu_s<0$ eigenvalues. 
\end{Lm}
 %Therefore, the Lagrange fixed point $O_1$ of the F3BP has three dimensional stable manifold $W_3^s(O_1)$ (containing $\gamma_I$) and two dimensional unstable manifold $W_3^u(O_1)$ (containing $\gamma_O$).

%For the local map, the second diagonal is the same as $d(\widehat{\cR \mathbb L})$.  

\begin{Not}
 Denote by $E_3^s(\gamma_I(\tau)),\ \tau\in [0,\infty)$ the tangent vector $($normalized to have length $1)$ of $W_\therefore^s(\gamma_I(\tau))$ at $(0,0)$ in the $(w_0,\mathsf g_{01})$-plane $($for the local map$)$. Similarly, we define $E_3^u(\gamma_O(\tau))$, $\tau\in (-\infty,0]$.
\end{Not}

We next consider the F4BP taking into further consideration of $Q_2$.  Assuming the total angular momentum vanishes, we use coordinates $Y$ in addition to the blowup coordinates for the I3BP. Again by the diagonal structure Lemma \ref{LmDiagonal}, so we can talk about the second diagonal block  by projecting to the $Y$-components.  The following lemma describes the linear dynamics of the F4BP around the Lagrangian fixed point. Its proof is given in Proposition 8.2 of \cite{GHX}. 
 
 \begin{Lm}\label{LmDF4BPLoc}\begin{enumerate}
\item  The linearized dynamics at the lower Lagrange fixed point, when projected to the $Y$-component, has four eigenvalues $\mu_u>-\frac{v_*}{2}>0>\mu_s$, where $\mu_u$ and $\mu_s$ are that in Lemma \ref{LmDF3BPLoc}. Denote by $\mathbf e_1,\mathbf e_2,\mathbf e_3,\mathbf e_4\in \R^4$ the  eigenvectors  respectively, then we have that the projection of $\mathbf e_1$ and $\mathbf e_4$ to the $(w_0,\mathsf g_{01})$-plane are eigenvectors for the linearized dynamics in Lemma \ref{LmDF3BPLoc}. Moreover, the vector $\mathbf e_2$ has nontrivial projection to the $G_2$ component, so does its pushforward along $\gamma_O$ by the tangent dynamics. 
\item The tangent dynamics along $\gamma_I$ and $\gamma_O$ commutes with the projection to the $(w_0,\mathsf g_{01})$-plane $($the F3BP$)$. 
\end{enumerate}
 \end{Lm}
%The proofs of the last two lemmas can be found in Section 8 of \cite{GHX}. 

   With the last lemma, we can construct the two-dimensional plane $E^s_{4}(\gamma_I(0))\subset T_{\gamma_I(0)}\cS_-$ satisfying that its pushforward by the tangent flow along $\gamma_I(\tau)$ as $\tau\to\infty$ is $\mathrm{span}\{\mathbf e_3,\mathbf e_4\}$. Similarly, we construct the 2 dimensional plane $E^u_{4}(\gamma_O(0))\subset T_{\gamma_O(0)}\cS_+$ whose pushforward by the tangent flow along $\gamma_O(\tau)$ as $\tau\to-\infty$ is $\mathrm{span}\{\mathbf e_1,\mathbf e_2\}$. Moreover, by Lemma \ref{LmDF4BPLoc}, we have that the projection to the $(w_0,\mathsf g_{01})$-component  of the pushforward of $\mathbf e_1$ to $T_{\gamma_I(0)}\cS_-$(respectively  the pushforward of $\mathbf e_4$ to $T_{\gamma_O(0)}\cS_+$) is $E^s_3(\gamma_I(0))$ (respectively $E^u_3(\gamma_O(0))$).

After the application of the local map, we perform a renormalization on the section $\cS_+$. The definition of the renormalization is to extend the Definition \ref{renormalize-map} to the four-body problem.  The renormalization does not change the values of the variables $(w_0,\mathsf{g}_{01},w_2,\mathsf{g}_{21})$, but we need to use the variables $(G_0,\mathsf{g}_{01},G_2,\mathsf{g}_{21})$ in order to apply the global map.% By the definition of the renormalization, the $r$-variable is rescaled to a scalar of order $1/\eps$, so the value of $G_i$ is about $1/\eps^{1/2}$ times of that of $w_i,\ i=0,2$. 

%Next, we consider  the F3BP that is a small perturbation of the I3BP with nonisosceles variables $(G_0,g_0)$. Linearizing at the lower Lagrangian fixed point, we get two eigenvalues $\mu_u>0>\mu_s$. We define $E^s_\therefore(\gamma_I(0))$ and $E_\therefore^u(\gamma_O(0))$ similarly as above in the present setting. 
%We denote by $d(\widetilde{\cR \mathbb L})$ the projection/restriction of $d\cR\mathbb L$ to the  $(G_0,g_0)$-components in the F3BP when evaluated along an orbit of the I3BP. 

In the following propositions, without danger of confusions, when talking about the vectors $\bu_1,\bl_1$ that were defined in Proposition \ref{PropDL1} for the I4BP, we naturally embed them into the phase space of the F4BP by adding zeros to the extra dimensions. Similarly for the subspaces $E^{u/s}_{4}$, $E^{u/s}_{3}$, etc. Moreover, though the basepoint of the subspaces $E^{s}_{4}$, $E^{s}_{3}$ is $\gamma_I(0)$, we naturally parallel transport these subspaces to nearby points of $\gamma_I(0)$. Similarly for $E^{u}_{4}$, $E^{u}_{3}$. 
\begin{Prop}\label{PropDL2} 
For all $\eta>0$, there exists a small number $d_0>0$ such that the following holds: Let $\gamma:\ [0,T]\to \cM$ be an orbit of the F4BP in the blowup coordinates with $\gamma(0)\in \cS_-,\ \gamma(T)\in \cS_+$ and $\gamma(0)$ is sufficiently close to $\gamma_I(0)$ such that
$|Y|<d_0 e^{-\mu_u T},\  \chi>d_0^{-1}e^{\mu_u T}$ along $\gamma$, then we have $$d(\cR \mathbb L)=(e^{-T v_*} \bu_1\otimes \bl_1)\oplus d(\widehat{\cR \mathbb L})+O(d_0),$$
\begin{equation}\label{EqdRL}
d(\widehat{\cR \mathbb L})P\subsetneq \cC_\eta(E^u_{4}(\gamma_O(0))),\quad \forall\ \textrm{2-plane}\ P\  \mathrm{with}\ P\cap \cC_\eta(E^s_{4}(\gamma_I(0)))=\{0\}. \end{equation}
where 
$\bu_1$ and $\bl_1$ are as in Proposition \ref{PropDL1}. Moreover, for all $v\in P$, we have$$  \|d(\widehat{\cR \mathbb L})v\|\geq d_0e^{-Tv_*/2}\|v\|.$$
\end{Prop}
%We refer the reader to  \cite[Section 8]{GHX} for a proof of the above statements. 

\subsection{The global map}In this section, we study the derivative of the  global map. 

\begin{Prop} \label{PropDG2}  
For all $\eta>0$, there exist $\chi_0,\ \eps_0, \nu_0$ such that for all $\chi>\chi_0,\eps<\eps_0, \nu<\nu_0$,  we have the following: Let $\gamma:\ [0,T]\to\cM$ be an orbit of the F4BP with $\gamma(0)\in \cS_+$, $\gamma(T)\in 
\cS_-$ and $|Y|\leq \nu$ along the orbit.  Then there exist vector fields $\bar\bu_i,\ i=1,2,3,$ and $\bar\bl_i$,  $i=1,2$, such that we have 
$$d\mathbb G=f_{\mu,\lambda}\cdot\chi \bar \bu_1\otimes \bar\bl_1+\chi^2 \bar \bu_2\otimes \bar\bl_2+O((\eps+\nu+\lambda^{-1})f_{\mu,\lambda}\cdot\chi), $$
and 
\begin{equation}\label{EqDG}d\mathbb G \cC_{\eta}(\bu_1, E_{4}^u(\gamma_O(0)))\subsetneq\cC_{2\eta}(\bar\bu_1, \bar\bu_2,\bar \bu_3),\end{equation}
and for each vector $v\in \cC_{\eta}(\bu_1, E_{4}^u(\gamma_O(0)))$, we have
\begin{equation}\label{EqDG>}
\|d\mathbb G v\|\geq \frac12\|v\|,\end{equation}
where $\bu_1$ is as in Proposition \ref{PropDL1}, and $\bar\bu_1$ and $\bar\bl_1$ are as in Proposition \ref{PropDG1}. These vectors have definite limits in the isosceles limit $Y\to 0$ and $\chi\to\infty.$ Moreover, we have $dG_2 \cdot \bar \bu_2\neq 0$. 
\end{Prop}
We give a proof of this statement in  Appendix \ref{SG0g0}, which is slightly different from that of \cite{GHX}. The proof involves a careful analysis of the structure of the $o(f_{\mu,\lambda}\cdot\chi)$ error term, since  there is no expansion in the $\bar\bu_3$ direction (this is why we have \eqref{EqDG>}),  we have to show that the $o(f_{\lambda,\mu}\cdot\chi)$ error in $d\mathbb G$ does not spoil $\bar\bu_3$. 
All the above vectors can be given explicitly.  We shall be interested in controlling the plane $\mathrm{span}(\bar\bu_2,\bar\bu_3)$ being pushforward by the flow, of which, one direction is used to control the $(G_0,\mathsf{g}_{01})$, i.e. the motion of be binary and the other is used to control $(G_2,\mathsf{g}_{21})$, i.e. to make sure that the particle $Q_1$ comes to a near triple collision configuration with the binary $Q_3$-$Q_4$ after application of the global map.
\subsection{The transversality condition}
When we compose $\mathcal R\mathbb L$ and $\mathbb G$ to get the Poincar\'e map $\mathbb P=\mathbb G\mathcal R\mathbb L$, we hope that $d\mathbb P$ preserves the cone $\cC_{2\eta}(\bar\bu_1, \bar\bu_2,\bar \bu_3)$ and expands each vector inside. For this purpose, we have to verify that $\mathrm{span}(\bar\bu_2,\bar \bu_3)$ satisfies \eqref{EqdRL} for a plane $P$, i.e. $\mathrm{span}\{\bar\bu_2,\bar\bu_3\}\cap \cC_\eta(E_{4}^s(\gamma_I(0)))=\{0\}. $ 

Let $\mathbf e_3',\mathbf e_4'$ be two vectors in $E_{4}^s(\gamma_I(0))$ spanning unit area, then the transversality condition is given explicitly as the nonvanishing of the determinant $\det(\bar\bu_2,\bar\bu_3, \mathbf e_3',\mathbf e_4')\neq 0$. We have the following lemma.
\begin{Lm} \label{LmTrans} The determinant  $\det(\bar\bu_2,\bar\bu_3, \mathbf e_3',\mathbf e_4')$ is an analytic function
of the masses $m_1,m_2$, so it is nonvanishing for generic masses if it is so for one choice of $(m_1,m_2)$.  In particular, for such masses we have $\text{span}\{\bar{\bu}_2,\bar{\bu}_3\}\cap\mathcal{C}_\eta(E^s_{4}(\gamma_I(0)))=\{0\}$ for some $\eta>0$. 
\end{Lm}
%The lemma is similar to Lemma 4.7 of \cite{GHX}. However, since the model here is different from that of \cite{GHX}, the lemma has to be verified in the present setting. 

Note that for the I3BP, we have $G_0=0$, and double collision between the binary occurs when $\psi=\pm\pi/2$. Thus to avoid double collision between the binary, it is enough to make sure that $G_0$ and $\psi$ do not vanish simultaneously, which is guaranteed by the following proposition. 

\begin{Not}\label{NotEs}
\begin{enumerate}
\item Denote by $\mathcal{B}_{O}\subset \mathbb{R}$ the set of all the instant $t$ where the $w_0$ or $G_0$ component of $E^u_3(\gamma_O(t))$ vanishes. Similarly, we define $\mathcal{B}_{I}$ to be that of $E^s_3(\gamma_I(t))$.
\item The sets $\Psi_O$ and $\Psi_I$ are, respectively,  the set of all the instant $t$ when the $\psi$-component of $\gamma_O$ and $\gamma_I$  has absolute value $\pi/2$.
\end{enumerate}
\end{Not}
%The transversality conditions in the following two propositions are necessary when we compose local and global maps.
\begin{Prop}\label{PropTrans1} For the choices of masses $m_1=m_3=m_4=1$ the following transversality conditions are satisfied:
 The sets $\mathcal{B}_O$ and $\mathcal{B}_I$ have finitely many elements, independent of $\chi$, and moreover, there exists $c>0$, independent of $\chi$, such that
$$\min\{\mathrm{dist}(\mathcal{B}_O, \Psi_O),\; \;\mathrm{dist}(\mathcal{B}_I, \Psi_I)\}>c.$$
 \end{Prop}
We present a numeric verification of the above statements in Appendix \ref{app-numeric}. 
\subsection{Proof of Theorem \ref{ThmF4BP}}

In this section, we give the proof of Theorem \ref{ThmF4BP}. 

We fix $(m_1,m_2)$ and $\lambda$ as in Theorem \ref{ThmI4BP} and $\lambda_i$, $i=1,\dots$ as in \eqref{delta-i}. We shall construct a genuine super-hyperbolic orbits (or non-collision singularities)  close to those obtained there. 

{\bf Step 1.} {\it the choice of admissible cube.}

We next introduce the notion of {\it admissible cube.} Let us fix a small  number $\eta>0$ throughout the proof. 
 
\begin{Def}\label{DefAdm} We say that a three-dimensional cube  $\mathsf C_i$ on $\cS_-$ is \emph{admissible} for the $i$-th iteration, if 
\begin{enumerate}
\item  $|Y|<\bar\nu_i:=e^{\mu_s \sum_{k<i}T_k}$ for each point $(X,Y)$ on $\sC_i$;
\item  $T\mathsf C\subset \mathcal C_{\eta}(\bar\bu_1,\bar\bu_2,\bar\bu_3)$;%, so we can parametrize the admissible surface by $\ell_0$ and $G_2$ since the corresponding entries in the vectors $\bu_1$ and $\bu_2$ respectively are nonvanishing;
%\item The $\ell_0$-component contains a interval of length $\dt e^{-T_1}$ around the $\ell_0$ value of $\sx$, and the $G_2$-component contains an interval of length $\dt^2$ around 0.
\item When parametrized by the variables $\ell_0,G_0,G_2$, we have that for each $(G_0,G_2)$, the $\ell_0$-curve intersects $W^s_3$ and has length at least $\eps^2$. For each $(\ell_0,G_2)$, the curve on the $(G_0,\sf g_{01})$-plane has length between $\eta\bar\nu_i$ and $\bar\nu_i$, and its two endpoints are on the boundary of the cone $\mathcal C_\eta(E_3^s(\gamma_I(0)))$. Note that the curve does not intersect $(G_0,\mathsf{g}_{01})=(0,0)$ for $\eta$ small and  its tangent vectors satisfy the cone condition in \eqref{EqdRL}. For each $(\ell_0,G_0)$, the curve on the $(G_2,\mathsf g_{21})$-plane has $G_2$ coordinates in the interval $[-1,1]\eta\bar\nu_i.$
\end{enumerate}
\end{Def}
%{\color{red}The choice of $\bar \nu_i$ cannot be correct. It is too large. Indeed, we have $\bar\nu_i=\lambda_i^{-D}$. For superhyperbolic orbit $\lambda_i$ is fixed for all $i$. However, the size of $|Y|$ has to shrink. We have $|Y_i|\sim e^{\lambda_s(\sum T_i)}=e^{\lambda_sTi}$ decays exponentially by the contraction of the hyperbolic fixed point. So the assumption $(3)$ cannot be satisfied. Choosing the value of $\bar\nu_i$ is actually a delicate matter. In \cite{GHX} we introduced sojourn time. We might also need it here. }

%Item (3) is given by the transversality in Lemma \ref{LmTrans}.

{\bf Step 2.} {\it the local map and the point-deleting procedure. }

We pick an  admissible cube $\mathsf C_i$, then apply the local map followed by the renormalization. By Proposition \ref{PropDL2} and the transversality condition Lemma \ref{LmTrans}, we know that the plane span$\{\bar \bu_2,\bar \bu_3\}$ gets expanded by $d\mathbb L$ by at least $\min \{e^{\mu_u T_i},e^{-\frac{v_*}{2} T_i}\}=e^{-\frac{v_*}{2} T_i}$ since we have $\mu_u>-\frac{v_*}{2}$ (c.f. Lemma \ref{LmDF4BPLoc}).  Thus, the projection to the $Y$-components of the cube $\mathsf C_i$ is strongly expanded by at least $e^{-T_iv_*/2}$ with $T_i=O(\log\lambda_i)$. %Thus we have $\sC_{i+1}\subset \bP \mathsf C_i$. 

We next examine the dynamics of the $(w_0,\mathsf g_{01})$-components closely.  The Lagrange fixed point has two eigenvalues $\mu_u(>0>)\mu_s$ when restricted to the $(w_0,\mathsf g_{01})$-plane. By  the cone condition  in \eqref{EqdRL}  and by the $\lambda$-lemma, a curve on $\sf C_i$ with constant $\ell_0,G_2$ component will approach  $W^u_{3}(\gamma_O(0))$ exponentially fast (by $e^{T_i \mu_s}$) and get stretched by $e^{T_i \mu_u}$. We keep only the part of the curve that stays within the cone $\mathcal C_\eta(E^s_3(\gamma_I(\tau))),\ \tau\geq 0$ and   $\mathcal C_\eta(E^u_3(\gamma_O(\tau))),\ \tau\leq 0$,  along the orbit. Therefore the remaining segment has length estimated as $e^{\mu_s T_i}$.

%If we consider all the four $Y$-variables, the point deleting procedure in the $Y$-components gives us a segment of length estimated as $e^{\mu_s T_1}$ and its $\sf g_{01}$-component contains a neighborhood of 0. 

%So we introduce a point-deleting procedure on $\bP\sC_i$ to keep only those points in $\sC_{i+1}$. 
%Such a point-deleting procedure induces one on the preimage $\sC_0$ with the remaining set of points $\bP^{-i}\sC_i.$

So we introduce a point-deleting procedure on $\cR\bL \sC_i$ as follows
\begin{enumerate}
	\item In the $X$-component, we keep only a subsegment resulting from Proposition \ref{PropLocalMM}. 
	\item In the $(G_0,\mathsf{g}_{01})$-component, along the flow, we keep only the part that stays within the cone $\mathcal C_\eta(E^u_3(\gamma_O(\tau))),\ \tau<0$. 
	\item In the $G_2$-component, we keep only the part with $G_2\in [-1,1]\eta e^{\mu_sT_i}\bar \nu_i$. 
\end{enumerate}
After the point-deleting procedure we obtain a cube $\bar \sC_i$ on $\cS_+. $
%{\it 3. The $(G_0,\sf g_{02})$ dynamics. }
%{\bf Here I have to define the vectors $\bar\bu_2,\bar\bl_2$ explicitly. }

{\bf Step 3.} {\it the perturbation of the nonisoscelesness to the $X$-variables.}

We next show  that the perturbation of the nonisoscelesness $|Y|\neq 0$ to the $X$-variables ($X=(r,v,w,\psi)$ in blowup coordinates) in the I4BP is negligible. 

 For this purpose, we need the following notion of sojourn time to give an upper bound on the time defining the local map. 

\begin{Def}[Sojourn time]\label{DefSojourn}
	Let $\{T_n\}$ be a sequence satisfying $T_n>C_1^{-1}\log \dt_0^{-1}$, where $\delta_0$ is as in Proposition~\ref{PropLocalMM}. We say that $\{T_n\}$ is a sequence of sojourn times of type-$(c,C),\ c>0,C>0$ if it satisfies in addition
	$$T_{n+1}\leq c\sum_{i=1}^n T_i+C,\ \forall \ n\geq 1. $$
\end{Def}
In the following, we assume that the above sojourn time estimates are satisfied for each step of the local map. 
Let us denote by $X'=F(X)$ the equation of motion for the I4BP in the blowup coordinates   and by $\bar X'=F(\bar X)+\mathcal E$ the equation of the F4BP, where the error $\mathcal{E}$ is estimated as $O(Y^2)$ (see \cite[Appendix A]{GHX}). Then denoting $\dt X=X-\bar X$, we get $(\dt X)'=DF \dt X+\mathcal E$ and we estimate the divergence of the orbits $$|\dt X(\tau)|\leq e^{M\tau}\int_0^\tau e^{-Ms} \mathcal E(s)ds\leq \bar C e^{M\tau } (Y^2(0))$$ by Gronwall, where $M$ bounds $DF$ which is estimated by $\mu_u$.   By  the point deleting process, we see that $Y^2(0)$ is estimated by $e^{2\mu_s(\sum_{i=1}^{n-1} T_i)}$, inductively obtained from the last step if $n> 1$. By our assumption on the sojourn time $\tau\leq T_n$ we thus get an upper bound $|\dt X(\tau)|\leq \bar C e^{-\bar DT_n}$ for $\tau\in[0,T_n]$ and $\bar D=-\frac{Mc+\mu_s}{c}\gg1$ if we choose $c$ small.   For the $r$-variable, we note that the $r$-equation is $r'=rv$, which does not involve the $Y$-variables, thus the first row of $\mathcal E$ vanishes.  For the other variables variables $(v,w,\psi)$ of $X$, most of the time the orbit stays close to the Lagrangian fixed point thus the estimate of the fundamental solution of the equation $\dt X'=DF(X)\dt X$ is dominated by exponentiating $DF$ restricted at the Lagrangian fixed point for time $T_n$.  When $\bar D$ is much larger than $|\mu_u|$ and $|\mu_s|$, we see that the $O(e^{-\bar DT_n})$-perturbation to the $X$-variables is negligible.

{\bf Step 4.} {\it the global map and the iteration. }

After the point-deleting procedure during the local map, we obtain a cube $\bar{\mathsf C}_i$ on $\cS_+$ whose tangent lies in the cone $\mathcal C_\eta(\bu_1,E^u_{4}(\gamma_O(0)))$. %Moreover, since the $\mathsf g_{21}$-component  of $\bar\bu_2$ is nonvanishing, we see that the $\mathsf g_{21}$-component of the admissible surface contains a neighborhood of $\mathsf g_{21}=0$ that is independent of $\chi$. In the image of the global map, due to the strong expansion and the fact that $\bu_2$ has nonvanishing $G_2$ and $\mathsf g_{21}$ components, we get that in the image $\mathbb P(\mathsf C)$, we can choose a piece of admissible surface. 
We next apply the global map. We use $d\mathbb G$ in Proposition \ref{PropDG2} to control the difference between an orbit issued from $\bar{\mathsf C}_i$ and an orbit of the I4BP. The assumption $|Y|\leq \nu$ of Proposition \ref{PropDG2} is verified by the  following lemma whose proof could be found in \cite[Section 9.5]{GHX}.
\begin{Lm} \label{gG0-O1}There exist $C>0$, $\tilde C>0$ and $\nu_0>0$ independent of $\chi$ such that 
	\begin{enumerate}
		\item for the escaping piece of an orbit, if on $\cS_+$, 
		$$|G_0|=\nu_1\in(0,\nu_0],\;|\mathsf{g}_{01}|,\;|\mathsf{g}_{21}|,\; |G_2|\leqslant C|G_0|,$$ 
		then we have 
		$$G_0(t)\in[(1-\epsilon)\nu_1,(1+\epsilon)\nu_1]\text{ and }|\mathsf{g}_{01}(t)|,\;|\mathsf{g}_{21}(t)|,\; |G_2(t)|\leqslant \tilde C\nu_1;$$
		\item for the returning piece of an orbit,  if on $\cS_-$, $$|G_0|=\nu_2\in(0,\nu_0],\;|\mathsf{g}_{01}|,\;|\mathsf{g}_{21}|,\; |G_2|\leqslant C|G_0|,$$ 
		then we have 
		$$G_0(t)\in[(1-\epsilon)\nu_2,(1+\epsilon)\nu_2]\text{ and }|\mathsf{g}_{01}(t)|,\;|\mathsf{g}_{21}(t)|,\; |G_2(t)|\leqslant \tilde C\nu_2.$$
	\end{enumerate}
\end{Lm}

Note that the point-deleting procedure during each local map decreases the value of $|Y|$ by $e^{\mu_sT_i}. $ Thus to verify $|Y|\leq \nu$ for the global map with the last lemma, it is enough to choose $T_1$ large so that $e^{\mu_s T_1}<\nu$. This can be done by choosing $\chi_0$ and $\lambda_*$ large. 

The last lemma also implies that when projected to the $(G_0,\mathsf g_{01})$ plane, the orbit never intersect the piont $(0,0)$. Next,  estimate \eqref{EqDG>} in Proposition \ref{PropDG2} implies that $\mathbb G(\bar{\mathsf C}_i)$, when projected to the plane span$\{\bar\bu_2,\bar \bu_3\}$,  covers a disk of radius $\frac14e^{\mu_s T_i}$ centered at zero. We perform another point-deleting procedure so that we get an admissible cube  $\sC_{i+1}$ from $\mathbb P({\mathsf C}_i)$ according to item (3) of Definition \ref{DefAdm}. %We remark that though the global map has no expansion in the $\bar\bu_3$ direction,  the transversality in Lemma \ref{LmTrans} and the cone condition in  \eqref{EqdRL} show that  it will be expanded by the local map. 

Thus, we can iterate the procedure for infinitely many steps. After pulling back to the initial admissible cube $\mathsf C_i$ along the flow, the point-deleting procedure gives rise to a point-deleting procedure on the initial admissible  cube $\bP^{-i}(\sC_i)\subset \mathsf C_0$, hence after infinitely many steps,  we get a Cantor set as initial conditions.

{\bf Step 5.} {\it No binary collision. } 

For the local map, it may happen that $G_0$ vanishes at some isolated points whose number does not depend on $\chi$ when  near triple collision. However, when this happens, by Proposition \ref{PropTrans1}, we have $\psi\neq \pm\pi/2$, so the binary collision actually does not occur for the local map piece of orbit. 

We next consider the global map. First, by Lemma \ref{gG0-O1},  for the whole orbit defined by the global map at step $i$, we have $G_0\in [(1-2\epsilon)\tilde\nu,(1+2\epsilon)\tilde\nu]$ with $\tilde\nu\in[
	\eta\bar\nu_i,\bar\nu_i]$, since by the point-deleting procedure during the local map, when projected to  the $(G_0,w_{01})$, we keep only the part within the cone $\cC_\eta(E_3^s(\gamma_O(0)))$, which is almost parallel to $E^u_3(\gamma_O(0))$ by the $\lambda$-lemma. In particular, $G_0$ never vanishes, which implies that the double collision between the binary $Q_3$-$Q_4$ does not occur.

It remains to exclude the possibility of  binary collision between  $Q_1$ and $Q_2$. We invoke the following proposition whose proof can be found in Section \ref{avoid-dc}.
\begin{Prop}\label{lm-avoid-dc} 
There is no double collision between $Q_1$ and $Q_2$ for the global map piece of orbits returning to an admissible cube on $\cS_-$. 
\end{Prop}
 To summarize, we have excluded all the possible collisions for an orbit issued from the Cantor set. 
 
For any initial point in the Cantor set, the estimate of the time during each return  is the same as that in Theorem \ref{ThmI4BP}, thus we get a super-hyperbolic orbit or a non-collision singularities depending on the choice of sequence $\{\lambda_i\}$. We hence prove Theorem \ref{ThmF4BP}.
\qed
% Now we can repeat the above argument by replacing $Q_1$-$Q_3$-$Q_4$ by $Q_1$-$Q_3$-$Q_4$. The procedure can be iterated for infinitely many steps. In the limit, we will result in a Cantor set on $S_0$ whose points can be iterated for infinitely many steps. 

\section{Excluding double collisions}\label{avoid-dc} %We consider only the orbits as in Proposition \ref{choiceofm2} and divide the orbits into three pieces: the escaping piece traveling  from $\cS_+$ to $\cS^M$, the intermediate piece from  $\cS^M$ to $\cS^M$, and the returning piece from $\cS^M$ to $\cS_-$. 
We are now in position to prove Proposition \ref{lm-avoid-dc} on excluding the double collision between $Q_1$ and $Q_2$. 
\subsection{Proof of Proposition \ref{lm-avoid-dc}}
\begin{proof}[Proof of Proposition \ref{lm-avoid-dc}] 
Suppose a piece of orbit $\gamma:\ [0,T]\to \cM$ defining the global map is issued from the section $\cS_+$ and returns to $\cS_-$ and experiences a double collision between $Q_1$ and $Q_2$ at time $T'$.

Let us suppose that after the near triple collision on the section $\cS_+$, we have $G_1^L=\bar\nu$ that is almost the number $\bar \nu_{n}$ after the $n$-th step of renormalization by the above point-deleting procedure. 

In the following, we shall  show that 

{\bf Claim: }{\it Suppose there is a double collision between $Q_1$ and $Q_2$, then there is a constant $c$ independent of $\chi$ and $\bar\nu$, such that we have $|G_2^L(T)|\geq c\bar\nu$ when the orbit comes back to the section $\cS_-.$}

Note that the conclusion of this claim does not meet the item (3) in our definition of admissible cube (c.f. Definition \ref{DefAdm}) if we choose $\eta$ sufficiently small therein. Thus our choice of admissible cubes excludes the possibility of collisions between $Q_1$ and $Q_2$ for any orbit issued from the Cantor set constructed in the last section. 

%In order for $Q_1$ to come to the next near triple collision with the binary, we need that the total angular momentum of the triple $-G_2^L=G^L_0+G^L_1=O(\nu)$. On the other hand, we shall prove below that the total angular momentum  $G_0^L+G_1^L$ is   order $O(\chi \nu)$, if double collision between $Q_1$ and $Q_2$ occurs. This gives a contradiction. Thus, no double collision between $Q_1$ and $Q_2$ is possible if we choose initial condition in the Cantor set that we constructed in Theorem \ref{ThmF4BP}.
 %We first outline the idea of the proof. Indeed, suppose when $Q_1$ is of $O(\chi)$-distance away from $Q_2$ and the binary, then the angular momentum $G_1^L=O(\nu)$  and $G_2^L=O(1/\chi)$. This means that $\angle(\bx_1^L,\bp^L_1)=O(1/\chi)$ and $\angle(\bx^L_2,\bp^L_2)=O(1/\chi^2)$.  Note that $|\bp_2^L|\ll|\bp_1^L|$ after the renormalization. Modeling the binary collision between $Q_1$ and $Q_2$ as an elastic collision, we get that after collision $Q_2$ is going to gain a velocity almost proportional to $\bp_1^L$ before collision. In particular, we get after collision $\angle(\bx^L_2,\bp^L_2)=O(1/\chi)$, which implies tha that the angular momentum $G_2^L$ after collision becomes $O(1)$ as claimed. 

We next give the proof of the claim.  The heuristic idea is that a collision between $Q_1$ and $Q_2$ will transfer a big amount of angular momentum of $G_1^L$ to $G_2^L$, so we can bound $G_2^L$ from below.  It is important to introduce a new section $\cS_d:=\{|x_1^R|=d\chi\}$ for some sufficiently small $d$ to be determined below, which is a section $d\chi$-distance away from the double collision of $Q_1$-$Q_2$. 

The section $\cS_d$ break the piece of orbit $\gamma$ into three pieces:
\begin{enumerate}
	\item The segment from $\cS_+$ to first intersecting $\cS_d$;
	\item The segment from first intersecting $\cS_d$ to second intersecting $\cS_d$;
	\item The segment from leaving $\cS_d$ to $\cS_-$. 
\end{enumerate}

For the first and third pieces, we shall work in the left Jacobi coordinates and show that $G_1^L$ and $G_2^L$ do not oscillate much, but for the second segment, we shall work in the right Jacobi coordinates and show that all the Delaunay coordinates do not  oscillate much, but when converted to the left Jacobi coordinates, there would be a big amount of angular momentum transferred to $G_2^L$.

For the piece (1), we apply Lemma \ref{gG0-O1} to get that when arriving at the section $\cS_d$, we have the estimate $|G_0|\geq \frac{9}{10}\bar\nu.$

We next consider piece (2).  It would be  convenient to introduce the angle $\al:=\angle(\bx_1^R,\bx_2^R)$ (see Figure \ref{col-f}) measured at the time when the orbit first cross the section $\cS_d$.%=(1+O(d))\angle(\bx_1^L-\bx^L_2,\bx_2^L)$ where the second $=$ follows from the definition of the right Jacobi coordinates. 
\begin{figure}[htb]
\includegraphics[scale=0.3]{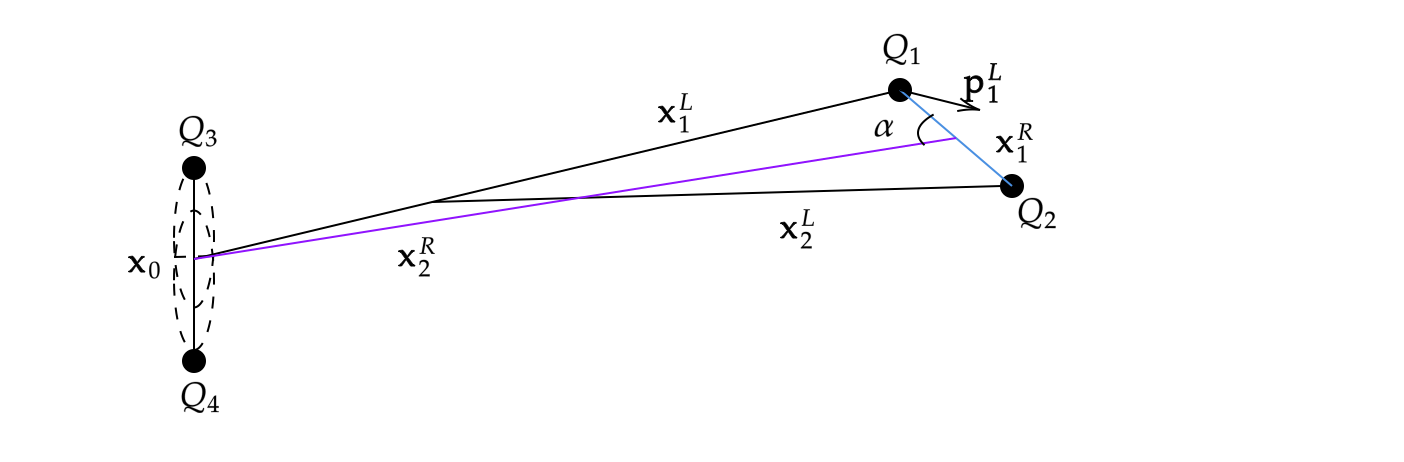}
\caption{The collision configuration}
\label{col-f}
\end{figure}
We shall use the following lemma whose proof is postponed to the next subsection. 
\begin{Lm}\label{LmAlpha}
	The angle $\al:=\angle(\bx_1^R,\bx_2^R)$ is estimated as $C^{-1}\frac{\bar\nu}{\chi}\leq |\al|\leq C\frac{\bar\nu}{\chi}$. 
\end{Lm}
%We choose a small number $d>0$ to be determined and use the right Jacobi coordinates in the region $\{|\bx_1^R|\leq d\chi\}$. 

Collision occurs if $\bx_1^R=0$. We let $\gamma_{in}$ be the piece of orbit coming to the collision with time reversed and let $\gamma_{out}$ be the piece of orbit ejected from the collision. After a time translation, we let $T_{in}$ (respectively $T_{out}$) be the time such that $\gamma_{in}$ (respectively $\gamma_{out}$) intersects the section $\cS_d$ and $0$ be the collision time. We shall compare the orbit parameters of the two orbits such as $\dt G^R_1=G^R_{1}(T_{out})-G^R_{1}(T_{in})$ and similarly we apply the $\dt$-notation to other variables. Note that at the collision time $t=0$, we have $\dt G^R_1(0)=0,\ \dt g^R_1(0)=0.$

We shall integrate the Hamiltonian equations derived from Hamiltonian \eqref{3pairHam} over time $d\chi$. Let us consider the $\dot G_1^R$ equation first. 
We have $\dot G_1^R=-\frac{\partial H}{\partial g^R_1}=-\frac{\partial U_{12}}{\partial g^R_1}$ where the potential $U_{12}$ is in \eqref{3pairHam}. Note that  by assumption we have $|\bx_0^R|=O(1),\ |\bx_1^R|=O(d\chi)$, $|\bx_2^R|\sim \chi$ and $\al=\angle(\bx_1^R,\bx_2^R)=O(\frac{\bar\nu}{\chi})$ by Lemma \ref{LmAlpha}. Thus we estimate 
\begin{equation}\label{EqG1R}
|\dot G_1^R|\leq  C\frac{|\bx_2^R||\bx_1^R|\al }{|\bx_2^R|^3}=C\frac{d\bar\nu}{\chi^2}.
\end{equation}

For both $\gamma_{in}$ and $\gamma_{out}$, we shall integrate the equation over time of order $d\chi$ to go from the boundary of $\{|\bx_1^R|\leq d\chi\}$to the double collision and back. This gives the estimate $G_{1}^R(T_{in}), G_{1}^R(T_{out})=O(\frac{d^2\bar\nu}{\chi})$ since we have $G_1^R=0$ at collision. 

By similar argument,  we get 
\begin{Lm}\label{LmInBall}
	Suppose there is a collision between $Q_1$ and $Q_2$. Then we have the quantities 
	$\dt G_1^R, \dt G^R_{2}, \dt g_1^R$ and  $\dt g_2^R$  are all estimated by  $O(\frac{d^2\bar\nu}{\chi})$. 
\end{Lm}

 By definition of $\al$, we have $(g_1-g_2)^R(T_{in})=\al$, thus we have $(g_1-g_2)^R(T_{out})=(1+O(d))\al$.

On the section $\cS_d$, we convert to the left Jacobi coordinates. From \eqref{R-L-transit}, we get
\begin{equation}
	\begin{aligned}
		G_2^L&= x_2^L\times p_2^L=\left(-\frac{m_1(2+m_1+m_2)}{(2+m_1)(m_1+m_2)}\bx_1^R+\frac{2}{2+m_1}\bx^R_2\right)\times \left(-\bp_1^R+\frac{m_2}{m_1+m_2}\bp_2^R\right)\\
		&=\frac{m_1(2+m_1+m_2)}{(2+m_1)(m_1+m_2)}G_1^R+\frac{2m_2}{(2+m_1)(m_1+m_2)}G_2^R\\
		&-\frac{2}{2+m_1}\bx^R_2\times \bp_1^R-\frac{m_1m_2(2+m_1+m_2)}{(2+m_1)(m_1+m_2)^2} \bx_1^R\times \bp_2^R\\
		&:=c_1G_1^R+c_2G_2^R-c_3\bx^R_2\times \bp_1^R-c_4\bx_1^R\times \bp_2^R,
	\end{aligned}
\end{equation}
where $c_i,\ i=1,2,3,4$ are positive numbers whose values are not important for us. 
This gives
\begin{equation}\label{EqdG}
\dt G^L_2=c_1\dt G_1^R+c_2\dt G_2^R-c_3\dt (\bx^R_2\times \bp_1^R)-c_4\dt(\bx_1^R\times \bp_2^R).	
\end{equation}

The first two terms on the RHS are estimated in the last lemma. We shall analyze the two latter terms in more details. It turns out that the main contribution is given by the term $\dt (\bx^R_2\times \bp_1^R)$. We have the following lemma, whose proof is postponed to the next subsection. 
\begin{Lm}\label{LmDeltaG} Restricted to the section $\cS_d$, 
we have $|\dt G^L_2|\geq c \bar\nu$ for some constant $c>0$ that is  independent of $\bar\nu$ or $\chi$. 
\end{Lm}

Finally, we consider piece (3). 
After leaving the section $\cS_d$, we shall use the left Jacobi coordinates. We next consider the evolution of $G_2^L$ derived from the Hamiltonian \eqref{EqHamGlob}. From the Hamiltonian equation, we get $\dot G_2^L=-\frac{\partial H}{\partial g_2^L}=-\frac{\partial U_2}{\partial g_2^L}$. To estimate the RHS of the equation, we note that $$\bx_0^L=O(1), \ \bx_1^L=O(\chi),\ \bx_2^L=O(\chi),\quad \angle(\bx_1^L,\bx_2^L)=O(1/\chi). $$  This gives the estimate 
$|\dot G_2^L|\leq C\frac{|\bx_1^L||\bx_2^L|1/\chi}{|\bx_2^L-\bx_1^L|^3}=C\frac{\chi}{d^3\chi^3}$. Integrating over time of order $\mu^{-1/2}\lambda_*^{1/2}\chi$ by Proposition \ref{choiceofm2} using the fact that $|Q_1(t)-Q_2(t)|$ like $\mu^{1/2} \lambda_*^{-1/2}t$,  we get 
\begin{equation}\label{EqdeltaG2L}
G_2^L(T')-G_2^L(T_{out})=O(\lambda_*^{1/2}/(\mu^{1/2}d^3\chi))	
\end{equation}
 where $T'$ is a final time when $Q_1$ is close to the binary. To proceed, we need the following Sublemma whose proof is postponed to the next subsection

\begin{SLm}\label{Lmlambda}For $m_1$ in a neighborhood of 1, we have the estimate $-\mu_s<-v_*$ for the I3BP.  	
\end{SLm}
\begin{proof}
These quantities can be explicitly computed using the formulas from \cite[Lemma D.1]{GHX}. For $m_1=1$, direct calculation gives $\mu_s=- 0.7848$ and $v_*=-\sqrt{6}$, which  verify the assertion. 
\end{proof}
Note that $\chi=\chi_0\prod_{i<n}\lambda_i=\chi_0e^{-v_*\sum_{i<n}T_i}$ after the $n$-th step of renormalization and $\bar\nu=e^{\lambda_s\sum_{i<n}T_i}$.  Thus, we get by the Sublemma that $1/\chi\ll \bar\nu$, hence for large $\chi_0$,  the estimate \eqref{EqdeltaG2L} implies that $G_2^L(T')$ remains of order $\bar\nu$ when $Q_1$ gets close to the binary. This completes the proof of the claim.  
%Assume the opposite, then there exists $t_0\in(0,T)$ such that $G_1^R(t_0)=0$ and $\mathsf{g}^R_{21}(t_0)=0$. We claim that on the $\cS^M$, 
%\begin{equation}\label{abound-Gg}G_1^R(T
\end{proof}

\subsection{Proof of the lemmas}
It remains to prove the lemmas used above. 

\begin{proof}[Proof of Lemma \ref{LmAlpha}]
	
%	Then by Lemma \ref{gG0-O1}, we get that $G_1^L$ remains almost $\bar\nu$ when arriving at the section $\cS^M$ as well as reaching the boundary of the region $\{|\bx_1^R|=d\chi\}$. 
	
	First, from the triangle formed by $\bx_1^L$ and $\bx_2^R$ in the figure, we know that the angle $\al\sim \frac{1}{d}\angle(\bx_1^L,\bx_2^R)$. Note also that the two angles $\al$ and $\angle(\bx_1^L,\bx_2^R)$ add up to $\angle(\bx_1^L, \bx_1^R)$, since the latter is the outer angle of the former two. Thus for small $d$, the angle $\al$ is almost the same as $\angle(\bx_1^L, \bx_1^R)$. 

We next estimate $\angle(\bx_1^L, \bx_1^R)$. By triangle inequality, we have 
\begin{equation}\label{Eqtriangle}
|\angle(\bx_1^L, \bp_1^L)|-|\angle(\bx_1^R, \bp_1^L)|\leq |\angle(\bx_1^L, \bx_1^R)|\leq |\angle(\bx_1^L, \bp_1^L)|+|\angle(\bx_1^R, \bp_1^L)|.	
\end{equation}
	In the following we shall show that $|\angle(\bx_1^R, \bp_1^L)|\ll |\angle(\bx_1^L, \bp_1^L)|$, thus $ \al=|\angle(\bx_1^L, \bx_1^R)|$ is almost $|\angle(\bx_1^L, \bp_1^L)|$, which is estimated as $C\bar\nu/\chi$ using the angular momentum $G^L_1\sim \bar\nu$. 
	
For the term $|\angle(\bx_1^R, \bp_1^L)|$ on the RHS of \eqref{Eqtriangle}, we further split 
\begin{equation}\label{Eqtriangle1}
|\angle(\bx_1^R, \bp_1^L)|\leq \angle(\bx_1^R, m\bp_1^R)|+|\angle(\bx_1^L, -m\bp_1^R+\bp_1^L)|,	
\end{equation}
	where the constant $m$ is chosen according to \eqref{L-R-transit} such that $-m\bp_1^R+\bp_1^L$ is a multiple of $\bp_2^L$, thus $$|\angle(\bx_1^L, -m\bp_1^R-\bp_1^L)|=|\angle(\bx_1^L, \bp_2^L)|\leq |\angle(\bx_1^L, \bx_2^L)|+|\angle(\bx_2^L, \bp_2^L)|.$$ 
	From Figure \ref{col-f}, we have  $|\angle(\bx_1^L, \bx_2^L)|\sim d|\angle(\bx_1^R, \bx_2^L)|<d\al$, which is much smaller than $\al$. Next,  to estimate $|\angle(\bx_2^L, \bp_2^L)|$, we need an estimate of $G_2^L$ on the section $\cS_d$. We shall show below that $G_2^L=O(1/(d^3\chi))$, which combining $\bx_2^L=O(\chi),\ \bp_2^L=O(\lambda^{-1/2})$ gives $|\angle(\bx_2^L, \bp_2^L)|\leq C  \sqrt\lambda/(d^3\chi^2)\ll\al$. It remains to estimate $G_2^L$.  We have the same  estimate as \eqref{EqdeltaG2L} for the oscillations of $G_2^L$ from the section $\cS_+$ to $\cS_d$. Thus to get that $G_2^L=O(1/(d^3\chi))$ on the section $\cS_d$, it is enough to prove a similar bound on $\cS_+$. This follows from the $\chi^2$-expansion in the estimate of $d\bG$ in Proposition \ref{PropDG2} and $dG_2\cdot \bar \bu_2$ therein. Indeed, since by the construction of the admissible cube in the last section, we require $G_2^L$ to be of order $\bar \nu$ on the section $\cS_-$ after application of $\bP$, thus its value on the section $\cS_+$ before applying $\bP$ is of order $\bar \nu/\chi^2$ by the $\chi^2$-expansion of $d\bG$. Thus we have proved the estimate $G_2^L=O(1/(d^3\chi))$ hence complete the estimate $|\angle(\bx_1^L, -m\bp_1^R-\bp_1^L)|\ll\al. $

The angle $\angle(\bx_1^R, m\bp_1^R)|=|\angle(\bx_1^R, \bp_1^R)|$ on the RHS of \eqref{Eqtriangle1} can be estimated using the angular momentum $G_1^R$ and the collision assumption. At the collision time we have $G_1^R=0$. Next assuming $\al<C\nu/\chi$, we have the estimate $\dot G_1^R=O(d\bar\nu/\chi^2)$ by \eqref{EqG1R}, which implies $G_1^R$ on the boundary of $\cS_d$ is $O(d^2\bar\nu/\chi)$ by integrating over time $O(d\chi)$. Since $G_1^R=\bx_1^R\times \bp_1^R$, we get $|\angle(\bx_1^R, \bp_1^R)|<d\bar\nu/\chi^2\ll$ since $|\bx_1^R|=d\chi$ and $\bp_1^R=O(1)$. Since the estimate of $|\angle(\bx_1^R, \bp_1^R)|$ is much less than the $C\bar\nu/\chi$ estimate of $\al$, thus \eqref{EqG1R} is always valid and the above argument is consistent. 

Thus we have proved that the term $|\angle(\bx_1^R, \bp_1^L)|$ much smaller than $\al,$ hence completes the proof. 

%In order to have a double collision between $Q_1$ and $Q_2$, then we have $\angle (\bp_1^L,\bx_2^L)\gg \angle(\bx_1^L,\bx_2^L)$. This implies $\nu\ll\bar\nu$.

\end{proof}

\begin{proof}[Proof of Lemma \ref{LmDeltaG}]
	To estimate $\dt G^L_2$, we analyze the RHS of \eqref{EqdG} term by term. We first have $\dt G_1^R,\dt G_2^R=O(d^2\bar\nu/\chi)$ by Lemma \ref{LmInBall}. 
	
%	For the rest of the terms, we first invoke \eqref{L-R-transit} to relate the right Jacobi coordinates to the left ones. Arguing as in the proof of Lemma \ref{LmAlpha}, we have both $\bp_2^R$ (the momentum of the mass center of $Q_1$ and $Q_2$) and $\bp_1^R$ (the momentum of the relative motion of $Q_1$ and $Q_2$) are given by multiples of  $\bp_1^L$ up to an error of order $\sqrt\lambda^{-1}$ in the direction of $\bp_2^L$. 
	
	For the remaining two terms, we next consider the term 	$\dt (\bx^R_2\times \bp_1^R)=O(\bar\nu)$ bounded away from zero and independent of $\chi$, which is the main contribution to $\dt G_2^L$. Indeed, we have that  $\bp_1^R$ has reversed the direction by $\pi +O(d^2\bar\nu/\chi)$, $|x_2^R|\sim \chi$ by the estimate of $\dt g_1^R$ in Lemma \ref{LmDeltaG} and the direction of $\bx_2^R$ has changed by $O(d^2\bar\nu/\chi)$ by the estimate of $\dt g_2^R$ in Lemma \ref{LmDeltaG}. Thus we have 
	$$\dt(\bx_2^R\times \bp_1^R)=(\bx_2^R\times \bp_1^R)(T_{out})-(\bx_2^R\times \bp_1^R)(T_{in})=-2(\bx_2^R\times \bp_1^R)(T_{in})+O(d^2\bar\nu).$$
		Moreover,  we have the estimate of the angle $$|\angle (\bx_2^R, \bp_1^R)|\geq |\angle (\bx_2^R, \bx_1^R)|-|\angle (\bx_1^R, \bp_1^R)|\geq \al-O(d^2\bar\nu/\chi)\geq c\frac{\bar\nu}{\chi}$$ for some constant $c>0$, where the $O(d^2\bar\nu/\chi)$ estimate of $\angle (\bx_1^R, \bp_1^R)$ comes from the estimate $G_1^R=O(d^2\bar\nu/\chi)$ preceding Lemma \ref{LmInBall} due to the collision. Thus we prove the estimate $\dt (\bx^R_2\times \bp_1^R)\geq c\bar\nu.$ %Note that we have treated $\bar\nu$ as a number of order  $1$ much larger than $d$. 
	
	Arguing similarly, we get $\dt(\bx_1^R\times \bp_2^R)=O(d\bar \nu)$. The estimate is different from $\dt (\bx^R_2\times \bp_1^R)=O(\bar\nu)$ is because $|\bx_1^R|=O(d\chi)$ and $|\bx_2^R|=O(\chi)$ for both the initial and final moments.

	This completes the proof. 
\end{proof}

\appendix
\section{Delaunay coordinates}\label{App-delaunay}
In this appendix, we give the relations between Cartesian coordinates and Delaunay coordinates. Consider the Kepler problem $H(P,Q)=\frac{|P|^2}{2m}-\frac{k}{|Q|},$ $(P,Q)\in \R^2\times \R^2$.
\subsection{Elliptic motion}
We have the following relations which explain the physical and geometrical meaning of the Delaunay coordinates.
$a=\frac{L^2}{mk}$ is the semimajor axis, $b=\frac{LG}{mk}$ is the semiminor axis, $E=-\frac{k}{2a}$ is the energy, $G=Q\times P$ is the angular momentum, and $e=\sqrt{1-\left(\frac{G}{L}\right)^2}$ is the eccentricity. Moreover, $g$ is the argument of apapsis and $\ell$ is the mean anomaly.  We can relate $\ell$ to the polar angle $\psi$ through
the equations
$\tan\frac \psi 2 = \sqrt{\frac{1+e}{1-e}}\cdot\tan\frac u 2,\quad u-e\sin u=\ell.$
We also have Kepler's law $\frac{a^3}{T^2}=\frac{k/m}{(2\pi)^2}$ which relates the semimajor axis
$a$ and the period $T$ of the elliptic motion.

Denoting the body's position by $Q=(q_1, q_2)$ and its momentum by $P=(p_1,p_2)$ we have the following formulas
in the case $g=0$
\begin{equation}
\label{DelEll}
\begin{aligned}
q_1=\frac{L^2}{mk}\left(\cos u-\sqrt{1-\frac{G^2}{L^2}}\right), \quad &
q_2=\frac{LG}{mk} \sin u,
\\
p_1=-\frac{mk}{L}\frac{\sin u}{1-\sqrt{1-\frac{G^2}{L^2}} \cos u}, \quad &
p_2=\frac{mk}{L^2}\frac{G\cos u}{1-\sqrt{1-\frac{G^2}{L^2}}\cos u},
\end{aligned}
\end{equation}
where $u$ and $\ell$ are related by $u-e\sin u=\ell$. This is an ellipse whose major axis lies on the horizontal axis, with one focus at the origin and the other focus on the negative horizontal axis. The body moves counterclockwise if $G>0$. 
%
%\[e_{G}=\frac{-\frac{G}{L^2}}{\sqrt{1-(\frac{G}{L})^2}}=\frac{-\frac{G}{L^2}}{e}, e_{GG}=-\frac{1}{L^2e}+\frac{\frac{G}{L^2}}{e^{2}}e_G.\]
%So we have $e_G|_{G=0}=0$, $e_{GG}|_{G=0}=-\frac{1}{L^2}$.
%\[u_{G}-e_{G}\sin u-e\cos u u_{G}=0.\]
%\[u_{GG}-e_{GG}\sin u-e_{G}\cos u u_G-e_{G}\cos u u_{G}+e\sin u (u_G)^2-e\cos u u_{GG}=0.\]
%So we have $u_G|_{G=0}=0$, $(u_{GG}+\frac{1}{L^2}\sin u- u_{GG}\cos u)|_{G=0}=0$.
%
%We have the partial derivative calculations
%\begin{equation}
%\begin{aligned}
%u_G(1-e\cos u)&=e_G\sin u,\quad u_G|_{G=0}=e_G|_{G=0}=0,\\
%u_{GG}(1-e\cos u)|_{G=0}&=e_{GG}\sin u|_{G=0},\quad
%e_{GG}|_{G=0}=-\frac{1}{L^2}.
%\end{aligned}
%\end{equation}
\subsection{Hyperbolic motion}
For Kepler hyperbolic motion we have similar expressions for the geometric and physical quantities
\[a=\frac{L^2}{mk}, \ b=\frac{LG}{mk},\ E=\frac{k}{2a},\ G=Q\times P,\ e=\sqrt{1+\left(\frac{G}{L}\right)^2}\]
as well as the parametrization of the orbit
\begin{equation}
\begin{aligned}
q_1=\frac{L^2}{mk}
\left(\cosh u-\sqrt{1+\frac{G^2}{L^2}}\right), \quad
& q_2=\frac{LG}{mk} \sinh u, \\
p_1=-\frac{mk}{L}\frac{\sinh u}{1-\sqrt{1+\frac{G^2}{L^2}} \cosh u}, \quad
& p_2=-\frac{mk}{L^2}\frac{G\cosh u}{1-\sqrt{1+\frac{G^2}{L^2}}\cosh u}.
\end{aligned}\label{eq: delaunay4}
\end{equation}
where $u$ and $\ell$ are related by
\begin{equation}u-e\sinh u=\ell, \text{ and } e=\sqrt{1+\left(\frac{G}{L}\right)^2}.
\label{eq: hypul}
\end{equation}
%Differentiating this relation, we have
%\[u_{G}-e_{G}\sinh u-e\cosh u u_{G}=0, u_{GG}-e_{GG}\sinh u-e_{G}\sinh u u_{G}-e_{G}\cosh u u_G-e\sinh u (u_G)^2-e\cosh u u_{GG}=0.\]
%So we have
%\begin{equation}
%e_G|_{G=0}=0, e_{GG}=\frac{1}{L^2},\
%u_G(1-\cosh u)|_{G=0}=0,\quad (u_{GG}-\frac{1}{L^2}\sinh u-u_{GG}\cosh u)|_{G=0}=0.\end{equation}
This hyperbola is symmetric with respect to the $x$-axis, open to the right, and the particle moves counterclockwise on it when $u$ increases ($l$ decreases) in the case when minus the angular momentum $G=Q\times P>0$.
The angle $g$ is defined to be the angle measured from the positive $x$-axis to the symmetric axis. There are two such angles that differ by $\pi$, depending on the orientation of the symmetric axis. %This $\pi$ difference disappears in  the hamiltonian equation after taking derivative so it does not matter which angle we choose. Here $g$ does not appear in \eqref{DelEll} and \eqref{eq: delaunay4} because the argument of apoapsis or perigee is chosen to be zero or $\pi$. 
In the general case, we need to rotate the $(q_1,q_2)$ and $(p_1,p_2)$ using the matrix $\mathrm{Rot}(g)=\left[\begin{array}{cc}\cos g & -\sin g \\\sin g & \cos g\end{array}\right]$, that is, we have the following formula for $\bx=(q_1,q_2)$ and $\bp=(p_1,p_2)$, 
\begin{equation}\label{EqCarDel}
\begin{aligned}
\bx&=\frac{1}{mk}\Big(-\cos(g)L^2(\cosh u-e)+\sin(g)LG\sinh u,-\sin(g)L^2(\cosh u-e)-\cos(g)LG \sinh u\Big),\\
\bp&=\frac{mk}{1-e\cosh u}\Big(\frac{1}{L}(\sinh u)\cos g-\frac{G}{L^2}(\sin g)\cosh u,\frac{1}{L}(\sinh u)\sin g+\frac{G}{L^2}\cos(g)\cosh u\Big).	
\end{aligned}	
\end{equation}
%Direct calculation gives us the following formula for the  product $\bp_1\times \bx_2$,
%\begin{equation}\label{crossproduct-1}
%\begin{split}\bp_1\times\bx_2
%&=\frac{m_1k_1}{(1-e_1\cosh u_1)m_2k_2}\Big(\frac{1}{L_1}\sinh u_1L_2^2(\cosh u_2-e_2)\sin(g_1-g_2)\\
%&\quad\quad\quad\quad+\frac{1}{L_1}\sinh u_1L_2G_2\sinh u_2\big(-\cos(g_1-g_2)\big)\\
%&\qquad\qquad +\frac{G_1}{L_1^2}\cosh u_1L_2^2(\cosh u_2-e_2)\big(\cos(g_1-g_2)\big)\\
%&\qquad \qquad+\frac{G_1}{L_1^2}\cosh u_1L_2G_2\sinh u_2\big(\sin(g_1-g_2)\big)\Big).
%\end{split}\end{equation}
%Switching the subindex 1 and 2, we have a formula of the cross product $\bp_2\times \bx_1$. 
\section{The derivative of the global map in the I4BP}\label{L-Derivatives}
In this section, we give the proof of Proposition \ref{PropDG1}. The proof is very similar to Proposition 3.11 of \cite{GHX}. with the  difference as follows. The estimate in \cite{GHX} consists of two segments: from the left $\cS_+$ to the middle section $\cS^M$ (where we exchange the left and right Jacobi coordinates) and from  $\cS^M$ to the right $\cS_-$ (the right to the left case is similar). However, the estimate here consists of three segments: from $\cS_+$ to $\cS^M$, from $\cS^M$ to the double collision of $Q_1$ and $Q_2$ and back to $\cS^M$, and from  $\cS^M$ to the next $\cS^+$.  It looks the situation here is slightly complicated, but the proof and the result are completely analogous to that of \cite{GHX}. Indeed, all the derivative matrices composing $d\bL$ was already computed in \cite{GHX}, only the way of composing them is different. In the following, we give the proof and refer readers to \cite{GHX} for some detailed calculations. 

We first recall a formula from   \cite[Section 7]{X1} computing the derivative of the Poincar\'e map defined by cutting the flow $\mathcal V'=\mathcal F(\mathcal V,t)$
 between the sections $S^i$ and $S^f.$ We use $\mathcal V^i$ to denote the values of variables $\mathcal V$ restricted on the {\it initial} section $S^i$, while $\mathcal V^f$ means values of $\mathcal V$ on the {\it final} section $S^f$, and we use $t^i$ to denote the initial time and $t^f$   the final time.  We want to compute the derivative
$\cD:=\frac{D \mathcal V^f}{D\mathcal V^i}$ of the Poincar\'{e} map.  We have
\begin{equation}
\cD
=\left(\Id-\mathcal F(t^f)\otimes \frac{D t^f}{D\mathcal V^f}\right)^{-1}\frac{D\mathcal V(t^f)}{D\mathcal V(t^i)}\left(\Id-\mathcal F(t^i)\otimes \frac{D t^i}{D \mathcal V^i}\right).
\label{eq: formald4}
\end{equation}
Here the middle term $\frac{D\mathcal V(t^f)}{D\mathcal V(t^i)}$ is the fundamental solution to the variational equation $(\dt \mathcal V)'=\nabla_{\mathcal V}\mathcal F(\mathcal V, t)\dt \mathcal V$ and the two terms $\left(\Id-\mathcal F(t^f)\otimes \frac{D t^f}{D\mathcal V^f}\right)^{-1}$ and $\left(\Id-\mathcal F(t^i)\otimes \frac{D t^i}{D \mathcal V^i}\right)$ are called boundary contributions taking into account the issue that different orbits take different times to travel between the two sections. 

We  only consider orbits from Proposition \ref{choiceofm2}.  Recall that we introduce a middle section $\cS^M=\{|\bx_1^R|=\beta\chi\}$ and  divide the orbit into three pieces, the escaping piece (from $\cS_+$ to $\cS^M$), the intermediate piece (from $\cS^M$ to $\cS^M$) and the returning piece (from $\cS^M$ to $\cS_-$).  
For the escaping and returning pieces of the orbit, we use the left Jacobi coordinates and for the intermediate piece, we use the right Jacobi coordinates.  For all the three pieces, we shall take $\ell_1$ to be the new time and eliminate the variable $L_1$ using the energy conservation. Thus the remaining $X$-variables are $(L_0,\ell_0,L_2,\ell_2).$ The proof of the Proposition \ref{PropDG1} consists the following three steps. 
\begin{itemize}
\item Step 1: estimating the fundamental solutions to the variational equation of all the three pieces.
\item Step 2: the boundary terms $\mathrm{Id}+  JD_XH\otimes \nabla_X \ell_1$ for all the three pieces. 
\item Step 3: the coordinate change from left to right and from right to left on the middle section $\cS^M$.
\end{itemize}
{\bf Step 1: the variational equation.} The Hamiltonian written in the Delaunay coordinates  for the I4BP in the left Jacobi coordinate system has the form  (see \eqref{EqHamGlob})
  \begin{equation}
\begin{aligned}
H_{:\cdot\cdot}&=-\frac{M_0k_0^2}{2L_0^2}+\frac{M_1k_1^2}{2L_1^2} +\frac{M_2k_2^2}{2L_2^2} +U,\quad U=U_{01}+U_2,
\end{aligned}
\end{equation}
where $U_2$ and $U_{01}$ are the restriction of that in \eqref{EqHamGlob} to the isosceles case: %$k_0=1,\; k_1=2m_1,\; k_2=(m_1+2)m_2$,
  $$U_{01}= \frac{k_1}{r_1}-\frac{2m_1}{(r_1^2+r_0^2/4)^{1/2}}, \quad U_2=\frac{k_2}{ r_2}-\frac{m_1m_2}{r_2-\frac{2r_1}{m_1+2} } -\frac{2m_2}{\big((r_2+\frac{m_1r_1}{m_1+2})^2+r_0^2/4\big)^{1/2}}.$$
Note that  when $ r_1\gg r_0$, we have
$U_{01}=\frac{k_1}{r_1}(1-\frac{1}{(1+r_0^2/4r_1^2)^{1/2}})=\frac{k_1r_0^2}{8r_1^3}+O(\frac{r_0^4}{r_1^5}).$
We have the following estimates for the Hamiltonian equations  for the escaping piece of the oribts (see equation (7.5) of \cite{GHX}),
\begin{equation}\label{EqHamLimit}
\begin{cases}
\dot L_0&=\left(\frac{\partial U_{01}}{\partial r_0} +\frac{\partial U_2}{\partial r_0}\right)\frac{\partial r_0}{\partial \ell_0}=O\left(\frac{1}{\ell_1^3}+\frac{1}{\chi^3}\right)r_0\frac{\partial r_0}{\partial \ell_0}\\%=-\frac{1/2}{(r_1^2+r_0^2/4)^{3/2}}\frac{L_0^4}{(M_0k_0)^2}\sin u_0+\frac{1}{\chi^3},\\
\dot \ell_0&=\frac{M_0k_0^2}{L_0^3}+\left(\frac{\partial U_{01}}{\partial r_0} +\frac{\partial U_2}{\partial r_0}\right)\frac{\partial r_0}{\partial L_0}=\frac{M_0k_0^2}{L_0^3}+O\left(\frac{1}{\ell^3_1}+\frac{1}{\chi^3}\right) r_0\frac{\partial r_0}{\partial L_0}\\%=\frac{M_0k_0^2}{L_0^3}+\frac{r^2_0}{L_0(r_1^2+r_0^2/4)^{3/2}}+\frac{1}{\chi^3},\\
\dot L_1&=\left(\frac{\partial U_{01}}{\partial r_1} +\frac{\partial U_2}{\partial r_1}\right)\frac{\partial r_1}{\partial \ell_1}=O\left(\frac{1}{\ell_1^4}+\frac{1}{\chi^2}\right)\frac{\partial r_1}{\partial \ell_1}\\%=\left(\frac{2}{r_1^3}-\frac{4}{(r_1^2+r_0^2/4)^{3/2}}\right)\left(\frac{L_1^2}{M_1k_1}\right)^2(-\sinh u_1)+\frac{1}{\chi^2},\\
\dot \ell_1&=-\frac{M_1k_1^2}{L_1^3}+\left(\frac{\partial U_{01}}{\partial r_1} +\frac{\partial U_2}{\partial r_1}\right)\frac{\partial r_1}{\partial L_1} =-\frac{M_1k_1^2}{L_1^3}+O\left(\frac{1}{\ell_1^4}+\frac{1}{\chi^2}\right)\frac{\partial r_1}{\partial L_1}\\
%=-\frac{M_1k_1^2}{L_1^3}+\left(\frac{2}{r_1^2}-\frac{4r_1}{(r_1^2+r_0^2/4)^{3/2}}\right)\frac{2 r_1}{L_1}+\frac{1}{\chi},\\
\dot L_2&=-\frac{\partial U_2}{\partial r_2}\frac{\partial r_2}{\partial \ell_2}=O(\frac{1}{\beta^2\chi^2}),\\
\dot \ell_2&=\frac{M_2k_2^2}{L_2^3}+\frac{\partial U_2}{\partial r_2}\frac{\partial r_2}{\partial L_2}=\frac{M_2k_2^2}{L_2^3}+O(\frac{\beta}{\chi}).%=\frac{M_2k_2^2}{L_2^3}+\frac{1}{\chi}.
\end{cases}\end{equation}
Due to Lemma \ref{L01-bound}, we learn that $0<\frac{1}{C'}\leqslant |\dot\ell_1|\leqslant C'$ for some suitable $C'>1$.  So  we  treat $\ell_1$ as the new time variable. Then we  remove the equation for variable $ L_1$ from the system of equations~\eqref{EqHamLimit}  and divide all the equations for the remaining variables $(L_0,\ell_0, L_2,\ell_2)$ by $\dot\ell_1$. Denote by $X=(L_0,\ell_0,L_2,\ell_2)$, by $D_X=(\frac{\partial }{\partial L_0},\frac{\partial }{\partial \ell_0},\frac{\partial }{\partial L_2},\frac{\partial }{\partial \ell_2} )$ the partial derivarives and by $\nabla_X$ the covariant derivative, $\nabla_X=D_X+(D_X L_1) \partial_{L_1}$, which takes into account of the dependence on $X$ through $L_1$.
The Hamiltonian equation now has the form
$\frac{dX}{d\ell_1} =\frac{1}{\dot \ell_1}J D_X H$ and the variational equation has the form (see equation (C.2) of \cite{GHX})
\begin{equation}\label{EqVar}
\begin{aligned}
\frac{d}{d\ell_1} (\dt X)&=\nabla_X\left(\frac{1}{\dot \ell_1}J D_X H\right)\dt X= \left(\frac{1}{\dot \ell_1}J \nabla_X(D_X H)-\frac{1}{(\dot \ell_1)^2}J D_X H\otimes \nabla_X\dot \ell_1\right)\dt X \\
 &=\left(\frac{1}{\dot \ell_1}J (D^2_X H+D_{L_1}D_X H\otimes D_X L_1)-\frac{1}{(\dot \ell_1)^2}J D_X H\otimes \nabla_X\dot \ell_1\right)\dt X
\end{aligned}
\end{equation}
 After a lengthy but direct calculation (see equation (C.6) of \cite{GHX} for details), we get the following form of the coefficient matrix (using the fact $L_0=O(1),\ L_1=O(1)$ and $L_2=O(\beta^{-1})$),
\begin{equation}\label{coefficients-M1}\nabla_X\left(\frac{1}{\dot \ell_1}J D_X H\right)=O\left[
\begin{array}{cccccc}
\frac{1}{\ell_1^3}& \frac{1}{r_0\ell_1^3}&\frac{\beta}{\chi^2}+\frac{\beta^3}{\ell_1^3}&\frac{1}{\chi^2\beta}\\
1& \frac{1}{\ell_1^3}&\beta^3&\frac{1}{\chi^2\beta^2}\\
\frac{1}{\beta^2\chi^2}&\frac{1}{\beta^2\chi^3}&\frac{1}{\chi^2\beta}
&\frac{1}{\chi^3\beta^4}\\
\beta^3&\frac{\beta^3}{\ell_1^3}+\frac{1}{\beta^2\chi^2}&\beta^4&\frac{1}{\chi^2\beta}
\end{array}\right].
\end{equation}
Using Picard iteration (which stabilizes quickly and two steps are enough, see equation (C.7) of \cite{GHX}),
we have the following estimate for the fundamental solution of the corresponding variational equation (see equation (C.8) of \cite{GHX}),
\begin{equation}\label{D-X-L}
\cD_E= \Id+O\left[
\begin{array}{cccc}
 \epsilon  & \epsilon ^2 & \beta ^3 & \frac{1}{\beta  \chi } \\
 \chi  & \epsilon ^2 \chi  & \beta ^3 \chi  & \frac{1}{\beta } \\
 \frac{1}{\beta ^2 \chi } & \frac{\epsilon ^2}{\beta ^2 \chi } & \frac{1}{\beta  \chi } & \frac{1}{\beta ^3 \chi ^2} \\
 \beta ^3 \chi  & \beta ^3 \epsilon ^2 \chi  & \beta ^4 \chi  & \beta ^2 \\
\end{array}
\right].
\end{equation}

For the returning piece of the  orbit going from   $\cS^M$ to $\cS_-$, we have for almost the whole time $L_0=O(1), \ L_2=O(1)$ and  $L_1=O(f_{\mu,\lambda})$ (c.f. Proposition \ref{PropDG1} for the definition of $f_{\mu,\lambda}$), since $L_1\sim |\bp_1|^{-1}$ and $|\bp_1|$ is estimated in Proposition \ref{choiceofm2}, hence $\dot{\ell_1}\sim \frac{1}{L_1^3}=O(f^{-3}_{\mu,\lambda})$. Thus, in \eqref{coefficients-M1}, we replace $\beta$ by 1 since $L_2$ has changed from $O(\beta^{-1})$ to $O(1)$, and  multiply an additional factor $O(f^3_{\mu,\lambda})$ due to dividing $\dot \ell_1$. Next, by equation \eqref{EqCarDel}, we get $|\bx_1|\sim L_1^2\ell_1$, thus when $|\bx_1|$ ranges in from $\eps^{-1}$ to $O(\chi)$, we get that $\ell_1$ ranges from $f^{-2}_{\mu,\lambda} \cdot\eps^{-1}$ to $f^{-2}_{\mu,\lambda} \cdot\chi$.   Then integrating  the variational equation over with respect to the time $\ell_1$, we have the following estimates for its fundamental solution,
\begin{equation}\label{D-X-R}
\cD_R=\Id+f_{\mu,\lambda}\cdot O\left[
\begin{array}{cccc}
   \epsilon ^2\chi  & \epsilon ^2 &   \epsilon ^2\chi  & \frac{1 }{\chi } \\
   \chi  & \epsilon  &    \chi  &1 \\
 \frac{1}{\chi } & \frac{\epsilon }{\chi ^2} & \frac{1}{\chi } & \frac{1 }{\chi ^2} \\
 \chi  & \epsilon  &   \chi  & 1\\
\end{array}
\right].
\end{equation}

For the intermediate piece of the orbit going from $\cS^M$ back to itself, we consider the right Jacobi coordinate system.  In the Delaunay coordinates, the Hamiltonian function \eqref{3pairHam} for the I4BP reads
$$H_{:\cdot\cdot}=-\frac{M_0(k_0)^2}{2L_0^2}+\frac{M_1^R(k_1^R)^2}{2(L_1^R)^2}+\frac{M_2^R(k_2^R)^2}{(L_2^R)^2}+U_{12}$$
where we have
$$U_{12}=\frac{k_2^R}{r_2^R}-\frac{2m_1}{(\frac{1}{4}r_0^2+(r_2^R+\frac{m_2}{m_1+m_2}r_1^R)^2)^{1/2}}-\frac{2m_2}{(\frac1 4 r_0^2+(r_2^R-\frac{m_1}{m_2+m_1}r_1^R)^2)^{1/2}}.$$
In particular, when $r_2\gg r_1$, we have 
$U_{12}=\frac{k_2^R}{r_2^R}-\frac{2m_1+2m_2}{(\frac{1}4r_0^2+(r_2^R)^2)^{1/2}}+O(\frac{(r_1^R)^2}{(r_2^R)^3}).$
We have the following estimates for the system of equations as $\chi\to\infty$
\begin{equation}\label{RightHam}
\begin{cases}
\dot L_0=-\frac{\partial U_{12}}{\partial \ell_0}=O(\frac{1}{\chi^3}),\;
\dot \ell_0=\frac{M_0k_0^2}{L_0^3}+\frac{\partial U_{12}}{\partial L_0}=O(1),\\
\dot L^R_1=-\frac{\partial U_{12}}{\partial\ell^R_1}=O(\frac{\ell^R_1}{\chi^3}),\;
\dot\ell^R_1=-\frac{M_1^R(k_1^R)^2}{(L_1^R)^3}+\frac{\partial U_{12}}{\partial L^R_1}=O(1),\\
\dot L_2^R=-\frac{\partial U_{12}}{\partial \ell^R_2}=O(\frac{\ell^R_2}{\chi^3}),\;
\dot \ell_2^R=-\frac{M_2^R(k_2^R)^2}{(L_2^R)^3}+\frac{\partial U_{12}}{\partial L^R_2}=O(1).
\end{cases}\end{equation}
Again we treat $\ell_1^R$ as the new time and remove $L_1^R$ from the equation. Then using \eqref{EqVar} and repeating the calculations for \eqref{coefficients-M1}, we have the following form for the coefficient matrix of the corresponding variational equation,
\[\nabla_X\Big(\frac{1}{\dot{\ell_1}^R}JD_XH\Big)=O\left[\begin{array}{cccc}\frac{1}{\chi^3}&\frac{1}{r_0\chi^3}&\frac{1}{\chi^3}&\frac{1}{\chi^4}\\
1&\frac{1}{\chi^3}&1&\frac{1}{\chi^2}\\
\frac{1}{\chi^2}&\frac{1}{\chi^4}&\frac{1}{\chi^2}&\frac{1}{\chi^3}\\
1&\frac{1}{\chi^3}&1&\frac{1}{\chi^2}\end{array}\right]\]
The timespan of the intermediate piece of the orbit is of order $O(\beta\chi)$, therefore, we have the following estimate for the fundamental solution
\begin{equation}\label{D-X-I}\cD_I=\Id+O\left[\begin{array}{cccc}\frac{1}{\chi^2}&\frac{1}{\chi^2}&\frac{1}{\chi^2}&\frac{1}{\chi^3}\\
\beta\chi&\frac{1}{\chi^2}&\beta\chi&\frac{1}{\chi^1}\\
\frac{1}{\chi}&\frac{1}{\chi^3}&\frac{1}{\chi}&\frac{1}{\chi^2}\\
\beta\chi&\frac{1}{\chi^2}&\beta\chi&\frac{1}{\chi}\end{array}\right].
\end{equation}
{\bf Step 2, the boundary contributions.}

The estimate of the boundary contributions is almost identical to that of \cite{GHX} Appendix C, Step 2. The boundary terms $(\Id+\nabla_X\ell_{1}\otimes JD_XH)$ are calculated on the sections $\cS_+$, $\cS_-$ and $\cS^M$.  For the middle section $\cS^M$, we need to compute the boundary contribution before and after  changing  Jacobi coordinates from the left one  to the  right one, and vise versa.
We first consider the left side of the  section $\cS^M=\{ |\bx_1^R|=\beta\chi\}$ for the escaping piece of the orbit. Recall that $X=(L_0,\ell_0,L_2^L,\ell_2^L)$. For the Hamiltonian equation, we have $JD_XH=O\left(\frac{1}{\chi^3},1,\frac{1}{(\beta\chi)^2} ,\beta^3\right).$ Using \eqref{L-R-transit}, We get $ \bx_1^R=\frac{2}{m_1+2}\bx_1^L-\bx_2^L=\beta\chi$ on $\cS^M$. Converting the Cartesian variables to Delaunay variables using formula from \cite[Lemma B.1]{GHX} and taking derivatives,  we have 
$$\frac{2}{m_1+2}\frac{1}{M_2^Lk_2^L}(2L_2^L\ell_2^LdL_2^L+(L_2^L)^2d\ell_2^L)-\frac{1}{M_1^Lk_1^L}(2L_1^L\ell_1^LdL_2^L+(L_1^L)^2d\ell_1^L)=0.$$ Then by the implicit function theorem, we get $D_X\ell_1^L=O(0,0,\beta\chi,\frac{1}{\beta^2}),\; \partial_{L_1^L}\ell_1^L\simeq\chi.$
So, from the left side of the section $\cS^M$, we have the boundary contribution (see equation (C.10) of \cite{GHX}),
\begin{equation}\cD^E_M=\Id+O\left(\frac{1}{\chi^3},1,\frac{1}{(\beta\chi)^2} ,\beta^3\right)\otimes O\left(\chi,\frac{1}{\chi^2}, \beta\chi,\frac{1}{\beta^2}\right) \end{equation}
Similarly, for the returning piece of the orbit, on the left side of $\cS^M$,  we have $D_X\ell_1^L=O(0,0, \mu^{1/2}\beta\chi, \mu^{1/2}\beta)$ and $\partial_{L_1^R}\ell_1^R=O(\mu^{1/2}\beta\chi)$. Then  the boundary contribution takes the form,
\[D_M^R=\Id+O(\frac{1}{\chi^3},1,\frac{1}{\chi^2},1)\otimes(\mu^{-1}\beta^{-2}\chi,\frac{1}{\chi^2},\mu^{-1}\beta^{-2}\chi,\mu^{1/2}\beta).\]
Next, we work from the right side of the  section $\cS^M$. On the right side, $X=(L_0,\ell_0,L^R_2,\ell^R_2)$. For the Hamiltonian vector filed, we have
$JD_XH=O\left(\frac{1}{\chi^3},1,\frac{1}{\chi^2} ,1\right).$
Once again we substitute the formulas from  \cite[Lemma B.1]{GHX} and take differentiation. We obtain,
$\frac{1}{M^R_1k^R_1}(2L_1^R\ell_1^RdL_1^R+(L_1^R)^{2}d\ell_1^R)=0.$
Using  implicit function theorem we have
$D_X\ell_1^R=(0,0,0,0),\; \partial_{L_1^R}\ell_1^R\sim \beta\chi$. 
Therefore, from the right side of the section $\cS^M$, we have the boundary contribution which is valid for both the escaping piece and returning piece of orbits (see Step 2 of Appendix C of \cite{GHX})
\[\cD^I_M=\Id+O\left(\frac{1}{(\beta\chi)^3},1,\frac{1}{\chi^2} ,1\right)\otimes O\left(\beta\chi,\frac{1}{(\beta\chi)^2},\beta \chi,\frac{1}{\chi}\right).\]

Finally, we consider the boundary contributions on the sections $\cS_\pm$. Both were wokred out in Step 2 of Appendix C of \cite{GHX}, and we cite the result here directly. On the section $\cS_+$,  we have the boundary contribution 
\[\cD_+=\Id+O(\eps^3,1,\frac{1}{\beta^2\chi^2},\beta^3)\otimes O\left(\eps^{-1},1,\eps^{-1}\beta^3,\frac{\eps^{-1}}{(\beta\chi)^2}\right),\]
and on the section $\cS_-=\{r=\epsilon^{-1}\}$,
 the boundary contribution is,
\[
\cD_-=\Id+O(\eps^3,1,\frac{1}{\chi^2},1)\otimes O\left(\eps^{-1},\eps^{2},\eps^{-1},\frac{\eps^{-1}}{\chi^2}\right).
\]

{\bf Step 3, transition from the left to the right. }

At the middle section $\cS^M$, we change the Delaunay coordinates on the left first to the left Jacobi coordinates, next from the left Jacobi to the right Jacobi coordinates, and finally from the right Jacobi to Delaunay coordinates on the right, and vice versa. This was done in Step 3 of Appendix C of \cite{GHX}. 
We have   the  following transition matrix on the section $\cS^M$ (see equation (C.11) of \cite{GHX})
\begin{equation}\label{middle-d-l}\cM_R=\frac{\partial(L_0,\ell_0,L_2,\ell_2)^R|_{\cS^M}}{\partial(L_0,\ell_0,L_2,\ell_2)^L|_{\cS^M}}=O\left[\begin{array}{cccc}
1&0&0&0\\
0&1&0&0\\
1&\frac{1}{\chi^3}&\beta&\frac{1}{\beta^2\chi}\\
\chi&\frac{1}{\chi^2}&\beta\chi&\beta^{-2}
\end{array}\right].
\end{equation}
Similarly, for  the returning piece leaving $\cS^M$,  we have the following transition matrix from the right Jacobi coordinates to the left one,
\begin{equation}\label{middle-d-l}\cM_L=\frac{\partial(L_0,\ell_0,L_2,\ell_2)^L|_{\cS^M}}{\partial(L_0,\ell_0,L_2,\ell_2)^R|_{\cS^M}}=O\left[\begin{array}{cccc}
1&0&0&0\\
0&1&0&0\\
\beta&\frac{\beta}{\chi^3}&\beta&\frac{\beta}{\chi}\\
\beta\chi&\frac{\beta}{\chi^2}&\beta\chi&\beta
\end{array}\right].
\end{equation}

After completing the three steps, using the formula \eqref{eq: formald4}, we multiply all the matrices together $(\cD_-\cD_R\cD_M^R)\cM_L(\cD_M^I\cD_I\cD_M^I)\cM_R(\cD_M^E\cD_E\cD_+)$ and obtain \eqref{DGG}.  Hence, we complete the proof of Proposition \ref{PropDG1}. \qed

\section{The derivative of the global map in the F4BP}\label{SG0g0}
In this section, we study the $C^1$ estimate of the global map in  the F4BP by allowing the variables $G_i,g_i$, $i=0,1,2$, to be non-vanishing. Again the proof is very similar to that of \cite{GHX} and the difference is as discussed at the beginning of the last appendix. 

We shall first work out the estimate of $d\mathbb G$ in the isosceles limit $Y\to0$, then handle the nonisoscelesness $Y\neq 0$.

\subsection{The $C^1$-dynamics of the variables $G_i$, $g_i$, $i=0,1,2$ in the isosceles limit}\label{SSGlob12}

In this section, we compute the derivative of the global map for the variables $(G_0,g_0, G_1,g_1,G_2,g_2)$.  As we have done in Appendix \ref{L-Derivatives},  we break the orbit from $\cS_+$ to $\cS_-$ into three pieces by the middle section $\cS^M$. For each piece, we apply equation \eqref{eq: formald4} to compute the derivative, then we also compute the derivative of the coordinate change from left to right and right to left on the middle section $\cS^M$.  Since we have the Hamiltonian equations $\dot Y_i=J D_{Y_i}H=0$ along an orbit of the I4BP, where  $J=\left[\begin{array}{cc}0&-1\\
1&0\end{array}\right]$ and $Y_i=(G_i,g_i),\ i=0,1,2$, we get that the boundary contributions in \eqref{eq: formald4} are all identity, so in the following, we only work on the fundamental solutions to the variational equation as well as the left-right transitions. We stress that the contribution from the latter is essential for the $(G_1,g_1,G_2,g_2)$-component.

We split the proof into two steps:
\begin{enumerate}
	\item  the fundamental solution of the variational equation,
	\item  the transition from  left to right and right to left.
	\end{enumerate}
It turns out that the estimate in item (1) is $O(1)$ and the $O(\chi^2)$ expansion in the statement is given by item (2). %To get the $O(1)$ estimate for item (1) is easy, but in the following, we 

For the orbits   under consideration, they have parameters $L_1^L=O(1)$, $L_2^L=O(\lambda^{1/2})=O(\beta^{-1})$ for the escaping pieces, $L^R_1=O(1)$, $L^R_2=O(1)$ for the intermediate pieces and $L_1^L =O(f_{\mu,\lambda})$ (see  Appendix C, Step 1), $L_2^L=O(1)$ for the returning pieces.  

{\bf Step 1,  the variational equation. }

Consider the variational equation evaluated on orbits between  $\cS_+$ and $\cS^M$: the escaping piece. Denoting $Y_0=(G_0,g_0)$, $Y_1=(G_1,g_1)$, $Y_2=(G_2,g_2)$ and $Y=(Y_1,Y_2)$,  then the Hamiltonian equation can be written as $\dot Y_0=JD_{Y_0}Y_0$, $\dot Y_1=JD_{Y_1}H,\ \dot Y_2=JD_{Y_2}H$ 
and the variational equation can be written as
\begin{equation}\label{VarEq-012}
\left[\begin{array}{c}
\dot{\dt Y}_0\\
\dot{\dt Y}_1\\
\dot{\dt Y}_2
\end{array}\right]=\left[\begin{array}{ccc}
JD_{Y_0}^2H&JD_{Y_0Y_1}^2H&JD_{Y_0Y_2}H\\
JD_{Y_1Y_0}^2H&JD_{Y_1}^2 H &JD_{Y_2Y_1}^2 H\\
JD_{Y_2Y_0}^2H&JD_{Y_1Y_2}^2 H&JD_{Y_2}^2 H
\end{array}\right]\left[\begin{array}{c}
\dt Y_0\\
\dt Y_1\\
\dt Y_2
\end{array}\right].
\end{equation}

By straightforward estimates we know that the contributions from the entries $JD^2_{Y_0Y_i}$, $i=0,1,2$, to the solution to the variational equation are negligible (see equation (9.3) of \cite{GHX}). 
The main contribution comes from the block 
$\left[\begin{array}{cc}
JD_{Y_1}^2 H &JD_{Y_2Y_1}^2 H\\
JD_{Y_1Y_2}^2 H&JD_{Y_2}^2 H
\end{array}\right]$ (see equation (9.2) of \cite{GHX}). % whose leading term is of rank 1, taking the form $O(\frac{r_1}{\chi^2})(1,-\frac{1}{L_1^L},1,-\frac{1}{L_2^L})^T\otimes(\frac{1}{L_1^L},1,-\frac{1}{L_2^L},-1)$ (see equation (9.2) of \cite{GHX}). 

Between the sections $\cS_+$ and $\cS^M$ the quantity $r_1$ increase almost linearly with respect to time from the order $\epsilon^{-1}$ to the order $\chi$. Then we integrate the  equation \eqref{VarEq-012} over a time interval of order $\chi$, to get the following estimate of the fundamental solution (see equation (9.8) of \cite{GHX}),
\begin{equation}\label{D-G-L}
\left[\begin{array}{cc}\mathbb{I}_{2\times2}+O_{2\times2}(\epsilon^2)&O_{2\times4}(\epsilon^2)\\
O_{4\times2}(\epsilon^2)&\mathcal{D}_G^L\end{array}\right].
\end{equation}
where the matrix $\mathcal{D}_G^L$ is the fundamental solution of the variational equation given by $\left[\begin{array}{cc}
	JD_{Y_1}^2 H &JD_{Y_2Y_1}^2 H\\
	JD_{Y_1Y_2}^2 H&JD_{Y_2}^2 H
\end{array}\right]$ and it has some structure that enables us to work it out explicitly (see equation (9.5) and (9.6) of \cite{GHX}). 
%where $\cD_L^G=\mathbb{I}_{4\times4}+A(1,-\frac{1}{L_1^L},1,-\frac{1}{L_2^L})^T\otimes(\frac{1}{L_1^L},1,-\frac{1}{L_2},-1)+O_{4\times4}(\epsilon^2)$ with $A\neq0$ being a constant of order 1 that can be computed explicitly. 

Similarly, for the intermediate piece of the orbits from the section $\cS^M$ to itself, we integrate the estimate of the variational equations in equation (9.2) of \cite{GHX} over time of order $\beta\chi$ to get the following form for the  fundamental solution of the corresponding variational equation,
\begin{equation}\label{D-G-I}\left[\begin{array}{cc}\mathbb{I}_{2\times 2}+O_{2\times2}(\frac{1}{\chi})&O_{2\times4}(\frac{1}{\chi})\\
O_{4\times2}(\frac{1}{\chi})&\cD_G^I\end{array}\right],
\end{equation}
where $\cD_I^G=\mathrm{Id}+O(\beta_{4\times 4})$. 
	
Finally, it follows from similar calculations that for the returning piece of the orbits between the sections $\cS^M$ and $\cS_-$, integrating the estimate of the variational equations in equation (9.2) of \cite{GHX} over time $\ell_1$ ranging from $f_{\mu,\lambda}^{-2}\eps^{-1}$ to $f_{\mu,\lambda}^{-2}\chi$ (see the derivation of  \eqref{D-X-R}), we get  the following estimates for its fundamental solution,
\begin{equation}\label{D-G-R}
\left[\begin{array}{cc}\mathbb{I}_{2\times2}+O_{2\times2}(\epsilon^2f_{\mu,\lambda})&O_{2\times4}(\epsilon^2f_{\mu,\lambda})\\
O_{4\times2}(\epsilon^2f_{\mu,\lambda})&\cD_G^R\end{array}\right],
\end{equation}
where the leading terms in $\cD_G^R$ are given by $f_{\mu,\lambda}\cD_G^L.$
%where $\cD_{R}^G=\mathbb{I}_{4\times4}+A'(1,-\frac{1}{\bar L_1^L},1,-\frac{1}{\bar L_2^L})^T\otimes(\frac{1}{\bar L_1^L},1,-\frac{1}{\bar L_2^L},-1)+O_{4\times4}(\epsilon^2)$ with $A'\neq0$ being a constant that can be computed explicitly. 

{\bf Step 2, Transition from the left to the right.}

We next consider the coordinate changes on the section $\cS^M$ induced by switching between left and right Jacobi coordinates. For the escaping orbit passing through $\cS^M$, we shall first convert the $Y$-variables into the left Jacobi coordinates, next convert from the left to right Jacobi coordinates using \eqref{L-R-transit}, then convert the right Jacobi coordinates back to the $Y$ varriables. We reverse the procedure for the returning orbit passing through $\cS^M$. It turns out that each of the transition gives us a factor $\chi$, that is why we have $O(\chi^2)$ in the statement. 
\begin{Lm}\label{LmTransitionGg}
We have the following derivative estimate on the section $\cS^M=\{|\bx_1^R|=\beta\chi\}$:

(1) When the escaping orbit arrives at $\cS^M$, for the left-right transition, we have 
$$\frac{\partial(G_0,g_0,G_1,g_1,G_2,g_2)^{R}}{\partial(G_0,g_0,G_1,g_1,G_2,g_2)^{L}}\Big|_{\cS^M}=\left[\begin{array}{cc}\mathbb{I}_{2\times2}&0\\
0&\mathbf{v}^R\otimes dG_1^R+O(1)\end{array}\right],$$
where $\mathbf{v}^R=(1,-\frac{1}{L_1^R},-1,\frac{1}{L_2^R})^T$ and
$dG_1^R:=\frac{\partial G_1^R}{\partial Y^L}\simeq\chi\cdot(\frac{1}{L_1^L},1,-\frac{1}{L_2^L},-1)$ with $L_1^L,L_1^R,L_2^R=O(1)$ and $L_2^L=O(\beta^{-1})$;

(2) When the returning orbit leaving $\cS^M$, for the right-left transition, we have  (we add a bar to the variables to distinguish them from the last item)
$$\frac{\partial(\bar G_0,\bar g_0, \bar  G_1,\bar g_1,\bar G_2,\bar g_2)^{L}}{\partial(\bar G_0,\bar g_0,\bar G_1,\bar g_1,\bar G_2,\bar g_2)^{R}}\Big|_{\cS^M}=\left[\begin{array}{cc}\mathbb{I}_{2\times2}&0\\
0&\mathbf{v}^L\otimes d\bar G_1^L+O(1)\end{array}\right],$$
where $\mathbf{v}^L=(1,-\frac{1}{\bar{L}_1^L},-1,\frac{1}{\bar{L}_2^L})^T$ and
$d\bar G_1^L:=\frac{\partial \bar G_1^R}{\partial \bar Y^L}\simeq\chi\cdot(-\frac{1}{\bar L_1^R},1,-\frac{1}{\bar L_2^R},-1)$ with $\bar L_2^L,\bar L_1^R,\bar L_2^R=O(1)$ and $\bar L_1^L=O(f_{\mu,\lambda})$.
\end{Lm}
For a proof of this lemma, we refer the reader to \cite[Appendix E]{GHX} where the statement (1) is proved and the proof of statement (2) is almost identical. 

 Hence, after combing \eqref{D-G-L}-\eqref{D-G-R} and Lemma \ref{LmTransitionGg}, the derivatives of the global map for the variables $G_i,g_i$, $i=0,1,2$, between the sections $\cS_+$ and $\cS_-$ takes the following form,
\begin{equation}\label{gs-tensor1}
\frac{\partial(G_0,g_0,G_1,g_1,G_2,g_2)\big|_{\cS_-}}{\partial(G_0,g_0,G_1,g_1,G_2,g_2)\big|_{\cS_+}}=\left[\begin{array}{cc}\mathbb{I}_{2\times2}&0\\
0&\chi^2\bar{\mathbf{u}}_2\otimes\bar{\mathbf{l}}_2\end{array}\right]+O((\eps^2+\beta)f_{\mu,\lambda}\chi),
\end{equation}
where $\bar{\bu}_2= \cD_G^R\mathbf{v}^L+O(\epsilon^2)=(1,-\frac{1}{\bar{L}_1^L},-1,\frac{1}{\bar{L}_2^L})+O(\epsilon^2)$ and $\bar{\mathbf{l}}_2=\frac{1}{\chi}dG_1^R \cD_G^L=f_{\mu,\lambda}\cdot(\frac{1}{L_1^L},1,-\frac{1}{L_2^L},-1)+O(\epsilon^2)$. For this calculation, we need to know $\cD_G^L$ and $\cD_G^R$ explicitly and it turns out that due to their special structure, multiplying them does change the vectors $\mathbf v^L$ and $dG_1^R$ much. We refer readers to equations (9.5), (9.6) and (9.10) of  \cite{GHX} for details of this calculation. 

When taking into account the conservation of angular momentum $G_0+G_1+G_2=0$, in the reduced system of coordinates  $(G_0,\mathsf g_{01},G_2,\mathsf g_{12})$. The details of this reduction are given in Section 9.3 of \cite{GHX}, which gives   the derivatives of the global map in the isosceles limit $Y\to0$
\[\left[\begin{array}{cc}\mathbb{I}_{2\times2}&0_{2\times2}\\
0_{2\times2}&0_{2\times2}\end{array}\right]+\chi^2f_{\mu,\lambda}\bar\bu\otimes\bar\bl+O((\eps^2+\beta)f_{\mu,\lambda}\chi),\]
where $\bar\bu=(0,\frac{1}{\bar{L}_2^L},1,-\frac{1}{\bar L_1^L}-\frac{1}{ \bar L_2^L})+O(\epsilon^2+\beta)$ and $\bar\bl=(0,0,\frac{1}{L_1^L},1)+O(\epsilon^2+\beta)$. 
%This completes the proof of Proposition \ref{PropDG1}.
This gives the $O(\chi^2)$ terms in the form of $d\bG$ in  Proposition \ref{PropDG2}. 
\subsection{The non-isosceles case}\label{non-isosceles-sec}
We next include the perturbations coming from the non-isoscelesness and complete the proof of Proposition \ref{PropDG2}. We shall control $|Y|\leq \nu$ for some small $\nu$, and show that the non-isoscelesness only contribute the $O(\nu f_{\mu,\lambda}\chi)$ error term to $d\bG$. The main difficulty in the section is to prove \eqref{EqDG} and \eqref{EqDG>}, which needs a careful analysis of the special structure of the fundamental solution to the variational equation. 
\begin{proof}[Proof of Proposition \ref{PropDG2}]
Between the section $\cS_+$ and $\cS^M$, let us consider the non-isosceles escaping orbits such that 
\begin{equation}\label{small-gs}|G_0|,|G_1|,|G_2|,|g_0-g_1-\frac{\pi}{2}|,|g_1-g_2|\leqslant \nu, \quad 0<\nu\ll1.
\end{equation}
From \eqref{D-G-L} we know that the above conditions are satisfied if they are valid on the section~$\cS_+$. 
We compute the  variational equation of the variables $(X,Y)$ where $X=(L_0,\ell_0,L_2,\ell_2)$ and  $Y=(G_0,g_0,G_1,g_1,G_2,g_2)$ along such an orbit. Recall that we reduce the variable $L_1$ from the conservation of energy and treat the variable $\ell_1$ as the new time,  we then have in the leading term 
\begin{equation}\label{VarEq-XY}
\left[\begin{array}{c}\dot{\dt X}\\
\dot{\dt Y}\end{array}\right]=\left[\begin{array}{cc}K_{11}& K_{12}\\
K_{21}&K_{22}\end{array}\right]\left[\begin{array}{c}\dt X\\
\dt Y\end{array}\right],
\end{equation}
where $K_{11}$ is the coefficient matrix \eqref{coefficients-M1} of the variational equation for $X$, and $K_{22}$ is the coefficient matrix in \eqref{VarEq-012}. We obtain the following estimates for its fundamental solution, 
\begin{equation}\label{DGleft}\left[\begin{array}{cc}(1+O(\nu))\mathcal{D}_2^L&\mathcal{K}^L_{12}\\
O_{6\times 4}(\nu)&(1+O(\nu))\mathcal{D}_G^L\end{array}\right],
\end{equation}
where $\mathcal{D}_2^L$ is in \eqref{D-X-L}, $\mathcal{D}_G^L$ is in \eqref{D-G-L}, and $\mathcal{K}^L_{12}=\nu \eps^2\chi \bar \bu_1\otimes \hat\bl_1+ O(\nu(\eps^3+\eps^2\beta)\chi)$
with $\bar{\bu}_1$, $\bar{\bl}_1$ from Proposition \ref{PropDG1}, $\hat\bl_1=\bar\bl_1\int_0^\chi  K_{12}(s)\mathbb Y(s)\,ds$ and $\mathbb{Y}(s)$ being  the fundamental solution of \eqref{VarEq-012}. We refer readers to Step 1 of Section 9.4 of \cite{GHX} for details of this estimate. 

For the immediate piece of the orbits, which leave the section $\cS^M$ and then return to it, we have the the following estimate for the fundamental solution of the corresponding variational equation,
\begin{equation}
\label{DG-middle}\left[\begin{array}{cc}(1+O(\nu))\mathcal{D}^I&O_{4\times6}(\nu\beta)\\
O_{6\times 4}(\nu\beta)&(1+O(\nu))\mathcal{D}_G^I\end{array}\right],
\end{equation}
where $\cD^I$ is as in \eqref{D-X-I} and $\cD_G^I$ is as in \eqref{D-G-I}. 

For the  returning orbits between the section $\cS^M$ and $\cS_-$ satisfying the condition \eqref{small-gs}, we have the following estimate for the fundamental solution of the corresponding variational equation, 
\begin{equation}\label{DGright}\left[\begin{array}{cc}(1+\nu)\mathcal{D}_R&O(\nu\eps f_{\mu,\lambda})\\
{\mathcal{K}}_{21}^R&(1+O(\nu))\mathcal{D}_G^R\end{array}\right],\end{equation}
where $\mathcal{D}_R$ is in \eqref{D-X-R}, $\mathcal{D}_G^R$ is in \eqref{D-G-R},  and $\mathcal{K}_{21}^R=O(\nu\eps^2f_{\mu,\lambda}\chi).$  We refer readers to Step 2 of Section 9.4 of \cite{GHX} for details of this calculation, with modification being integrating $\ell_1$ over time interval $f^{-2}_{\mu,\lambda}\eps^{-1}$ to $f^{-2}_{\mu,\lambda}\chi$, for the same reason as the way we derive \eqref{D-G-R}. %Note that the estimates of the off-diagonal blocks in \eqref{DGright} are quite different from that in \eqref{DGleft}. In particular, the $\nu\chi$-dependent block in \eqref{DGleft} is the (1,2)-block while that in \eqref{DGright} is the (2,1)-block. The difference is caused by the fact that $r_1$ is increasing linearly in $t$ while $r_2$ is decreasing linearly in $t$. 
 
We next consider  the transition from the left to right on the section $\cS^M$. We get by equation (9.21) of \cite{GHX}
\begin{equation}\label{non-iso-transit}\frac{\partial (X^R,Y^R)}{\partial(X^L,Y^L)}\Big|_{\cS^M}=\left[\begin{array}{c|cc}\mathcal{M}_{4\times4}&\Big(\begin{array}{c}0\\O_{2\times2}(\frac{\nu}{\chi^2})\end{array}\Big)&\Big(\begin{array}{c}0\\ 
\nu_0\chi (O(1)_{2\times1})\otimes dG_1^R\end{array}\Big)\\
\hline
0&\mathbb{I}_{2\times2}&0\\
\Big(\begin{array}{ccc}O_{4\times2}(\frac{\nu}{\chi^2})&O_{4\times2}(\nu\chi)\end{array}\Big)&O_{4\times2}(\frac{\nu}{\chi^2})&\chi\mathbf{v}^R\otimes dG_1^R\end{array}\right]
\end{equation}
where $\mathcal{M}$ is the transition matrix for the $X^{L,R}$  given \eqref{middle-d-l} and $\mathbf{v}^R\otimes dG_1^R$ is in Lemma \ref{LmTransitionGg}. The same structure holds true the right-left transition with obvious modification.

To give the structure of $d\bG$ in the statement, we multiply the matrices obtained above. The nonisosceles $|Y|<\nu$ gives only $O(\nu\chi)$ contribution to $d\bG$ and does not spoil the $O(\chi^2)$ and $O(\chi)$ leading terms obtained in the isosceles limit. For \eqref{EqDG} and \eqref{EqDG>} in the statement, we refer readers to Step 4 of Section 9.4 of \cite{GHX} for the argument, which utilizes the special structure of the matrices \eqref{DGleft} and \eqref{DGright}. This completes the proof of the proposition. 
\end{proof}
\section{Numeric verification of the transversality conditions}\label{app-numeric}

Since the dynamics near triple collision is  nonperturbative, in particular, it is not easy to explicitly find $\gamma_I$ and $\gamma_O$ and integrate the variational equations along them, we choose to verify the transversality condition using a computer, as people have done in the literature (see \cite{SM} etc). We only consider the case with 3 equal masses, that is, $m_1=m_3=m_4=1$. So in this case $m_2=3+o(1)$ for being admissible as in Definition \ref{mass-admissible}. Then  the equations of motion and the variational equations are concrete ODEs. The numeric task is to solve for $\gamma_I$ and $\gamma_O$, then integrate along them the corresponding variational equation with the stable and unstable directions at the fixed point as initial conditions. Explicit formulas for the variational equations can be find in \cite[Section 8 and Appendix D]{ GHX}.

{\bf Step 1: } We numerically compute\footnote{We used the ODE solver 'ode113' in Matlab with 'RelTol'=$10^{-13}$ and 'AbsTol'=$10^{-15}$, and with the initial condition $O_1+10^{-14}\mathbf{e}^u$, where $\mathbf{e}^u$ is the unstable direction of the fixed point $O_1$ on the collision manifold.  }  $\gamma_{O}$ (going to the side $\psi=\frac{\pi}{2}$) by solving I3BP on the collision manifold  from  a very small neighborhood of a fixed point along the unstable direction to the section $\cS_+$, and integrate\footnote{We used the ODE solver 'ode113' in Matlab with 'RelTol'=$10^{-7}$ and 'AbsTol'=$10^{-8}$.} the corresponding variational equation  along this orbit for the variables $(dw_0,d\mathsf{g}_{01},dw_2,d\mathsf{g}_{12})$ with  initial conditions   $\mathbf{e}_1$ and $\mathbf{e}_2$, which are the unit eigenvectors corresponding to the two largest eigenvalues $\mu_u>0$ and $-\frac{v_*}{2}$, starting from a neighborhood of $O$. We perform the renormalization on the section $\cS_+$ and then change from the coordinate $(dw_0,d\mathsf{g}_{01},dw_2,\mathsf{g}_{12})$ to $(dG_0,d\mathsf{g}_{01},dG_2,d\mathsf{g}_{12})$. Numerically, we get $E^u_{4}(\gamma(0))=\text{span}\{\hat{\mathbf{e}}_1,\hat{\mathbf{e}}_2\}$, where $\hat{\mathbf{e}}_1=(0.5847,0.5608,0,0.5862)$ and  $\hat{\mathbf{e}}_2=( 0.4774,0.5282, -0.0701, 0.6987)$.

% {\bf Step 2:} {\it the piece of orbit shadowing $\gamma_I$. }

%  The transversality amounts to show that the vector $\hat{\be}_{0}$, when pushed forward along $\gamma^2_I$ approaching the Lagrange fixed point, does not coincides with the stable $E_0^s$ direction at the fixed point $O_2$ in the $(w_0,g_0)$-plane. To prove this, we first integrate the variational equation in Lemma~\ref{LmVar0} with reversed time from the fixed point to the section $\cS_-^R=\{r^R\simeq10^4\}$. The fundamental solution gives rise to a linear transform, sending the vector $E_0^s$ at the fixed point to a vector $\hat \be_1$ at $\cS^R_-$. Then on $\cS_-^R$, we change coordinates from $(w_0,g_0)$ to $(G_0,g_0)$. The derivative map of this change of coordinates on the section $\cS^R_-$ results in  multiplying the $\delta w_0$-component by $(r^R)^{-1/2}$. Then we get a vector $\hat{\mathbf{e}}_2$ in the $(G_0,g_0)$-coordinates.

 % {\bf The numeric result:}  
 
{\bf Step 2:} Next, we find the orbit $\gamma_I$ numerically. Note that for the I3BP, at the Lagrange fixed point $O$ with $v_*<0$, there are two linearly independent stable directions. One is on the collision manifold, which we denote as~$\mathbf{e}_1^s$, the other is $\mathbf{e}_2^s=(v_*,E,0,0)$, where $E$ is the total energy\footnote{By Proposition \ref{choiceofm2}, now   the total energy for the I3BP is $-1+o(1)$.}. So we  solve\footnote{We use the ODE solver `ode113' in Matlab with 'RelTol'=$10^{-13}$ and 'AbsTol'=$10^{-14}$.}  I3BP backward  along the stable direction with initial conditions of the form $$O+10^{-8}\big(-\cos(k10^{-3}\frac{\pi}{2})\mathbf{e}_2^s+\sin(k10^{-3}\frac{\pi}{2})\mathbf{e}_1^s\big), \quad k\in[-1000,1000]\cap\mathbb{Z}.$$ From the numeric result, there exist a consecutive set of  initial conditions ($k=496,\dots,504$) for which the corresponding backward solutions satisfy the energy condition  for $\gamma_I$ ($E_1$ is positive and of order $O(10^{-3})$) when the variable $r$ is of order $10^3$ and go to the side $\psi=\frac{\pi}{2}$.   We then  evaluate the corresponding variational equation   backward along the corresponding to all these orbits  with initial condition  $E_0^s(O)$ and then change from the coordinate $(dw_0,d\mathsf{g}_{01},dw_2,\mathsf{g}_{12})$ to $(dG_0,d\mathsf{g}_{01},dG_2,d\mathsf{g}_{12})$ on the section $\cS^R_-$. Numerically, we get that $E_{4}^s(\gamma_I(0))$ is spanned by $\{\hat{\mathbf{e}}_3,\hat{\mathbf{e}}_4\}$, where  $\hat{\mathbf{e}}_3=(0,0,0,1)$ and  $\hat{\mathbf{e}}_4=( -0.7864,0.4473,0,0.4261)$. We also have $E^s_3(\gamma_I(0))=\text{span}\{\bar{\mathbf{e}}_4:=( -0.8692,0.4944)\}$

  {\bf Step 3:} {\it Checking the cone conditions  in Proposition \ref{PropDG2}. } Since $\bar{\mathbf{l}}_2=(0,0,\frac{1}{L_1^L},1)+O(\beta)$, from the definition of $\bar{\mathbf{u}}_3$ and the numerics in Step 1, we have $\mathbf{u}=\frac{0.0701-0.6607L_1^L}{0.5862L_1^L}\hat{\mathbf{e}}_1+\hat{\mathbf{e}}_2+O(\beta)$ and setting the last two entries of $\mathbf{u}$ to zero we obtain $\bar{\mathbf{u}}_3$. Since by Proposition \ref{PropRenorm} $L_1^L=2/\sqrt{3}+o(1)$, we know that $\text{span}\{\bar{\mathbf{u}}_3\}\cap\mathcal{C}_{\frac{1}{100}}(E^s_{3}(\gamma_I(0)))=\{0\}$.
Since $\bar{\mathbf{u}}_2=(0,0,-1,-\frac{1}{\bar L_2^L})+o(1)$ and by Proposition \ref{choiceofm2} now $\frac{1}{\bar L_2^L}=2/(9\sqrt{3})+o(1)$,  with the numeric results from Step  2, we clearly have $\text{span}\{\bar{\mathbf{u}}_2,\bar{\mathbf{u}}_3\}\cap\mathcal{C}_{\frac{1}{100}}(E^s_{4}(\gamma_I(0)))=\{0\}$. 
  \begin{figure}[htb]
\begin{subfigure}{0.48\textwidth}
\includegraphics[height=4.5cm,width=7cm]{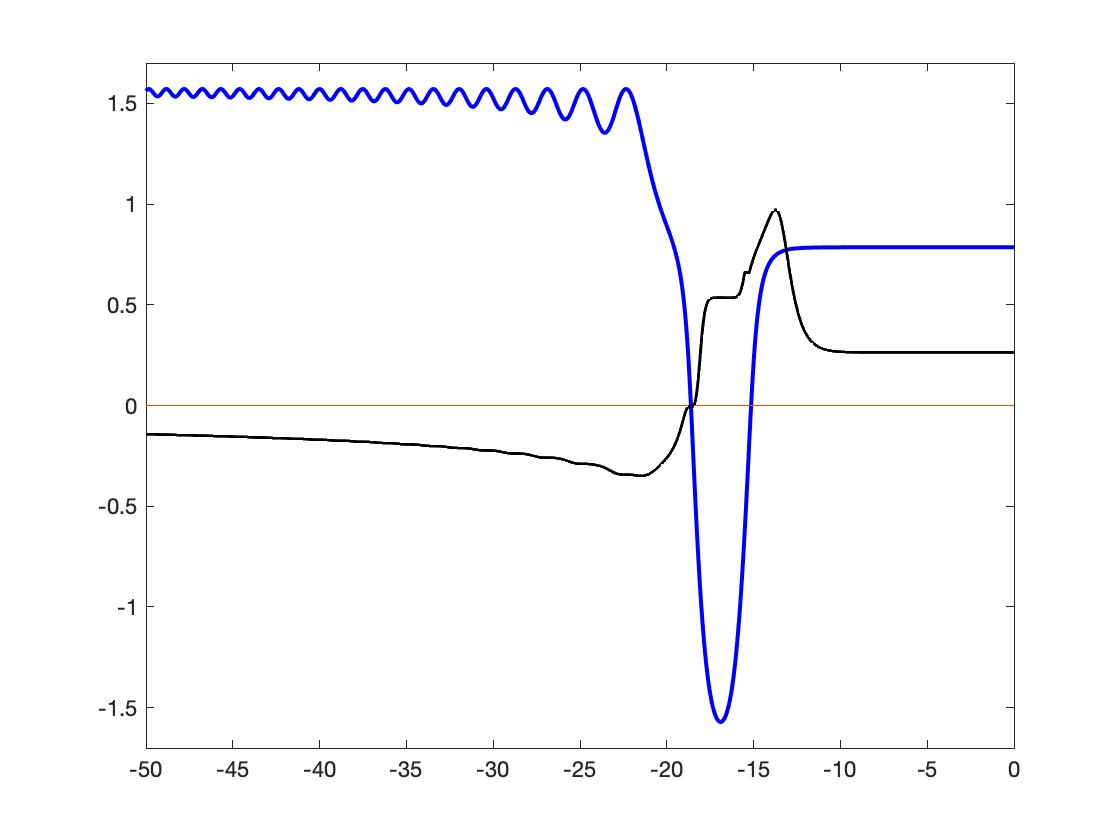}
\end{subfigure}
\begin{subfigure}{0.48\textwidth}
\includegraphics[height=4.5cm,width=7cm]{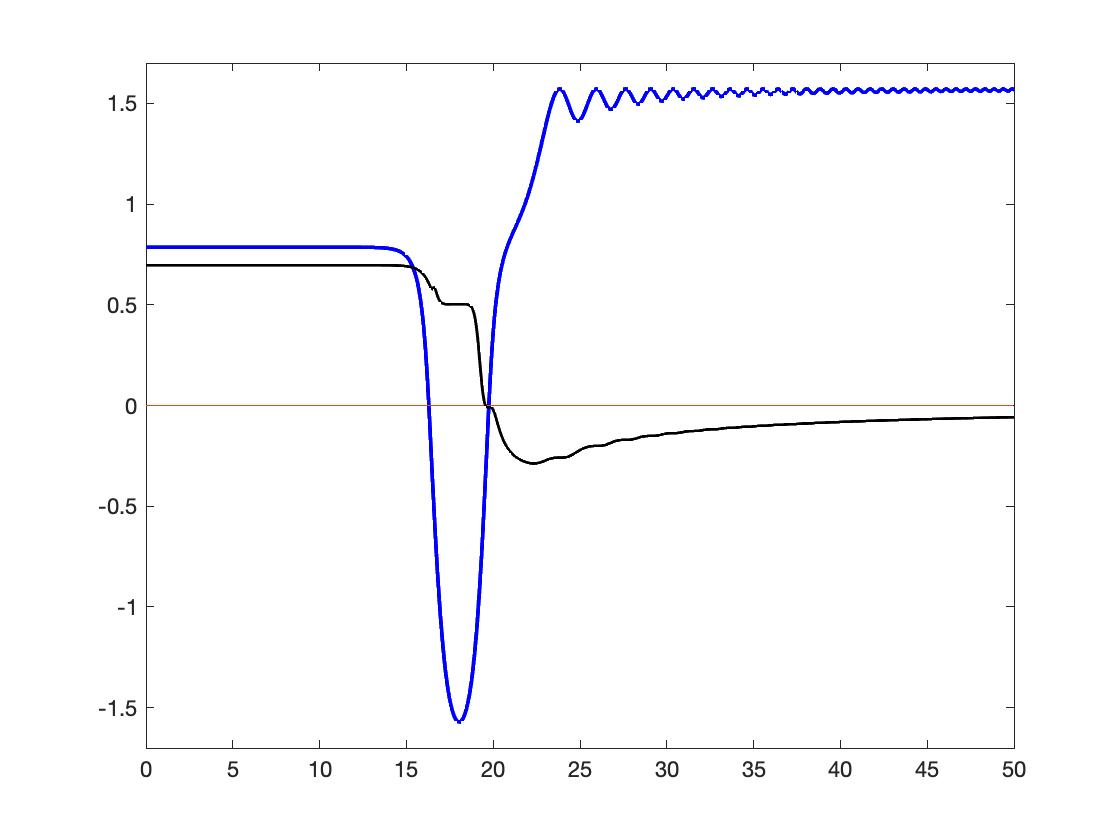}
\end{subfigure}
\caption{The left: backward solution along $\gamma_I$ and  the right: forward solution along $\gamma_O$; the bright curves are the evolution of the $\psi$-variable and the dark curves are that of the $dw_0$-variable. }
\label{fig-dw0}
\end{figure}

  {\bf Step 4:} {\it Verifying Proposition \ref{PropTrans1}. } We just keep track of $dw_0$ components of $E^u_0(\gamma_O(t))$, $t<0$ and $E^s_0(\gamma_I(t))$, $t>0$.  The evolution of $dw_0$ and $\psi$ is illustrated in Figure \ref{fig-dw0}. The values are close to nonzero constants when the orbits are in a neighborhood of the Lagrange fixed point as well as close to the sections $\cS_\pm$, so they can only become zero at finitely many points. On the other hand, along $\gamma_I$ and $\gamma_O$, the variable $\psi$ also becomes zero at finitely many points, in particular it does not vanish in a neighborhood of the Lagrange fixed point. This reduces the problem to checking the non-coincidence of the zeros of $\psi\pm\pi/2$ and $w_0$ on two compact time intervals. Hence we have numerically justified the non-degeneracy  in Proposition \ref{PropTrans1}.

%\section*{Acknowledgement}
%The authors are supported by the grant NSFC (Significant project No.11790273) in China. 

\end{document}